\documentclass[12pt]{memo-l}
\usepackage{hyperref}

\usepackage{tikz}

\usepackage{enumerate}

\usepackage{amssymb, mathrsfs,amsmath}
\usepackage{latexsym,array}
\usepackage{amsfonts}
\usepackage{shadow}

\def\e{\epsilon_N}
\def\C{\mathbb C}
\def\R{\mathbb R}
\def\D{\mathbb D}
\def\N{\mathbb N}
\def\H{\mathbb H}

\date{}

\newcommand{\bc}{\mathbb B\mathbb C}
\renewcommand{\k}{{\bf k}}
\renewcommand{\i}{{\bf i}}
\renewcommand{\j}{{\bf j}}
\renewcommand{\e}{{\bf e}}
\newcommand{\la}{\lambda}
\newcommand{\edag}{{\bf e^\dagger}}

\newtheorem{theorem}{Theorem}[chapter]

\theoremstyle{definition}

\theoremstyle{remark}
\newtheorem{remark}[theorem]{Remark}

\newtheorem{subtheorem}{Theorem}[section]

\newtheorem{subPn}[subtheorem]{Proposition}
\newtheorem{subCy}[subtheorem]{Corollary}

\newtheorem{subDn}[subtheorem]{Definition}
\newtheorem{subEx}[subtheorem]{Example}
\newtheorem{subremark}[subtheorem]{Remark}

\newtheorem{subParr}[subtheorem]{}

\numberwithin{section}{chapter}
\numberwithin{equation}{chapter}

\makeindex
\begin{document}
\frontmatter

\title
{
Basics of functional  analysis  with   bicomplex  scalars, and
bicomplex      Schur analysis }

\author{Daniel Alpay}
\address{(DA) Department of mathematics,
Ben-Gurion University of the Negev, P.O. Box 653, Beer-Sheva
84105, Israel} \email{dany@math.bgu.ac.il}

\author
{Maria Elena Luna-Elizarrar\'as}
\address{(ML-E) Escuela Superior
de Fisica y Matem\'aticas, Instituto Polit\'ecnico Nacional,
Mexico City, M\'exico.} \email{eluna@esfm.ipn.mx}

\author{Michael Shapiro}
\address{(MS) Escuela Superior
de Fisica y Matem\'aticas, Instituto Polit\'ecnico Nacional,
Mexico City, M\'exico.} \email{shapiro@esfm.ipn.mx}

\author{Daniele C. Struppa}
\address{(DS), Schmid College of Science and Technology, Chapman University, Orange,
CA 92866, USA}
\email{struppa@chapman.edu}

\pagestyle{plain}
\pagenumbering{arabic}

\subjclass[2000]{Primary }

\keywords{bicomplex numbers}

\maketitle

\tableofcontents

\abstract
oooooojjjjjj
\endabstract

\def\i{\mathbf i}
\def\j{\mathbf j}
\def\k{\mathbf k}

\chapter*{Introduction}

Bicomplex numbers have been studied for quite a long time, probably beginning with the work of the Italian school of Segre \cite{segre},  Spampinato  \cite{S1,S2},   and Scorza Dragoni \cite{scorzadragoni}. Their interest arose   from the fact that such numbers   offer a commutative alternative to the skew field of quaternions       (both sets  are  real  four--dimensional  spaces),     and that in many ways they generalize  complex numbers more closely and more accurately than quaternions do.    Of course, commutativity is gained at a price, and in this case the price is the fact that the ring of bicomplex numbers is not a field, since zero divisors arise to prevent such a possibility.    Anyway,  one  may  expect  that  bicomplex  numbers  can serve  as  scalars  both  in  the  theory  of  functions   and  in  functional  analysis,  being  at  least  a  reasonable  counterpart  for  the  quaternions.   \\

The most comprehensive study of analysis in the bicomplex setting    is certainly the book of G. B. Price \cite{Price}, and in recent years there has been a significant impulse to the study of the properties of those functions on the ring  $\bc$ of bicomplex numbers, whose properties suggest a notion of bicomplex holomorphy. Rather than giving an exhaustive list of references, we refer the reader to the  article \cite{LSSV} and to the forthcoming monograph \cite{LSSVbook}. As demonstrated in these references (as well as in the more specialized  papers on which those references rely), the fundamental aspects of   bicomplex analysis (the analysis of bicomplex holomorphic functions) are by now fairly well understood, and it is possible to study  some of the more delicate aspects of the theory of the modules   of such bicomplex holomorphic functions. The history of complex analysis indicates that progress in such an arena (the study of analytic functionals, for example) cannot be achieved without a strong grasp of functional analysis. Thus, the origin of this monograph. \\

With the goal of providing the foundations for a rigorous study of modules of bicomplex holomorphic functions, we  develop here a general theory of functional analysis with   bicomplex scalars.    Functional analysis in $\bc$ is an  essentially      new subject,      
and  it  seems  to  have  had     two  
independent starting points: the paper \cite{RT}, submitted and published in 2006, and the earlier paper \cite{LS2009-2}, originally submitted back in 2005.   It   may  be  instructive  to compare this with the case of functional analysis with the quaternionic  scalars  which   dates back    to the  papers of Teichmuller \cite{Te} and Soukhomlinoff \cite{So},  
and  which   is  actually  widely  known and studied;  for  some  recent  works see, e.g.,    \cite{AS2004, ALS2005, ALS2007, LS2005_1,LS2005, LS2009}.   \\

During the last several years, the initial ideas    of bicomplex functional analysis have been studied, see,  e.g.,    
\cite{LMR2010, LMR2011-1, LMR2011-2}   and  also   \cite{KKR}.    While  there  are   some  overlaps  with  our  work,   
 we present here a much more complete and exhaustive treatment of the theory, as well as many new ideas and results. In particular, we show how the general ideas that we develop for a bicomplex functional analysis can be directly employed to generalize the classical Schur analysis to the bicomplex setting. \\

Let us now describe in more detail the contents of each chapter in this monograph.  \\

Even though the basic properties of bicomplex number are   well known and widely available, our analysis requires some more delicate discussion of the various structures which are hidden in the ring of bicomplex numbers. For this reason we use Chapter 1 to study in detail the subset $\D$ of hyperbolic numbers, and we are able to establish why we claim that $\D$ plays, for bicomplex numbers, the same role that $\R$ plays for complex numbers. In particular, we introduce a new  partial order on $\D$, which has interesting connections with the Minkowksi space of special relativity, and that allows us to introduce a new hyperbolic-valued norm on the set of bicomplex numbers.  \\

In Chapter 2, these ideas are extended to the case of matrices with bicomplex entries. Even though no surprises occur here, the study of bicomplex matrices cannot be found elsewhere, and we     present here  a  detailed  description of the peculiarities  which arise here.  Besides,  we  include into this chapter a quick survey of  holomorphicity  in the  bicomplex   setting.   \\

In Chapter 3, we study $\bc$--modules, and we show their subtle nature as well as the many structures that one can impose on them. In particular we will show that every $\bc$--module decomposes into the direct sum of submodules,   which   we  call  the   idempotent decomposition of a $\bc$--module (a notion inherited from the idempotent decomposition for bicomplex numbers): such a decomposition will play a key role in the remainder   of the work. We also show how to use two complex linear spaces to generate a $\bc$--module of which they are the idempotent representation.  \\

The most surprising portion of the work takes place in Chapter 4, where we study inner products   and  norms   in bicomplex modules.    The matter is that  depending on the kind of the scalars an inner  product can take real,  complex,  or quaternionic values  (they can be even more general) but  the corresponding norm is always real--valued.  We  consider two kinds of norms on bicomplex modules: a real-valued norm (as one would expect), and a hyperbolic-valued norm. Interestingly enough, while both norms can be used to build the theory of normed bicomplex modules,    the hyperbolic-valued norm appears to be  much    better compatible with the structure of $\bc$--modules.   Note  also  that  we  apply  the developed    told  for  a  study of the ring   $\mathbb H (\C)$ of biquaternions  (complex quaternions)  seen as a  $\bc$--module.     Since  the  $\H (\C)$--valued  functions arise in a  wide  range   of    areas  we  believe  that   our  analysis  will  be  rather  helpful  for    constructing  the theory of linear space of such functions.   \\

Chapter 5 sets the stage for the study of linear functionals on bicomplex modules. The results in this chapter are sometimes surprising and will set the stage for any advanced function theory for bicomplex function spaces. Recent works have studied spaces of analytic functionals on spaces of bicomplex holomorphic spaces \cite{struppa, struppavv}, but in those articles the duality is defined by considering only the $\C$--linear structure. We expect the results from this chapter to allow a new understanding of those spaces of analytic functionals and a deeper description of the duality properties. \\

Finally, in Chapter 6, we describe a bicomplex version of the classical Schur analysis. This is   a    significant application of the theory developed in the first five chapters of this monograph. Schur analysis originates with the classical papers of Schur, Herglotz, and other   (see \cite{MR2002b:47144}
for a survey)     and studies holomorphic contractive functions in the open unit disc of $\C$. This analysis has remarkable connections and applicationt to interpolation problems, moment problems, and the theory of linear systems. It is because of these relations that Schur analysis has already been extended to many other settings, for example several complex variables \cite{agler-hellinger},    slice-holomorphic functions \cite{acs1},    and hyperholomorphic functions,   \cite{asv-cras,MR2124899,MR2240272,MR2275397}.    In this final chapter of our work, we show how these same ideas can be extended in a very successful way to the case of bicomplex spaces, and we consider, in this setting, Blaschke factors, the Hardy spaces, and the notion of Schur multipliers and their realizations.  \\

The Mexican  authors   have been partially supported by
  CONACYT projects as well by Instituto Polit\'ecnico Nacional in the
  framework of COFAA and SIP programs;   they are  also   
  grateful to Chapman University for the support offered in preparing
   this work.

\chapter{Bicomplex  and  hyperbolic numbers }

\section{Bicomplex numbers}
\label{intro bicomplex numbers}
The ring  of  bicomplex  numbers  is  the  commutative  ring
$\mathbb B \mathbb C$   defined as  fo\-llows:
$$    \displaystyle   \mathbb B \mathbb C  :=  \left\{   Z =  z_1  + z_2 \,
\mathbf j  \mid  z_1 , \, z_2   \in  \mathbb C ( \mathbf i ) \,  \right\}   .$$
where  $\mathbf i$  and  $\mathbf j $  are  commuting  imaginary
units, i.e.,
$$\mathbf i \, \mathbf j = \mathbf j \, \mathbf i,    \quad \mathbf i^2  = \mathbf j^2 = -1 ,   $$
 and $\mathbb C ( \mathbf i) $ is
the set  of  complex  numbers  with the imaginary  unit $\mathbf
i$.  Observe that  for the particular case  where      $z_1=x_1$ and
$z_2= x_2$ are  real numbers,  $Z=  x_1 + x_2 \, \mathbf j$ is a
complex number with the imaginary unit $\mathbf j$. Since the two
imaginary units $\mathbf i$ and $ \mathbf j$ coexist  inside
$\mathbb B \mathbb C$, in what follows  we will  distinguish  between   the
two sets of complex numbers  $ \mathbb C ( \mathbf i )$ and $
\mathbb C ( \mathbf j )$.  Note also that   since $(\mathbf i \mathbf j)^2=1,$
if $ z_1 = x_1 $  is real    and  $z_2 =  y_2 \,
\mathbf i$ is a purely  imaginary number,  one has that     $Z=  x_1 +  y_2 \,
\mathbf i \, \mathbf j$ is an element of the set $\mathbb D$
of hyperbolic numbers, defined  as
$$\mathbb D :   = \left\{ \, a + b \, \mathbf k \mid a,b \in \mathbb R , \; \mathbf k^2 =1 , \; \mathbf k \notin  \mathbb R \, \right\} . $$

Thus the    subset     $\left\{ \, x_1 + y_2 \, \mathbf i \, \mathbf j \mid x_1 , y_2 \in \mathbb R \, \right\}$   in   $\mathbb B \mathbb C$  is isomorphic    (as real algebra) to $\mathbb D$  and can be identified  with  it.  \\

We  will soon   see   that  in addition  to   containing  two copies  of the complex numbers, the whole set $\bc$  has many  deep  similarities with
the set  of complex numbers although, of course,   many  differences, sometimes quite striking,  arise also.  \\


\section{Conjugations and moduli.} \label{section conj and mod}

Any  bicomplex  number can be written in six different
ways, which are relevant for what follows.
\begin{equation}\label{different ways}
\begin{array}{rcl}
Z &   =  &   (x_1 + \mathbf i \, y_1 )  +  ( x_2  + \mathbf i \, y_2 ) \, \mathbf j  \; =:  \;  z_1 + z_2 \, \mathbf j
\\  &  &  \\  & =   &   (x_1 + x_2 \, \mathbf j )  + ( y_1 + y_2 \, \mathbf j ) \, \mathbf i  \; =  :   \;  \eta_1 + \eta_2 \, \mathbf i \,
\\   &   &  \\   & = &
(x_1  +  \i \j  \, y_2 )  + \j  (x_2  - \i \j  \, y_1)  =:  \mathfrak z_1  + \j \,  \mathfrak z_2
\\  &  &  \\  & = &
(x_1  + \i \j \, y_2 ) + \i ( y_1 - \i \j \, x_2)  =:  \mathfrak x_1  + \i    \, \mathfrak x_2  \,
\\  &  &  \\  & = &
(x_1  + \i  \, y_1 ) + \i \j  ( y_2 -    \i  \, x_2)  =:  \alpha_1  + \k      \, \alpha_2  \,
\\  &  &  \\  & = &
(x_1  +  \j \, x_2 ) + \i  \j   ( y_2 -     \j \, y_1 )  =:  \nu_1  + \k    \, \nu_2  \,  ,
\end{array}
\end{equation}

where $z_1 , \, z_2 , \, \alpha_1 , \, \alpha_2     \in \mathbb C ( \mathbf i )  $,  \,   $ \eta_1  \, \eta_2  , \nu_1 , \, \nu_2    \in \mathbb C  ( \mathbf j) $ \,    and $\mathfrak z_1$ ,  $  \mathfrak z_2$, \,$ \mathfrak x_1 , \,  \mathfrak x_2   \in \mathbb D$.  We present  all these representations of bicomplex numbers not only for completeness  but also because  all of  them   manifest themselves  in the study  of  the structure  of the modules with bicomplex scalars. We will later introduce one more representation, known as the idempotent  representation of a bicomplex number.  \\

Since $\bc$ contains two imaginary units whose square is $-1$,  and one hyperbolic unit whose square is 1, we can consider three   conjugations   for  bicomplex  numbers in analogy  with  the  usual  complex  conjugation:  \\

\begin{enumerate}

\item[({\sc i})]  $  \overline Z  :=    \overline z_1   +
\overline z_2  \, \mathbf j   \, $  \;    (the  bar-conjugation);

\medskip

\item[({\sc ii})]  $Z^\dagger   :=    z_1  - z_2  \, \mathbf j  $ \; \,   (the  $\dagger$--conjugation);

\medskip

\item[({\sc iii})] $ Z^\ast  :=   \left( \,  \overline Z  \,
\right)^\dagger  = \overline{  \left(  Z^\dagger \right) }  =
\overline z_1  - \overline z_2 \,  \mathbf j   \, $  \;  (the  $\ast$--conjugation), \\

\end{enumerate}

where  $ \overline{ z}_1 , \overline{z}_2 $    denote   the   usual   complex  conjugates  to   $z_1, z_2$ in   $\C (\i)$.   \\

Let us see  how these  conjugations  act  on the complex  numbers  in    $\C(\i)$ and   in   $\C(\, \j)$ and  on  the   hyperbolic numbers in $\D$.  If $Z= z_1  \in   \C (\i) $,  i.e.,  $z_2=0$,  then $Z=  z_1 =  x_1 + \i y_1$  and  one has:
$$  \overline{Z} = \overline{z}_1 =   x_1 -  \i y_1   =  z_1^\ast = Z^\ast  , \qquad \quad   Z^\dagger = z_1^\dagger = z_1  ,$$

that is,  both the  bar--conjugation  and  the   $\ast$--conjugation,  restricted  to   $\C(\i)$,  coincide  with  the  usual  complex  conjugation  there,  and  they both fix all  elements of  $\C(\i)$.  \\

If    $Z=  \eta_1$  belongs  to  $\C(\, \j)$,  that   is   $\eta_1     =  x_1  + x_2 \, \mathbf j $,  then      one  has:
$$    \overline{\eta_1}  = \eta_1 ,   \qquad  \quad  \eta_1^\ast =  x_1  - x_2 \, \mathbf j  =  \eta_1^\dagger   ,  $$

that  is,  both the   $\ast$--conjugation  and  the  $\dagger$--conjugation,      restricted  to   $  \mathbb C  ( \mathbf j ) $,    coincide  with  the  usual conjugation  there.  In  order  to  avoid   any   confusion  with the notation, from now on we will identify  the  conjugation  on  $  \mathbb C  ( \mathbf j ) $  with the   $\ast$--conjugation.   Note  also that  any  element  in   $\C(\, \j)$   is  fixed by the  bar--conjugation.  \\

Finally, if  $Z  =  x_1  + \i \j \, y_2 \in \mathbb D$,  that  is  $  y_1= x_2= 0 $, then
$$ \overline Z = x_1 - \i \j \, y_2 = Z^\dagger   ,  \qquad  \quad    Z^\ast =Z  ,    $$

Thus,  the  bar--conjugation  and  the  $\dagger$--conjugation   restricted  to  $\D$  coincide with the  intrinsic  conjugation  there.
We will  use   the  bar--conjugation  to  denote  the  latter.     Note  that  any  hyperbolic number is fixed b the $\ast$--conjugation.     \\











Each  conjugation  is an additive, involutive, and multiplicative operation on $\bc$:

\medskip

\begin{enumerate}

\item[({\sc i})]   $ \overline{  \left( Z + W  \right) }   =
\overline Z  +  \overline W$;    \;   $\left(  Z + W
\right)^\dagger  =   Z^\dagger  + W^\dagger$;    \\  \\   $\left(  Z +
W \right)^\ast  =   Z^\ast  + W^\ast$.

\bigskip

\item[({\sc ii})]  $\overline{ \overline Z }  = Z $;    \;  $
\left( Z^\dagger \right)^\dagger  = Z$;   \;  $\left( Z^\ast
\right) ^\ast =Z$.

\bigskip

\item[({\sc iii})]    $ \overline{  \left(  Z \cdot W \right) }
=  \overline Z  \cdot \overline W$;    \;   $ \left( Z \cdot W
\right)^\dagger  = Z^\dagger \cdot W^\dagger$;   \; $ \left( Z
\cdot W \right)^\ast  = Z^\ast \cdot W^\ast  $.

\end{enumerate}

Thus, each  conjugation  is  a
ring  automorphism  of  $\mathbb B \mathbb C$.  \\

In the complex case the modulus  of a complex number is   intimately  related  with     the complex conjugation:  multiplying  a  complex  number  by its conjugate one gets  the  square of its  modulus.    Applying this idea to each of the     three  conjugations, three  possible  ``moduli"  arise  in  accordance  with the formulas  for  their  squares:

\medskip

\begin{itemize}
\item  $  | Z |_{\mathbf i}^2  :=     Z \cdot   Z^\dagger   =
z_1^2  + z_2^2   $
\\   \\
\hbox{}   $     \qquad    =   \left(  \,   | \eta_1  |^2  -   | \eta_2 |^2  \right)    +   2 \, Re  \,
(  \eta_1  \, \eta_2^\ast \, ) \, \mathbf  i $
\\   \\
\hbox{}   $     \qquad  =   \displaystyle  \left(   |  \mathfrak z_1 |^2  + | \mathfrak z_2 |^2  \right)  + \j \,  \left(      \overline{ \mathfrak z}_1 \, \mathfrak z_2  -  \overline{ \, \overline{ \mathfrak z}_1 \, \mathfrak z_2   }  \, \right)  $
\\   \\
\hbox{}    $ \qquad      =   \displaystyle  \left(   |  \mathfrak x_1 |^2    -   | \mathfrak x_2 |^2  \right)  + \i \,  \left(  \mathfrak x_1 \,     \overline{ \mathfrak x}_2 +    \overline{ \,   \mathfrak x_1 \,  \overline{ \mathfrak x}_2 \, }  \, \right)   $
\\   \\
\hbox{}    $ \qquad      =   \displaystyle  \alpha_1^2  - \alpha_2^2 $
\\   \\
\hbox{}    $ \qquad      =   \displaystyle  \left(   |  \nu_1 |^2    -   | \nu_2 |^2  \right)  -   \i \,  2 \, Im \, \left(  \nu_1^\ast  \nu_2   \, \right)
  \in \mathbb C ( \mathbf i ) $;

\bigskip

\item  $  | Z |_{\mathbf j}^2  :=  Z \cdot  \overline Z      =
\left(  \,   | z_1  |^2  -   | z_2 |^2  \right)    +   2 \, Re  \,
(  z_1  \, \overline z_2 \, ) \, \mathbf  j  $
\\   \\
\hbox{}    \qquad   $ =   \;  \eta_1^2  + \eta_2^2  $
\\  \\
\hbox{}  \qquad   $  \displaystyle  =   \;   \left(   |  \mathfrak z_1 |^2  -  | \mathfrak z_2 |^2  \right)  + \j \,  \left(    \mathfrak z_1 \,    \overline{ \mathfrak z}_2 +  \overline{ \, \mathfrak z_1  \,   \overline{ \mathfrak z}_2 \,    }  \, \right)    $
\\   \\
\hbox{}    \qquad   $   =  \;    \displaystyle  \left(   |  \mathfrak x_1 |^2    +     | \mathfrak x_2 |^2  \right)  + \i \,  \left(  \overline{ \mathfrak x}_1  \,  \mathfrak x_2 \,    -     \overline{ \,   \overline{ \mathfrak x}_1 \,  \mathfrak x_2  }  \, \right) $
\\   \\
\hbox{}    \qquad   $   =  \;    \displaystyle  \left(   |  \alpha_1 |^2   -     | \alpha_2 |^2  \right)  + \k \,  \left(    \alpha_2 \,   \overline{ \alpha}_1     -  \alpha_1 \,      \overline{ \alpha}_2 \,   \right) $
\\   \\
\hbox{}    \qquad   $   =  \;    \displaystyle    \nu_1^2   -      \nu_2^2       \in \mathbb C (
\mathbf j ) $;

\bigskip

\item     $  | Z |_{\mathbf k}^2  :=  Z \cdot   Z^\ast     =        \left(
\,   |  z_1  |^2  +   | z_2 |^2  \right)    -   2 \, Im  \,  (
z_1  \, \overline z_2 \, ) \, \mathbf  k  $
\\     \\
\hbox{}    \qquad   $  =  \;
 \left(
\,   |  \eta_1  |^2  +   | \eta_2 |^2  \right)    -   2 \, Im  \,  (
\eta_1  \, \eta_2^\ast \, ) \, \mathbf  k     $
\\   \\
\hbox{}    \qquad   $  = \;    \mathfrak z_1^2  + \mathfrak z_2^2   $
\\   \\
\hbox{}    \qquad   $
=  \mathfrak x_1^2   + \mathfrak x_2^2   $
\\     \\
\hbox{}    \qquad   $  =  \;
 \left(
\,   |  \alpha_1  |^2  +   | \alpha_2 |^2  \right)    +  \k  \,  \left(    \alpha_2  \, \overline{ \alpha}_1   +      \alpha_1    \,   \overline{ \alpha}_2 \, \right) \,      $
\\     \\
\hbox{}    \qquad   $  =  \;
 \left(
\,   |  \nu_1  |^2  +   | \nu_2 |^2  \right)    +   \k   \,  \left(  \nu_1  \, \nu_2^\ast      + \nu_2  \, \nu_1^\ast   \, \right)
\;   \in \mathbb D ,   $

\end{itemize}

\medskip

where  for  a  complex  number $z$ (in $ \mathbb C (\mathbf i )$  or  $ \mathbb C (\mathbf j )$) we denote  by   $| z |$    its usual modulus and for a hyperbolic  number  $\mathfrak z  =  a + b \,  \k$ we use  the notation  $| \mathfrak z |^2  = a^2 - b^2 $.  \\

Unlike what happens in the complex case, these moduli are not $\mathbb R^+$-valued. The first two moduli  are complex-valued (in $\mathbb C( \mathbf i )$ and    $\mathbb C( \mathbf j )$   respectively),  while the last one is hyperbolic-valued.  These moduli nevertheless behave as expected with respect to multiplication. Specifically, we have

$$  |  Z \cdot W |_{\mathbf i}^2  =  | Z |_{\mathbf i}^2  \cdot
| W   |_{\mathbf i}^2 \, ;   $$

$$  |  Z \cdot W |_{\mathbf j}^2  =  | Z |_{\mathbf j}^2  \cdot
| W   |_{\mathbf j}^2 \, ;   $$

$$  |  Z \cdot W |_{\mathbf k}^2  =  | Z |_{\mathbf k}^2  \cdot
| W   |_{\mathbf k}^2 \, .  $$

\begin{remark}
The hyperbolic-valued modulus $| Z |_\k$ of a bicomplex number $Z$ satisfies
$$    | Z |_{\mathbf k}^2   =     \left(   | z_1 | ^2  +  |  z_2 |^2  \right)   +
\left(  -2 \, Im \, ( z_1 \, \overline z_2 ) \, \right) \,
\mathbf k  =:  a  + b \, \mathbf k ,  $$

where  $a$ and  $b$  satisfy  the  inequalities
$$  a^2  - b^2  \geq 0 \,   \qquad  {\rm and} \quad  a\geq 0  $$
(this is a consequence of the  fact  that  $ \displaystyle
| \, Im \,  ( z_1 \, \overline z_2 \, ) \, |  \leq | z_1 | \cdot
| z_2 | $).

\end{remark}

This remarks justifies the introduction of  the set of  ``positive"  hyperbolic  numbers:
$$   \displaystyle   \mathbb D^+  :=  \left\{  a  +  b \,
\mathbf k  \mid   a^2   - b^2  \geq 0  ,  \;  a \geq 0 \,
\right\} \, , $$

so that   $    | Z |_{\mathbf k}^2  \in  \mathbb D^+$.     Such  a  definition  of  ``positiveness"  for hyperbolic  numbers  does not  look  intuitively  clear  but  later  we  give  another  description of  $\D^+$  that clarifies    the  reason  for  such  a  name.   It  turns out  that  the  positive  hyperbolic numbers  play  with respect to  all  hyperbolic numbers   a   role  deeply similar  to that  of  real  non negative   numbers  with  respect  to    all   real  numbers.   \\

\section{The Euclidean norm  on   $\bc$}\label{subsection 2.3}

  Since none of the moduli above is real valued, we can consider  also the  Euclidean norm on  $\mathbb B \mathbb C$ seen  as  $\mathbb C^2( \mathbf i )  = \left\{ \, (z_1 , z_2) \mid  z_1 + z_2 \, \mathbf j \in  \mathbb B \mathbb C \, \right\} $,  as   $\mathbb C^2( \mathbf j )  = \left\{ \, (\eta_1 , \eta_2) \mid  \eta_1 + \eta_2 \, \mathbf i \in  \mathbb B \mathbb C \, \right\} $,  or  as  $\mathbb R^4  =  \left\{  \,    (x_1 , y_1 , x_2 , y_2) \mid  x_1 +  \right.  $     $\left. + \,   \i y_1 +    \j x_2 + \k y_2 \in \bc \,   \right\}$.    This  Euclidean norm  is connected to  the  properties
of    bicomplex  numbers  via  the  $\mathbb D^+$-valued  modulus as follows:
$$ \begin{array}{rcl}
 | Z |  &   =   &  \sqrt{x_1^2 + y_1^2 + x_2^2 + y_2^2 }   \,  =   \,       \sqrt{ | z_1 | ^2  +  |  z_2 |^2  \, }
 \\  &  &  \\   &   =  &      \sqrt{ | \eta_1 | ^2  +  |  \eta_2 |^2  \, }  \,    =   \,     \sqrt{ Re \,
\left(  | Z |_k^2   \, \right) \, }  \, .
\end{array}$$

It is easy to  prove  (using  the  triangle inequality) that, for  any   $Z$  and  $U$  in  $\bc$,
\begin{equation}\label{desigualdad norma euclideana}
| Z \cdot U |  \leq  \sqrt{2} \, |Z | \cdot  | U | \, .
\end{equation}

We can actually show that  if $U  \in  \mathbb B \mathbb C$ is arbitrary but $Z$ is either a complex   number    in  $\mathbb C ( \mathbf i )$     or   $\mathbb C ( \mathbf j )$,   or a  hyperbolic  number  then we can say something more:

\medskip

\begin{enumerate}

\item[ {\it a})]  if $Z \in \mathbb C ( \mathbf i )$     or   $\mathbb C ( \mathbf j )$   then   $ | Z \cdot U |  =  | Z | \cdot | U | $;

\medskip

\item[ {\it b})]  if $Z \in \mathbb D$  then

$$ | Z \cdot U |^2  =  | Z |^2 \cdot | U |^2   +   4 \, x_1 \, y_2 \, Re ( \mathbf i \, u_1 \, \overline u_2 ) $$

\noindent
where $Z =   x_1 +    \k y_2  $   and   $U  =  u_1 + \j  u_2$.

\end{enumerate}

We will  prove these  two    properties  in the next section.

\bigskip

\section{Idempotent  decompositions}

Since for any  bicomplex  number  $Z  =  z_1 + z_2 \j$ it is
$$ Z \cdot Z^\dagger = z_1^2 + z_2^2 \in \mathbb C ( \mathbf i ) ,  $$

it follows  that any   bicomplex  number  $Z$  with    $  | Z |_{\mathbf i} \neq 0$  is invertible,
and its inverse is given by
$$ \displaystyle  Z^{ -1 }   =  \frac{ Z^\dagger}{ | Z |_{\mathbf i}^2 } \,$$

If, on the other hand,  $Z \neq 0$  but  $   | Z |_{\mathbf i} =0$  then $Z$ is   a zero   divisor. In fact, there are no other   zero  divisors.\\

We  denote the set of zero divisors  by $\mathfrak S$, thus
$$\mathfrak S : =  \left\{  \, Z \mid    Z \not=0,  \;   z_1^2 + z_2^2 =0 \, \right\} .$$

It  turns out that there  are  two very special zero  divisors. Set
$$  \displaystyle   \mathbf e :=  \frac{1}{2}  \, \left( 1  +
\mathbf i \, \mathbf j \, \right)  ,   $$

then  its  $\dagger$--conjugate is
$$  \displaystyle      \mathbf e^\dagger
=   \frac{1}{2}  \, \left( 1  -   \mathbf i \, \mathbf j \,
\right)  \, .$$

It  is   immediate   to  check  that
$$  \mathbf e^2  = \mathbf e \, ;  \qquad   \left( \mathbf e^\dagger \,
\right)^2  =  \mathbf e^\dagger \, ; \qquad   \mathbf e  +  \mathbf e^\dagger  = 1\, ; $$

$$  \mathbf e^\ast = \mathbf e \, ,  \quad  \left( \mathbf e^\dagger \,
\right)^\ast  = \mathbf e^\dagger \, ; \qquad    \mathbf e \cdot   \mathbf e^\dagger =0 \, .$$

The  last  property  says  that  $\mathbf e $ and     $\mathbf e^\dagger $  are  zero  divisors,  and  the  first  three  mean  that  they are mutually  complementary idempotent  elements.   \\

Thus,  the  two  sets
$$  \mathbb B \mathbb C_{ \mathbf e} :=     \mathbf e \cdot      \mathbb B \mathbb C  \quad  {\rm and}  \qquad     \mathbb B \mathbb C_{ \mathbf e^\dagger} :=     \mathbf e^\dagger \cdot      \mathbb B \mathbb C  $$

are     (principal)   ideals in the  ring     $ \mathbb B \mathbb C$  and  they  have  the  properties:
$$    \mathbb B \mathbb C_{ \mathbf e}   \cap        \mathbb B \mathbb C_{ \mathbf e^\dagger}   = \{  0 \}  $$

and
\begin{equation}\label{descomposicion idempotente bicomplejos}
\mathbb B \mathbb C  =     \mathbb B \mathbb C_{ \mathbf e}  +       \mathbb B \mathbb C_{ \mathbf e^\dagger}  \, .
\end{equation}

\medskip

We  shall call  (\ref{descomposicion idempotente bicomplejos})  the  idempotent  decomposition  of   $ \mathbb B \mathbb C $,  and  we  shall  see  later  that   bicomplex  modules  inherit   from  their  scalars  a  similar  decomposition.  Of  course,  both  ideals   $  \mathbb B \mathbb C_{ \mathbf e}  $  and  $    \mathbb B \mathbb C_{ \mathbf e^\dagger}  $  are  uniquely  determined  but  their  elements  admit different representations. In fact,   every  bicomplex  number  $Z =  ( x_1 +  \mathbf i \, y_1 )  +      ( x_2 +  \mathbf i \, y_2 )  =  z_1  + z_2 \, \mathbf j $  can  be  written  as
\begin{equation}\label{repr_idempotent}
Z = z_1  + z_2 \, \mathbf j  = \beta_1  \, \mathbf e  + \beta_2 \,
\mathbf e^\dagger \, ,
\end{equation}

where
\begin{equation}
\label{beta12} \beta_1  =  z_1  -  \mathbf i \, z_2\quad{\rm
and}\quad\beta_2 =  z_1   + \mathbf i \, z_2 ,
\end{equation}

are complex  numbers  in  $\C(\i)$. On the other hand, $Z$ can also be written as
\begin{equation}\label{repr_idempotent_Cj}
Z = \eta_1  + \eta_2 \, \mathbf i  = \gamma_1  \, \mathbf e  + \gamma_2 \,
\mathbf e^\dagger \, ,
\end{equation}

where   $\eta_1 :=  x_1 + x_2 \, \mathbf j $,    $\eta_2 :=  y_1 + y_2 \, \mathbf j $,     $ \gamma_1 : =  \eta_1  -  \mathbf j \, \eta_2 $,  $ \gamma_2
=  \eta_1   + \mathbf j \, \eta_2 $   are  complex  numbers  in   $\C( \, \j)$.    Each  of  the  formulas    (\ref{repr_idempotent})   and    (\ref{repr_idempotent_Cj})    can  be  equally  called  the  idempotent  representation  of  a  bicomplex  number. More  specifically,    (\ref{repr_idempotent})  is  the  idempotent  representation for   $\mathbb B  \mathbb C$  seen  as  $\mathbb C^2 ( \mathbf i ) :=  \mathbb C ( \mathbf i ) \times   \mathbb C ( \mathbf i )$  and   (\ref{repr_idempotent_Cj})  is  the  idempotent  representation  for   $\mathbb B  \mathbb C$   seen  as
$\mathbb C^2 ( \mathbf j ) :=  \mathbb C ( \mathbf j ) \times   \mathbb C ( \mathbf j )$.  Usually,  (see   for  instance   \cite{RS},     \cite{RT},    \cite{Price})  only  the  representation   (\ref{repr_idempotent})  is   considered,  but  we  see  no  reason  to  restrict  ourselves  just to this  case:  the  consequences  are  similar  but  different.    \\

It is worth pointing out that
$$ \beta_1 \, \mathbf e = \gamma_1 \, \mathbf e \quad  {\rm  and }  \quad  \beta_2 \, \mathbf e^\dagger = \gamma_2 \, \mathbf e^\dagger$$

although  $\beta_1$  and  $\beta_2$  are   in   $ \mathbb C ( \mathbf i )$, while   $\gamma_1$ and $\gamma_2$  are  in   $ \mathbb C ( \mathbf j )$.  More  specifically,   given    $\beta_1  = Re \, \beta_1  + \mathbf i  \, Im \, \beta_1$,    the  equality
$$ \beta_1 \, \mathbf e = \gamma_1 \, \mathbf e$$
   is  true  if  and  only  if
$$ \gamma_1  =  Re \, \beta_1  - \mathbf j  \, Im \, \beta_1.$$ Similarly if   $\beta_2  = Re \, \beta_2  + \mathbf i  \, Im \, \beta_2$,    the  equality
$$ \beta_2 \, \mathbf e^\dagger = \gamma_2 \, \mathbf e^\dagger$$
 is  true  if  and  only  if
$$ \gamma_2  =  Re \, \beta_2  +  \mathbf j  \, Im \, \beta_2.  $$

   Altogether  decomposition          (\ref{descomposicion idempotente bicomplejos})  can be written  in  any of  the two  equivalent  forms:
$$  \mathbb B \mathbb C  =     \mathbb C ( \mathbf i ) \cdot  \mathbf e  +      \mathbb C    ( \mathbf i )   \cdot   \mathbf e^\dagger \, ;$$
$$  \mathbb B \mathbb C  =     \mathbb C ( \mathbf j ) \cdot  \mathbf e  +      \mathbb C    ( \mathbf j )   \cdot   \mathbf e^\dagger \, .$$

If  $  \mathbb B \mathbb C $  is  seen  as  a   $   \mathbb C    ( \mathbf i )  $--linear  (respectively,  a  $  \mathbb C    ( \mathbf j )  $--linear)  space  then  the  first    (respectively,  the  second)  decomposition  becomes  a  direct  sum.  \\

It  is usually  stated  that  the  idempotent  representation  is  unique  which  seems  to  contradict  to  our  two  formulas.  But  the  matter     is  that  each  of    (\ref{repr_idempotent})   and   (\ref{repr_idempotent_Cj})  is  unique  in  the  following  sense.    Assuming  that  a  bicomplex  number  $Z$   has   two  idempotent  representation  with  coefficients  in   $\C (\i)$:   $Z =  \beta_1  \, \mathbf e  + \beta_2 \,
\mathbf e^\dagger  =   \eta_1  \, \mathbf e  + \eta_2 \,
\mathbf e^\dagger$  with   $\beta_1 , \, \beta_2,   \, \eta_1, \,  \eta_2   \in  \mathbb C ( \mathbf i ) $ then it is  easy to show that   $\beta_1 = \eta_1$  and  $  \beta_2 =   \eta_2  $;  similarly      $Z =  \gamma_1  \, \mathbf e  + \gamma_2 \,
\mathbf e^\dagger  =   \xi_1  \, \mathbf e  + \xi_2 \,
\mathbf e^\dagger$  with   $\gamma_1 , \, \gamma_2,   \, \xi_1, \,  \xi_2   \in  \mathbb C ( \mathbf j ) $  implies  that     $\gamma_1 = \xi_1$  and  $  \gamma_2 =   \xi_2  $.  \\

Note that formula   (\ref{repr_idempotent})   gives  a  relation  between  the two  bases  $\{ 1 , \, \mathbf j \, \}$  and    $ \{ \mathbf e, \,  \mathbf e^\dagger \,
\}$  in  the    $
\mathbb C ( \mathbf i)$--linear  space         $\mathbb B  \mathbb C  =  \mathbb C^2( \mathbf i)$.  One has the following  transition  formula     using  the  matrix of change of
variables:
\begin{equation}
\left(     \begin{array}{l}  z_1  \\  \\  \\    z_2  \end{array}  \right)
=
\left(  \begin{array}{cc}    \displaystyle  \frac{1}{2}   &
\displaystyle  \frac{1}{2}  \\  \\
\displaystyle  -  \frac{1}{2 \,   \mathbf i}   &   \displaystyle
\frac{1}{2  \,  \mathbf i}
\end{array}    \right)
\cdot
\left(    \begin{array}{l}  \beta_1  \\  \\  \\    \beta_2  \end{array}
\right)   \, .
\end{equation}

Note  that   the  new  basis     $ \{ \mathbf e, \,  \mathbf
e^\dagger \, \}$   is  orthogonal  with  respect  to the
Euclidean  inner  product  in  $\mathbb C^2(   \mathbf i)$ given  by
$$  \displaystyle   \left\langle  \left( z_1 , \, z_2  \right)  \, , \,
\left(  u_1 ,  \,  u_2  \right)  \,  \right\rangle_{\mathbb C^2 ( \mathbf i) }
:=  z_1   \overline u_1  +  z_2 \overline u_2  \,  ,     $$

for all   $(z_1 ,z_2) , \;   (u_1 ,u_2) \in   \mathbb C ^2  ( \mathbf i) \, $.     Indeed,  since
$ \displaystyle  \mathbf e =  \left(  \frac{1}{2}
, \,   \frac{ \mathbf i}{2} \, \right)  $  and     $ \displaystyle
\mathbf e^\dagger =  \left(  \frac{1}{2} , \,  -  \frac{ \mathbf
i}{2} \, \right)  $  as  elements    in  $\mathbb C^2( \mathbf i)$,   we  have:
$$ \displaystyle     \left\langle  \mathbf e  , \, \mathbf e^\dagger  \,
\right\rangle_{\mathbb C ^2( \mathbf i)} =0  \,  \quad  {\rm and}  \quad
\left\langle   \mathbf e  , \, \mathbf e  \,
\right\rangle_{\mathbb C ^2( \mathbf i)}   =     \left\langle  \mathbf
e^\dagger  , \, \mathbf e^\dagger  \, \right\rangle_{\mathbb C ^2( \mathbf i)}
=   \frac{1}{2} \,  ,$$

and so   $ \{ \mathbf e, \,  \mathbf e^\dagger \,
\}$   is   an  orthogonal  but  not  orthonormal  basis  for
$\mathbb C ^2  ( \mathbf i)$.  As  a  consequence   one  gets:
\begin{equation}\label{estrella}
\displaystyle   |Z | =  \frac{1}{ \sqrt{2\, } \, } \sqrt{ | \beta_1 |^2 +
| \beta_2 | ^2 \, }  \, .
\end{equation}

\medskip

Formula (\ref{repr_idempotent_Cj}) can be interpreted analogously.  Now  we  have  $\mathbb B  \mathbb C  =  \mathbb C^2 ( \mathbf j)$
and  two  bases in  it  as  $  \mathbb C ( \mathbf j)$--linear   space  are  $ \{  1, \,   \mathbf i \, \}  $   and  again(!)   $ \{ \mathbf e, \,  \mathbf e^\dagger \,
\}$.

%

Note that   $\mathbf e$   and  $ \mathbf e^\dagger$  form  a  basis  in  both   $\mathbb C^2 ( \mathbf i)$  and  $ \mathbb C^2 ( \mathbf j)$.  The  basis    $ \{ \mathbf e, \,  \mathbf e^\dagger \,
\}$   remains orthogonal  in   $ \mathbb C^2 ( \mathbf j)$ with  transition   matrix     $ \displaystyle  \left(  \begin{array}{cc}    \displaystyle  \frac{1}{2}   &
\displaystyle  \frac{1}{2}  \\  \\
\displaystyle  -  \frac{1}{2 \,   \mathbf j}   &   \displaystyle
\frac{1}{2  \,  \mathbf j}
\end{array}    \right) $  and  if        $ Z =  \gamma_1 \, \mathbf e  + \gamma_2 \, \mathbf e^\dagger$,   $ \gamma_1 , \, \gamma_2 \in \mathbb C (\mathbf j) $,     then
\begin{equation}\label{relacion norma euclidian y componentes idemp Cj}
\displaystyle   |Z | =  \frac{1}{ \sqrt{2\, } \, } \sqrt{ | \gamma_1 |^2 +
| \gamma_2 | ^2 \, }  \, .
\end{equation}

A major advantage of  both  idempotent  representations  of  bicomplex  numbers  is  that  they   allow  to  perform  the operations  of addition, multiplication, inverse, square  roots, etc.  component-wise.  For  more  details  see,  e.g.,     \cite{RS}.

\medskip

Let us  come back to properties  {\it a})  and {\it b})  in Section  \ref{subsection 2.3}.  We  prove  first   {\it a}).  Indeed, take  $Z = z_1 \in \mathbb C ( \mathbf i ) $ and   $U = u_1 + u_2 \,  \mathbf j   =  ( u_1 - \mathbf i \, u_2 ) \, \mathbf e   +  ( u_1 + \mathbf i \, u_2 ) \, \mathbf e^\dagger   $, then
$$ \begin{array}{rcl}
| Z \cdot U |^2  & = &  | \, z_1 \, ( u_1 + u_2 \, \mathbf j ) \, |^2   =    | \, ( z_1 \,  u_1)   +   (z_1 \,  u_2) \, \mathbf j  \, |^2
\\  &  &  \\  & = &
 | \,  z_1 \,  u_1 \, |^2   +   | \, z_1 \,  u_2  \, |^2  \; = \; | \, z_1 \, |^2 \cdot | \, U \, |^2 \, ,
 \end{array} $$

and  the first  part  of   property  {\it a})  is  proved.  As to the second  part, take  $Z = x_1 + x_2 \, \mathbf j =  ( x_1 - \mathbf i \, x_2 ) \, \mathbf e  +    ( x_1 + \mathbf i \, x_2 ) \, \mathbf e^\dagger \in \mathbb C ( \mathbf j )   $,  then

\medskip

$ \begin{array}{lll}
| Z \cdot U |^2  & = &   | \,  \left(   ( x_1 - \mathbf i \, x_2 ) \, \mathbf e  +    ( x_1 + \mathbf i \, x_2 ) \, \mathbf e^\dagger  \right) \cdot
\\  & &  \\  & &   \qquad  \qquad  \cdot  \left(  ( u_1 - \mathbf i \, u_2 ) \, \mathbf e   +  ( u_1 + \mathbf i \, u_2 ) \, \mathbf e^\dagger  \right) \, |^2  \\  &  &  \\  & =  &
|    ( x_1 - \mathbf i \, x_2 )    ( u_1 - \mathbf i \, u_2 )  \mathbf e    +    ( x_1 + \mathbf i \, x_2 )   ( u_1 + \mathbf i \, u_2 )  \mathbf e^\dagger |^2
\\  &  &  \\  & =  &
\displaystyle \frac{1}{2}  \left(     |    x_1 - \mathbf i  x_2 |^2  \cdot  |   u_1 - \mathbf i     u_2 |^2     +    | x_1 + \mathbf i  x_2 |^2 \cdot |  u_1 + \mathbf i  u_2 |^2 \right)
\\  &  &  \\  & =  &  | Z |^2 \cdot | U  |^2 \, ,
 \end{array} $

thus  property  {\it a})  is  proved.  Finally,  in  order  to deal  with  property  {\it b}),     take  $ Z = x_1 + y_2 \, \mathbf i \, \mathbf j = ( x_1 + y_2) \, \mathbf e   +    ( x_1 - y_2) \, \mathbf e^\dagger \in \mathbb D$, then

\medskip
$ \begin{array}{lll}
| Z \cdot U |^2  & = &   | ( x_1 + y_2 \, \mathbf i \, \mathbf j  ) \cdot (  u_1 + u_2  \, \mathbf j ) |^2
\\  & & \\  & = &
\displaystyle \left|    ( x_1 + y_2 ) \cdot  (  u_1 -   \mathbf i   u_2)  \,  \mathbf e  +   ( x_1 - y_2 ) \cdot  (  u_1 +   \mathbf i  u_2)  \,  \mathbf e^\dagger   \right|^2
\\  & & \\  & = &
\displaystyle    \frac{1}{2} \,  \left(   ( x_1 + y_2 )^2  \cdot  |  u_1 -   \mathbf i   u_2 |^2     +   ( x_1 - y_2 )^2  \cdot  |  u_1 +   \mathbf i  u_2|^2     \right)
\\  & & \\  & = &
\displaystyle  | Z |^2 \cdot | U |^2   +     4   \, x_1 \, y_2 \, Re ( \mathbf i \, u_1 \, \overline u_2 )    \, .
 \end{array} $

\bigskip

  In the next sections we will    need the
mappings:
$$  \pi_{1, \mathbf i} , \;  \pi_{2, \mathbf i} :
\mathbb B \mathbb C  \to  \mathbb C( \mathbf i)   $$

 given  by
$$  \pi_{\ell , \mathbf i} (Z)  =  \pi_{\ell ,
\mathbf i}  ( \beta_1 \, \mathbf e  + \beta_2 \,
\mathbf e^\dagger \, )  :=   \beta_\ell  \in
\mathbb C ( \mathbf i)  \, .$$

\medskip

These maps are nothing but the    projections onto the ``coordinate axes"   in $\mathbb C^2(\mathbf i)$ with the  basis  $\{ \mathbf e , \, \mathbf e^\dagger \, \}$.  Completely analogously  one  has the mappings
$$  \pi_{1, \mathbf j} , \;  \pi_{2, \mathbf j} :  \mathbb B \mathbb C  \to  \mathbb C( \mathbf j)   $$

given  by
$$  \pi_{\ell , \mathbf j} (Z)  =  \pi_{\ell , \mathbf j}  ( \gamma_1 \, \mathbf e  + \gamma_2 \,
\mathbf e^\dagger \, )  :=   \gamma_\ell  \in  \mathbb C ( \mathbf j)  \, ,$$

which  are  now    the    projections onto the coordinate axes in $\mathbb C^2(\mathbf j)$ with the  same  basis  $\{ \mathbf e , \, \mathbf e^\dagger \, \}$.

\medskip

Finally,   using  the  cartesian   representations  for  any  $Z \in  \mathbb B \mathbb C$:
$$\begin{array}{rcl}
Z  &    =  &  (x_1 + \mathbf i  \, x_2 )  +  (y_1 + \mathbf i \, y_2 ) \, \mathbf j  =  z_1  +  z_2   \mathbf j
\\  &  &  \\   &  =  &
( x_1  +  \mathbf j \,  y_1  )  +  (  x_2  +  \mathbf j \, y_2 ) \, \mathbf i   =  \eta_1  +   \eta_2  \, \mathbf i  \, ,
\end{array}$$

we  define  two more projections
$$    \Pi_{1, \mathbf i} , \;  \Pi_{2, \mathbf i} :  \mathbb B \mathbb C  \to  \mathbb C( \mathbf i)  \qquad  {\rm and}  \qquad     \Pi_{1, \mathbf j} , \;  \Pi_{2, \mathbf j} :  \mathbb B \mathbb C  \to  \mathbb C( \mathbf j)   $$

given  by
$$  \Pi_{\ell , \mathbf i} (Z)  =  \Pi_{\ell , \mathbf i}  (   z_1  +  z_2   \mathbf j     \, )  :=   z_\ell  \in  \mathbb C ( \mathbf i)  \, ,$$

and
$$  \Pi_{\ell , \mathbf j} (Z)  =  \Pi_{\ell , \mathbf j}  (   \eta_1  +  \eta_2   \mathbf i     \, )  :=   \eta_\ell  \in  \mathbb C ( \mathbf j)  \, .$$

\medskip

The  properties  of  the  operations  with  bicomplex  numbers immediately imply the following equalities
for  any   $   \, Z , \, W \in \bc$:

\medskip

\begin{enumerate}

\item[(a)]  $ \pi_{ \ell , \i }  (Z W ) =  \pi_{ \ell , \i }  (Z)   \,   \pi_{ \ell , \i }   (W) $;

\medskip

\item[(b)]  $ \pi_{ \ell ,  \, \j }  (Z W) =  \pi_{ \ell , \, \j }   (Z)    \,  \pi_{ \ell , \,  \j }  (W) $;

\medskip

\item[(c)]  $ \forall \, \lambda \in \C ( \i) , \;   \;    \pi_{ \ell , \i }  (\lambda )  = \lambda ; \;   \;     \Pi_{ 1 , \i }  (\lambda)  = \lambda , \;   \;     \Pi_{ 2 , \i }  (\lambda )  =0$;

\medskip

\item[(d)]  $ \forall \, \mu \in \C ( \,  \j) , \;  \;     \pi_{ \ell , \,  \j }  (\mu )  = \mu ; \;    \;     \Pi_{ 1 , \,  \j }  (\mu)  = \mu   , \;    \;     \Pi_{ 2 , \,  \j }  (\mu )  =0$;

\medskip

\item[(e)]  $ \pi_{1 , \i } = \Pi_{1 , \i }  - \i  \,  \Pi_{2 , \i } \, ;    \;   \;      \pi_{2 , \i } = \Pi_{1 , \i }   +    \i   \, \Pi_{2 , \i } \, ;  $

\medskip

\item[(f)]  $ \pi_{1 ,   \, \j } = \Pi_{1 , \,  \j }  - \j    \,    \Pi_{2 , \j } \, ;    \;  \;      \pi_{2 ,  \,  \j } =    \Pi_{1 ,  \, \j }   +    \j     \,   \Pi_{2 ,  \,   \j } \, $;

\medskip

\item[(g)]   $  \displaystyle  \Pi_{1 , \i }  = \frac{1}{2}  \left(    \pi_{1, \i }  + \pi_{2 , \i } \right)  ; \; \;  \Pi_{2 , \i }  = \frac{\i }{2}  \left(    \pi_{1, \i }  -   \pi_{2 , \i } \right) $;

\bigskip

\item[(h)]   $  \displaystyle  \Pi_{1 , \,  \j  }  = \frac{1}{2}  \left(    \pi_{1, \,  \j }  + \pi_{2 , \j } \right)  ; \; \;  \Pi_{2 , \j }  = \frac{\j }{2}  \left(    \pi_{1,  \,  \j }  -   \pi_{2 , \,  \j } \right) $.

\end{enumerate}

\medskip

We  leave the verification   of  these  properties   to  the   reader.

\bigskip

\section{A  partial  order  on  $\D$  and a hyperbolic--valued  norm} \label{section hyperbolic norm}

Let us describe some properties that  hyperbolic numbers inherit  from    bicomplex ones.  First of all we note that  both  $\mathbf e $  and   $\mathbf e^\dagger $   are   hyperbolic numbers. More importantly, the idempotent representation of any hyperbolic number  $\alpha = a + b \, \mathbf k$  is
$$  \alpha =  \nu \, \mathbf e  +  \mu \, \mathbf e^\dagger \, , \quad
\nu  , \, \mu  \in  \mathbb R ,$$

with  $\nu = b + a$,  $\mu  = b-a$. It is also immediate  to  see  that
$$  \displaystyle  \mathbb D^+  =    \left\{  \nu  \, \mathbf e  +  \mu \,
\mathbf e^\dagger   \mid   \nu,  \,  \mu  \geq 0  \,  \right\}  \, .$$

Thus  positive  hyperbolic numbers   are  those  hyperbolic numbers  whose  both  idempotent  components  are  non negative,   that somehow explains the  origin of the name.

\medskip

\begin{figure}[h]
\centering
\includegraphics[trim = 0cm 10cm 0cm 1cm, clip, scale=0.4]{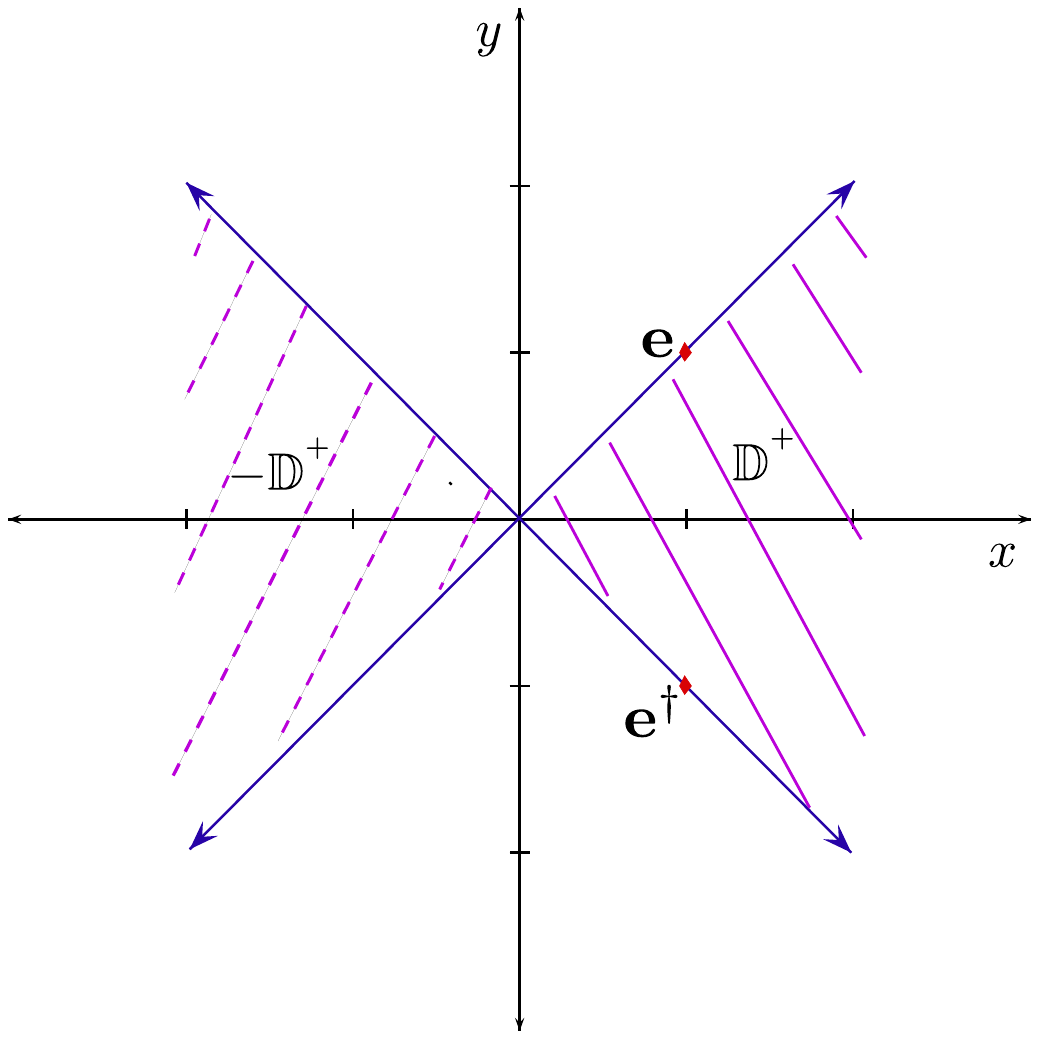}
\caption{}
\end{figure}

\bigskip

On Fig. 1  the points  $(x,y)$  correspond  to the  hyperbolic  numbers   $Z=x + \k y $.  One  sees that,  geometrically,  the  hyperbolic  positive  numbers  are  situated in  the     quarter  plane   denoted  by  $\D^+$.  The  quarter  plane   symmetric  to   it    with  respect  to  the   origin  corresponds  to  the  ``negative"  hyperbolic  numbers, i.e.,   to  those  which  have  both  idempotent  components  negative.   The  rest  of  points  corresponds   to  those  hyperbolic  numbers  which  can  not  be  called   either  positive  or  negative.      \\

Let  us now  define  in  $\D$     the  following  binary  relation:  given  $\alpha_1$,  $\alpha_2  \in  \D$,  we  write  $  \alpha_1  \lessdot    \alpha_2 $  whenever   $ \alpha_2  - \alpha_1  \in \D^+$.  It  is  obvious  that  this  relation  is  reflexive,  transitive  and  antisymmetric, and that therefore it  defines  a  partial  order  on  $\D$.   With  abuse  of  notation,  we  will say  that  $\alpha_2 $  is  $\D$--larger   ($\D$--less  respectively)  that  $\alpha_1$  if   $\alpha_1 \lessdot  \alpha_2$  ($\alpha_2 \lessdot \alpha_1$ respectively).  \\

Observe that  given  $r , \, s \in \R$,   then   $r \leq   s$  if  and  only  if   $r \lessdot s$,  that  is,     the partial order $\lessdot$   is  an  extension of the  total order  $  \leq  $  in  $\R$.  In  fact,   $\lessdot$   defines   a  total  order on any  line  through  the   origin  in  $\D$.  

\medskip

\begin{figure}[h]
\centering 
\includegraphics[trim = 0cm 10cm 0cm 1cm, clip, scale=0.4]{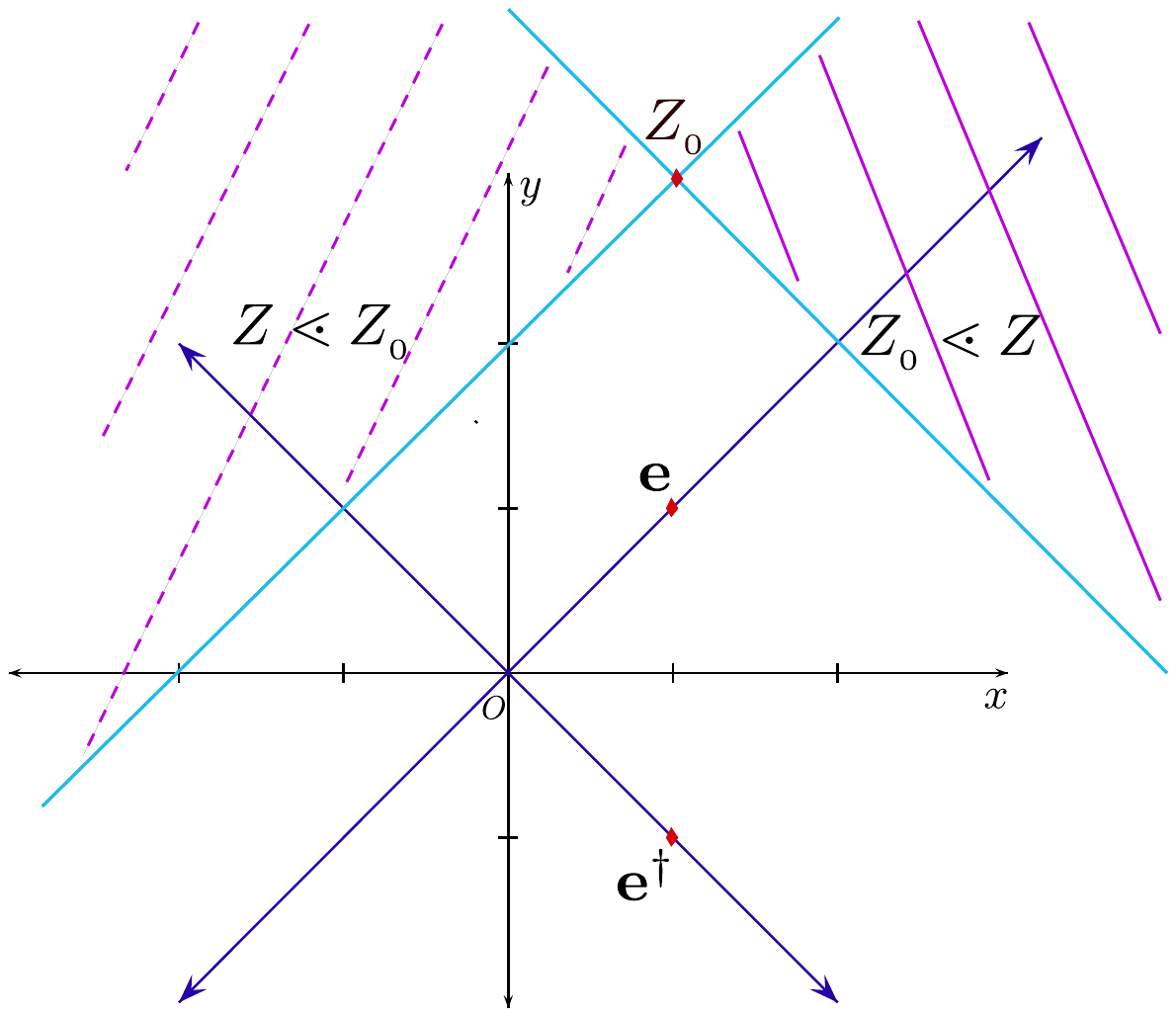}
\caption{}
\end{figure}

\bigskip

On Figure 2,   \,   $Z_0 =x_0 + \k y_0$  is  an  arbitrary  hyperbolic  number,  and  one can see that the entire  plane is divided into:  the  quarter  plane of hyperbolic  numbers  which are  $\D$--greater  that  $Z_0$    ($Z  \gtrdot    Z_0  $);  the  quarter  plane of  hyperbolic  numbers  which   are  $\D$--less  than  $Z_0$  ($  Z \lessdot  Z_0$);  and  the  two  quarter  planes   where  the  hyperbolic  numbers  are  not  $\D$--comparable  with  $Z_0$  (neither  $Z \gtrdot Z_0$ nor  $Z \lessdot Z_0$  holds).

\medskip

\begin{subremark}
While the notion of positive hyperbolic number may not appear immediately intuitive, one can see its relation to a
two-dimensional space time with the Minkowksi metric. Specifically, consider a one-dimensional time, say along the $x-$axis, and a one dimensional space along the $y-$axis. Embed both axes in the hyperbolic space $\mathbb D.$ Finally, assume the speed of light to be normalized to one. Now, the set of zero-divisors with $x>0$ represents the future of a ray of light being sent from the origin towards either direction. Analogously, the set of zero-divisors with $x<0$ represents the past of the same ray of light. Thus, positive hyperbolic numbers are nothing but the future cone, i.e. they describe the events that are in the future of the origin, while the negative hyperbolic numbers are the past cone, i.e. they are the events in the future of the origin. Within this interpretation, we see how the set of hyperbolic numbers larger than a given hyperbolic number $Z$ represent the events in the future of the event represented by $Z$.
\end{subremark}


 Given a  subset  $\mathcal A  $  in  $\D$, we can now    define a notion of  $\D$--upper  and  $\D$--lower  bound,  as well as the notions of  a set being  $\D$--bounded  from  above,  from  below,  and  finally of  a   $\D$--bounded  set. If $\mathcal A \subset \D$ is a set $\D$-bounded from above, we define the notion of its supremum, $sup_{\D} \mathcal A$ as usual to be the least upper bound for $\mathcal A.$ However, it is immediate to note that one can find a more convenient expression for it, as indicated in the following remark.

\begin{subremark}
Given  a    set   $\mathcal A \subset \D$,  $\D$--bounded  from  above,      consider    the sets  $ \mathcal A_1 := \{ \, a_1 \mid  a_1 \e + a_2 \edag \in \mathcal A \, \} $  and  $ \mathcal A_2 := \{ \, a_2 \mid  a_1 \e + a_2 \edag \in \mathcal A \, \} $.  Then  the  $\sup_\D  \mathcal A$  is given  by
$$  {\rm sup}_\D  \mathcal A  :=  \sup \mathcal A_1 \cdot \e  + \sup \mathcal A_2 \cdot \edag.$$
\end{subremark}

\medskip

We have previously defined the hyperbolic modulus $| Z |_\k^2$ of any  $Z  \in  \bc$   by the formula   $| Z |_\k^2 :=  Z \cdot Z^\ast$.  Thus,  writing   $Z  =  \beta_1 \e + \beta_2 \edag$,  one has that  $| Z |_\k^2  =   | \beta_1 |^2 \e + | \beta_2 |^2  \edag  $,  and  taking  the  positive  squares  roots  of  each  coefficient  we  have  that
$$   | Z |_\k  :=    | \beta_1 | \e + | \beta_2 |  \edag   . $$

That is, we have a  map
$$  | \cdot  |_\k :  \bc  \longrightarrow  \D^+$$

\medskip

with  the  following  properties:

\begin{enumerate}

\item[({\sc i})]   $  |Z |_\k  =0$  if  and  only  if  $Z=0$;

\medskip

\item[({\sc ii})]  $ |Z \cdot W |_\k  =  |Z |_\k \cdot |W |_\k $  for any  $Z,$  $W \in \bc$;

\medskip

\item[({\sc iii})]     $ | Z +W |_\k   \lessdot  | Z |_\k  +  | W |_\k$.

\end{enumerate}

The  first  two  of  them  are clear.  Let's  prove    ({\sc iii}).
$$  \begin{array}{rcl}
|Z +   W |_\k  & = &  |   ( \beta_1  + \nu_1)   \cdot \e  +  ( \beta_2  + \nu_2) \cdot \edag |_\k
\\  &  &  \\  & = &
| \beta_1  +   \nu_1 |  \cdot   \e   +   | \beta_2  +  \nu_2 |   \cdot  \edag  \\  &  &    \\  &   \lessdot  &
( | \beta_1 | + | \nu_1 | ) \cdot \e   +  ( | \beta_2 | + | \nu_2 | ) \cdot \edag   \\  &   &  \\   &  =  &          |Z |_\k   + | W |_\k .
\end{array}   $$

It is on the basis of these properties that we will say that  $ | \cdot |_\k$ is the  hyperbolic--valued  ($\D$--valued)   norm  on  the  $\bc$--module  $\bc$.

\medskip

It is instructive to compare   ({\sc ii})  with   (\ref{desigualdad norma euclideana})  where  the  norm of  the  pro\-duct  and  the  product  of  the  norms  are  related   with  an  inequality.  We  believe  that  one could say   that  the  hyperbolic  norm  of  bicomplex  numbers  is  better  compatible  with  the  algebraic  structure  of  the  latter  although, of course,  one has to allow hyperbolic values.

\begin{subremark}  \label{very usseful remark}
\end{subremark}
\begin{enumerate}

\item[(1)]   Since  for  any  $Z\in \bc$ it holds  that  
\begin{equation}\label{desig norma hyperb less Euclideana}
 | Z  |_\k  \lessdot  \sqrt{2}  \cdot |Z|
 \end{equation}

\noindent
with  $| Z|$  the  Euclidean norm of $Z$,  one has:
$$| Z \cdot W |_\k  \lessdot  \sqrt{2} \cdot |Z| \cdot |W |_\k  .$$

\noindent In  contrast   with  property   ({\sc ii}),   the  above  inequality  involves  both  the  Euclidean  and  the  hyperbolic  norms.

\medskip

\item[(2)]  Take   $\gamma $  and  $\nu $  in $\D^+$, then  clearly
\begin{equation}\label{desig hyperb positiv implies desig Eucl}
\gamma   \lessdot     \nu   \quad \Longrightarrow   \quad     |  \gamma  |  \leq | \nu | .
\end{equation}

\medskip

\item[(3)]  Note that the definition of hyperbolic norm for a bicomplex number $Z$ does not depend on the choice of its idempotent representation.
We  have  used,  for   $Z \in \bc$,  the  idempotent  representation  $ Z=    \beta_1 \e + \beta_2 \edag$,    with   $\beta_1$,  $\beta_2  \in \C ( \i)$.  If  we  would  have  departed  from the  idempotent  representation   $Z=  \gamma_1 \e  + \gamma_2 \edag$,  with  $\gamma_1$,  $\gamma_2  \in  \C ( \, \j)$ then we  would  have    arrived  at  the  same  definition   of  the  hyperbolic norm  since   $| \beta_1 | = | \gamma_1| $  and  $| \beta_2 | = | \gamma_2| $.

\end{enumerate}

\medskip

The comparison of the  Euclidean  norm   $|Z|$  and the $\D$-valued  norm  $|Z|_\k$ of a bicomplex number gives
\begin{equation}\label{relation hyperbolic Euclidean norms}
 \displaystyle  \left|  |Z|_\k \right|  =  \frac{1}{2} \, \sqrt{ \,  | \beta_1 |^2  + |  \beta_2 |^2 \,   }  =  |Z |
\end{equation}

where the left--hand  side  is the  Euclidean norm  of  a  hyperbolic  number.

%

\begin{subDn}
A  sequence  of  bicomplex  numbers  $\{ Z_n \}_{n \in \N}$  converges  to  the  bicomplex  number  $Z_0$  with  respect  the  hyperbolic--valued  norm  $|   \cdot |_\k$   if   for all   $    \varepsilon >0$  there  exists  $N \in  \N$  such that  for   all  $  \,  n \geq N$  there  holds:
$$  |  Z_n  - Z_0 |_\k  \lessdot \varepsilon .  $$
\end{subDn}

From inequality  (\ref{desig norma hyperb less Euclideana})  and  equality  (\ref{relation hyperbolic Euclidean norms}),  it     follows  that  a   sequence  $\{ Z_n \}_{n \in \N}$  converges  to  the  bicomplex  number  $Z_0$  with  respect  the  hyperbolic--valued  norm  if  and  only  if  it  converges  to  $Z_0$  with  res\-pect    to    the  Euclidean norm, and so even though the two norms cannot be compared, as they take values in different rings, one still obtains   the same notion of convergence.

\medskip

It  is  possible  to  give  a  precise  geometrical  description  of the set  of  bicomplex  numbers  with  a  fixed  hyperbolic  norm,   that  is   we  want  to  introduce  the  ``sphere  of  hyperbolic  radius $\gamma_0$ inside  $\bc$".      In other  words,  given    $\gamma_0 = a_0 \e + b_0 \edag$  (see  Figure 3)   we  are looking for the set
$$  \displaystyle \mathbb S_{\gamma_0}:=     \{  Z \in \bc \mid  |Z|_\k  = \gamma_0 \, \}   . $$

First  note  that,  if  one  of  $a_0$ or $b_0$ is  zero, let's say $\gamma_0=  a_0 \cdot \e  $, then
$$\mathbb S_{\gamma_0} =  \{Z=  \beta_1 \cdot \e \mid | \beta_1 | =  a_0 \, \},$$

and this set is  a  circumference  in  the  real  two--dimensional  plane   $\bc_\e$ with center at the origin  and radius  $\displaystyle  \frac{ a_0 }{\sqrt{2}}    =  |  a_0  \cdot    \e |   $.  Similarly,  if  $\gamma_0  =  b_0 \cdot \edag$,  the  set  $\mathbb S_{\gamma_0}$  is  a  circumference  in $\bc_\edag$    with center at   the origin  and radius  $\displaystyle  \frac{ b_0 }{\sqrt{2}}    =   |  b_0 \cdot \edag |  $.

\medskip

If    $a_0 \neq 0 $  and  $b_0 \neq 0$  then   the  intersection  of  the    hyperbolic  plane  $\D$  and  the  sphere $\mathbb S_{\gamma_0}$  consists  of  exactly   four  hyperbolic  numbers:       $\pm a_0 \e  \pm b_0 \edag$,  where  the  plane  $\D$   touches  the  sphere  tangentially.   What is more,  in this case  the  sphere   $\mathbb S_{\gamma_0}$  is  the  surface  of  a  three--dimensional   torus,  which is  obtained  by   taking  the  cartesian product  of  the  two  circumferences,  one  in  the  plane  $\bc_\e$ (see Figure 5), the other in the plane $\bc_\edag$  (see Figure 6), both centered at zero and with radii    $\displaystyle  \frac{ a_0 }{\sqrt{2}} $  and     $\displaystyle  \frac{ b_0 }{\sqrt{2}} $  respectively:
$$  \displaystyle  \mathbb S_{\gamma_0} =  \left\{   \beta_1  \cdot \e +  \beta_2  \cdot \edag  \;  :   |  \beta_1 | = a_0 , \ ,\;  | \beta_2 | = b_0 \,   \right\} \, .$$

See  Figures  below.

\bigskip

\begin{figure}[h]
\centering 
\includegraphics[trim = 0cm 19cm 0cm 1cm, clip, scale=0.75]{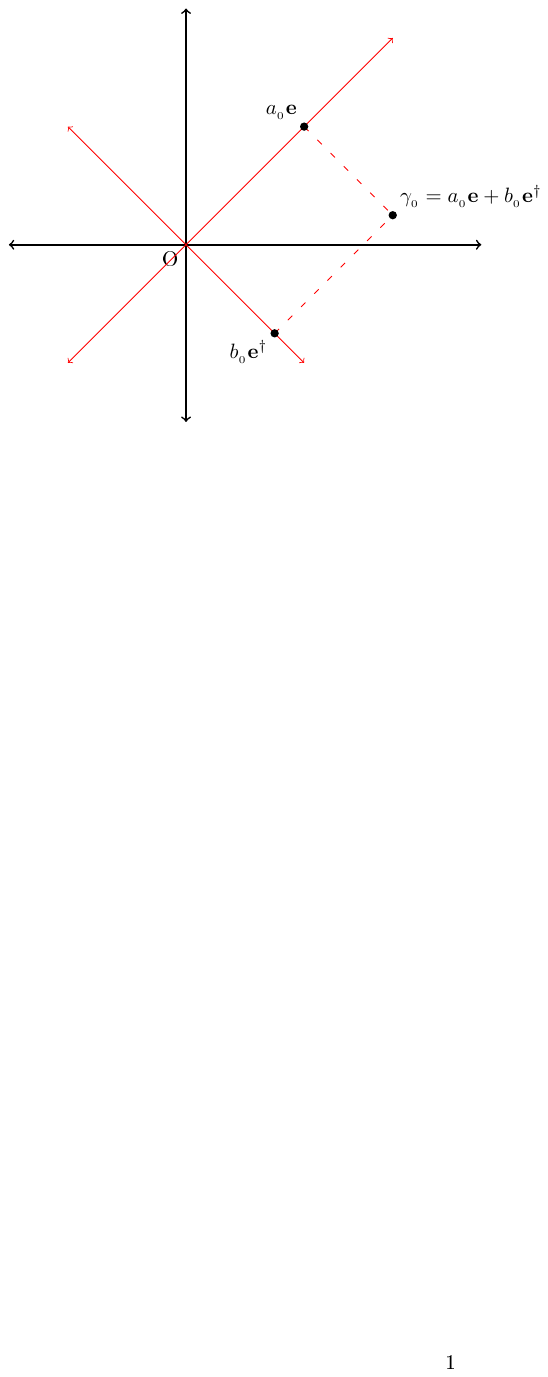}
\caption{}
\end{figure}

\bigskip

\bigskip

\begin{figure}[h]
\centering 
\vspace{3.5cm}
\includegraphics[trim = 0cm 21cm 0cm 1cm, clip, scale=0.75]{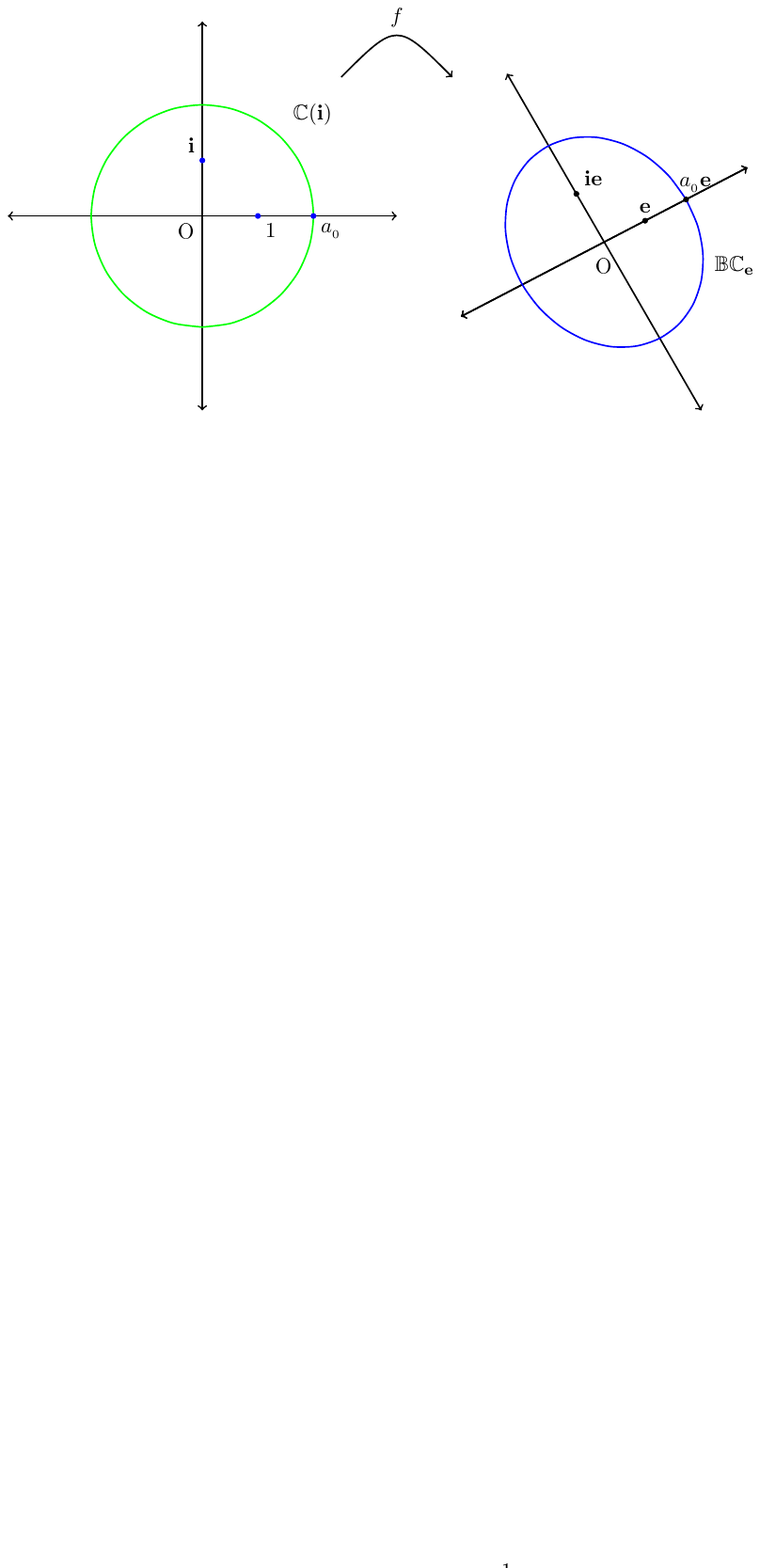}
\caption{}
\end{figure}

\bigskip

\bigskip

\begin{figure}[h]
\centering 
\includegraphics[trim = 0cm 21cm 0cm 1cm, clip, scale=0.75]{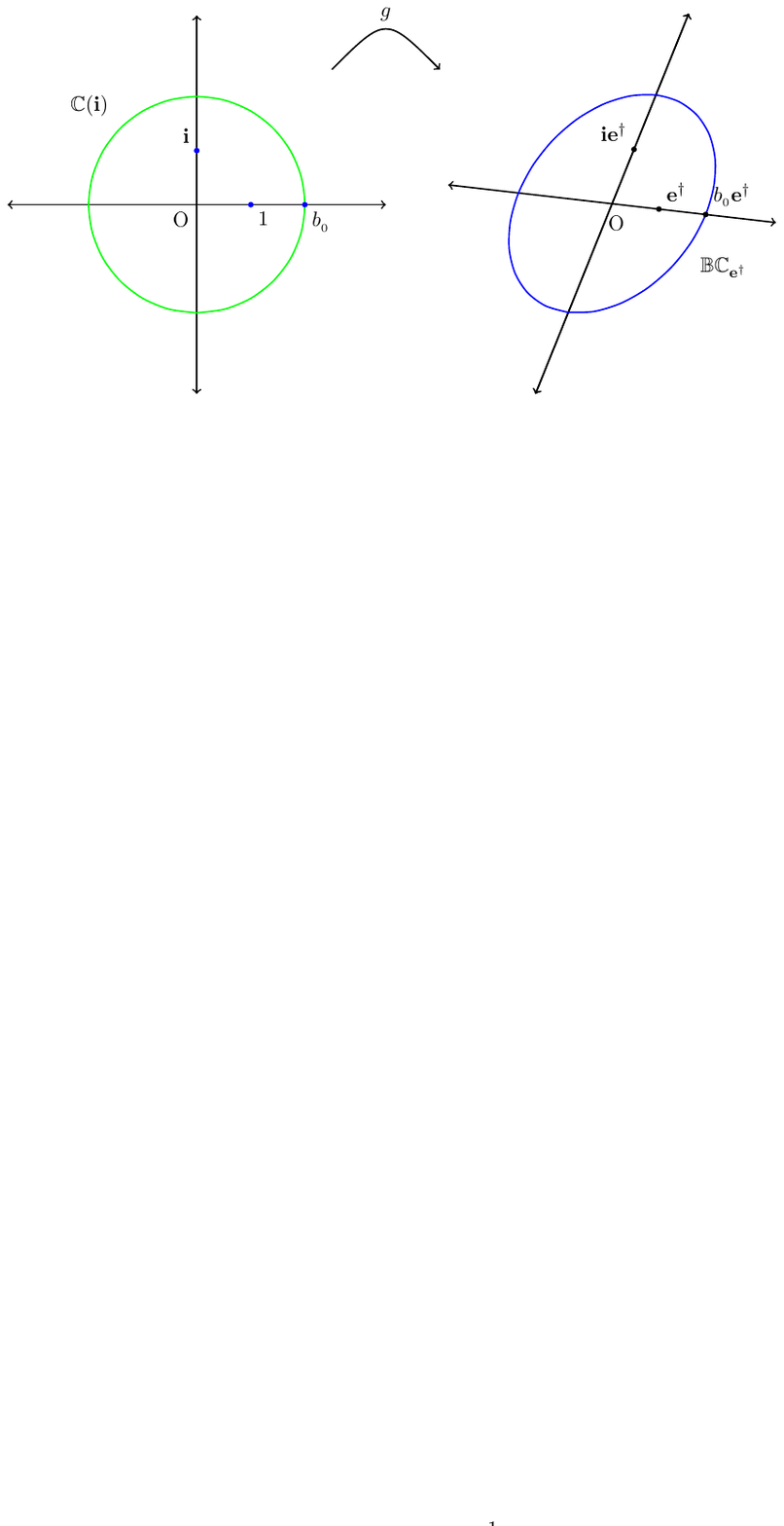}
\caption{}
\end{figure}

\bigskip

\bigskip

\begin{figure}[h]
\vspace{3.5cm}
\centering 
\includegraphics[trim = 0cm 0cm 0cm 0cm, clip, scale=2.4]{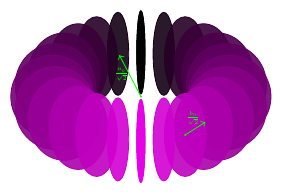}
\caption{}
\end{figure}

\bigskip

\newpage


 \chapter{Bicomplex  functions  and   matrices}
          \label{sec:pmat}


\section{Bicomplex holomorphic  functions}

The notion of  bicomplex holomorphic functions has been  introduced  a  long  ago  and  we  refer  to  the  reader  to the introduction of the book  \cite{Price}  and  the  book  by itself  can serve as  a  first  reading  on  this  subject.   For the latest  developments  see    \cite{LSSV2013}     and  \cite{LSSVbook}.   
In this Section we    present    a   summary    of the basic  facts  of  the  theory  which  will  be  used  in the last  chapter.  \\

The {\em derivative} $F'(Z_0)$ of the function $F : \Omega \subset \bc \to \bc$ at a point $Z_0\in \Omega$
  is   defined  to  be   the limit, if it exists,
  \begin{equation}
    \label{bicomplex_derivative}
    F'(Z_0):=     \lim_{\mathfrak{S}_0\not\ni H\to 0}\frac{F(Z_0+H)-F(Z_0)}{H} ,
  \end{equation}
  for $Z_0 $ in the domain of $F$ such that $H = h_1  + \j  h_2  $ is an invertible
  bicomplex number. In this case, the function $F$ is called {\em
    derivable} at $Z_0$.  \\
    
In analogy  with holomorphic functions in one complex variable  the  Cauchy--Riemann  type  conditions  arise where  now  the complex, not real,  
partial derivatives participate.   More explicitly, if we consider a bicomplex function $F=f_1+\j f_2$ derivable at $Z_0$,   then we have  that
the   complex partial derivatives
$$ \displaystyle  F'_{z_1}(Z_0)  = \lim_{ h_1 \to 0} \frac{ F( Z_0 + h_1 ) -   F (Z_0) }{  h_1}   ,   $$

$$ \displaystyle  F'_{z_2}(Z_0)  =  \lim_{ h_2 \to 0} \frac{ F( Z_0 + \j   h_2 ) -   F (Z_0) }{  h_2}   ,   $$

exist   and   verify the identity:
\begin{equation}
\label{complexi_derivative}
F'(Z_0) = F'_{z_1}(Z_0) = -\j F'_{z_2}(Z_0)\,,
\end{equation}
which is equivalent to the     {\em complex Cauchy-Riemann
system} for $F$ (at $Z_0$):
\begin{equation}
      \label{CR_complex}
      \left(  f_1 \right)'_{z_1}(Z_0) = \left(  f_2 \right)'_{z_2}(Z_0)\,,
      \qquad \left(  f_1 \right)'_{z_2}(Z_0) = -  \left( f_2 \right)'_{z_1}(Z_0) \,.
    \end{equation}

Of course,  if $F$ has bicomplex
derivative at each point of $\Omega$, we will say that $F$ is a {\em bicomplex holomorphic}, or
$\bc$-holomorphic, function.  \\  

For a $\bc$--holomorphic function $F$,  formulas   (\ref{CR_complex})   imply  that

\begin{equation}
\label{eq:baba1}
\frac{ \partial F}{  \partial \bar{z}_1} (Z_0) = \frac{ \partial F}{ \partial \bar{z}_2}(Z_0) = 0\,,
  \end{equation}
  i.e.,
  \begin{equation}
    \label{eq:baba2}
    \frac{ \partial f_1}{ \partial   \bar{z}_1}(Z_0) = \frac{   \partial f_1}{    \partial \bar{z}_2}(Z_0) =
    \frac{ \partial f_2}{    \partial  \bar{z}_1}(Z_0) = \frac{ \partial f_2}{   \partial \bar{z}_2}(Z_0) = 0\,,
  \end{equation}
  where the symbols $\displaystyle\frac{\partial }{\partial \bar{z}_1}$ and
  $\displaystyle\frac{ \partial }{  \partial \bar{z}_2}$ are the commonly used formal
  operations on    functions of $z_1$ and $z_2$.

\medskip

This means in particular   that $F   $ is     holomorphic    with
respect to $z_1$ for any $z_2$ fixed and $F$ is holomorphic with
respect to $z_2$ for any $z_1$ fixed. Thus, see for
instance~\cite[pages 4-5]{krantz}, $F$ is holomorphic in the classical
sense of two complex variables. This implies immediately many quite
useful properties of $F$, in particular, it is of class
$\mathcal{C}^\infty(\Omega)$.

\medskip

In one complex variable there exist two  mutually  conjugate  Cauchy--Riemann  operators  which  characterize  the  usual  holomorphy.  Since  in  $\bc$  there  exist  three conjugations,   it turns  out  that   there exist  four  bicomplex    operators  which   characterize in a  similar  fashion   the  $\bc$--holomorphy.  They are:

\begin{equation}\label{BC_diffop}
\begin{array}{lr}
\displaystyle    \frac{\partial }{ \partial Z}  :=\frac{1}{2}\left(   \frac{ \partial }{  \partial z_1}   -   \j   \frac{ \partial }{ \partial z_2}  \right),     &   \displaystyle    \qquad
\frac{ \partial }{ \partial Z^\dagger}   :=    \frac{1}{2}   \left(   \frac{\partial }{ \partial z_1}  +   \j    \frac{ \partial }{ \partial z_2}   \right)   ,  \\   \\
\displaystyle   \frac{  \partial }{ \partial \bar{Z}}      :=   \frac{1}{2}\left(   \frac{\partial }{  \partial \bar{z}_1}   -    \j   \frac{ \partial }{ \partial \bar{z}_2}  \right)   ,     &   \displaystyle   \qquad
\frac{\partial }{ \partial Z^\ast}   :=    \frac{1}{2}    \left(   \frac{ \partial }{  \partial \bar{z}_1}   +    \j   \frac{ \partial }{ \partial \bar{z}_2}\right)\,.
\end{array}
\end{equation}

\begin{subtheorem}\label{big}
    Given $F\in\mathcal{C}^1(\Omega,\bc)$,   it  is $\bc$--holomorphic  if  and  only  if
\begin{equation}\label{Oct3_2}
  \frac{ \partial F}{ \partial  Z^\dagger}(Z) = \frac{ \partial F}{  \partial \bar{Z}}(Z) = \frac{  \partial F}{ \partial Z^\ast}(Z) = 0\,
\end{equation}
  holds on $\Omega.$  Moreover, when it is true one has:   $  \displaystyle  F' (Z) = \frac{ \partial F}{ \partial Z } (Z) $.
\end{subtheorem}

\medskip

Let us see now  how   the  idempotent  representation  of  bicomplex  numbers  becomes crucial  in a  deeper  understanding of the nature of  $\bc$--holomorphic functions.    
Take a bicomplex function $F:\Omega\subset   \bc   \to  \bc$ with $\Omega$ being a domain. We write all the
bicomplex numbers involved in  idempotent form, for instance,
\begin{align*}
  Z &= \beta_1\e + \beta_2\edag = (\ell_1 +\i m_1)\e + (\ell_2 +\i m_2)\edag\,,\\
  F(Z) &= G_1(Z)\e + G_2(Z)\edag\,,\\
  H &= \eta_1\e + \eta_2\edag = (u_1 +\i v_1) \e + (u_2 +\i v_2)\edag\,.
\end{align*}

\noindent Let us introduce the sets
$$\Omega_1:=\{\beta_1\,\big|\,
\beta_1\e + \beta_2\edag \in \Omega \} \subset \C(\i)
$$
and
$$\Omega_2:=\{\beta_2\,\big|\, \beta_1\e + \beta_2\edag \in \Omega \}
\subset \C(\i).
$$
It is direct to prove that $\Omega_1$ and $\Omega_2$
are domains in $\C(\i)$.

\medskip

\begin{subtheorem}
  \label{7.3}
  A bicomplex function $F  =   G_1 \e  + G_2  \edag  :  \Omega\subset   \bc  \to   \bc$ of class
  $\mathcal{C}^1$ is $ \bc$--holomorphic   if and only if   the  following  two  conditions  hold:

  \begin{itemize}
  \item[(I)] The component $G_1$, seen as a $\C(\i)$--valued function
    of two complex va\-ria\-bles $(\beta_1,\beta_2)$ is holomorphic; what is more,
    it does not depend on the variable $\beta_2$ and thus $G_1$ is a
    holomorphic function of the variable $\beta_1$.

  \item[(II)] The component $G_2$, seen as a $\C(\i)$--valued function
    of two complex va\-ria\-bles $(\beta_1,\beta_2)$ is holomorphic; what is more,
    it does not depend on the variable $\beta_1$ and thus $G_2$ is a
    holomorphic function of the variable $\beta_2$.

  \end{itemize}
\end{subtheorem}

\medskip

\begin{subCy}
  Let $F$ be a $\bc$--holomorphic function in $\Omega$, then $F$  is  of  the  form     $F(Z) =  G_1 (\beta_1) \e   +  G_2 (\beta_2) \edag$   with   $Z=\beta_1\e + \beta_2\edag\in\Omega$  and  its  derivative  is  given  by
$$
  F'(Z) = \e \cdot G'_1(\beta_1) +    \edag  \cdot G'_2(\beta_2)\,   ,
  $$
or equivalently:
$$
  F'(z_1  + \j z_2 ) =     \e    \cdot  G'_1(z_1  - \i  z_2  ) + \edag   \cdot  G'_2(z_1 +  \i   z_2)\,  .
  $$

\end{subCy}

\bigskip


\section{Bicomplex matrices}

We  will denote  by   $\bc^{m \times n}$  the  set  of   $m \times n$  matrices  with  bicomplex  entries.   For any such matrix $A=(a_{\ell j})\in\bc^{m \times n}$, it is possible to consider both its  cartesian  and  idempotent  representation, which are obtained by accordingly decomposing each of its entries   so that, for example, we have
$$  A  =  \mathscr A_{1 , \i }  \e  + \mathscr A_{2 , \i }    \edag   =  \mathscr A_{1 , \,  \j }  \e  + \mathscr A_{2 , \,  \j }    \edag    $$

where   $\mathscr A_{ 1 , \i } $,  $ \mathscr A_{ 2 , \i }   \in  \C^{m \times n} ( \i  ) $    and
$\mathscr A_{ 1 ,   \, \j } $,  $ \mathscr A_{ 2 , \,   \j }   \in  \C^{m \times n}  ( \,  \j  ) $.   \\

Of  course  the  set   $\bc^{m \times n } $  inherits  many    structures  from   $\bc$.   It   is  obviously  a   $\bc$--module  (a concept that will be discussed in detail  in   the  next  chapter) such that
$$  \begin{array}{rcl}
\bc^{m \times n }  & = &   \C^{m \times n} ( \i  ) \cdot  \e   +     \C^{m \times n} ( \i  )  \cdot  \edag
\\  &  &   \\   &  =  &
 \C^{m \times n} ( \,   \j  )   \cdot  \e   +    \C^{m \times n} ( \,  \j  )  \cdot  \edag \, ,
 \end{array}  $$

 where  the  summands  are  $\bc$--submodules    of   $  \bc^{m \times n }  $.   In  particular,  given   $ \mathscr B   \in    \C^{m \times n} ( \i  ) $  then   $ \mathscr B \cdot \e  \in     \C^{m \times n} ( \i  )  \cdot \e  $   and   $  \edag  \cdot (  \mathscr B \cdot \e )  =  0_{m \times n } $.

\medskip

As  in  the  scalar  case,  the  operations  over  the  matrices   in  the  idempotent  decomposition  can be  realized  component-wise (though keeping in mind the  non--commutativity  of  matrix  multiplication).

\medskip

\begin{subPn}
Let $A$ be an $n \times n$ bicomplex matrix

$$  A  =  \mathscr A_{1 , \i }  \e  + \mathscr A_{2 , \i }    \edag   =  \mathscr A_{1 , \,  \j }  \e  + \mathscr A_{2 , \,  \j }    \edag    .$$

Then its determinant is given by

$$  det A  =  det  \mathscr A_{1 , \i }   \, \e   +  det \mathscr A_{ 2 ,  \i }   \, \edag   =    det  \mathscr A_{1 , \,  \j }   \, \e   +  det \mathscr A_{ 2 ,  \,  \j }   \, \edag . $$

\end{subPn}
{\sf Proof:}  \\
The proof can be done by induction on $n.$ In the case $n=2$ this is immediately demonstrated
by the following easy calculation:

$$\begin{array}{rcl}
det A   &    =   &        \displaystyle     det  \left(  \begin{array}{rr}
                                                                   a_{11}    &  a_{12}  \\  &   \\
                                                                   a_{21}    &   a_{22}
                                                                   \end{array}  \right)
=  a_{11}  a_{22}  -  a_{12}  a_{21}
\\  &   &  \\   &  =  &
\displaystyle    \left(    a_{11  }'  \e  + a_{11}'' \edag \right)   \left(   a_{22}' \e  + a_{22}'' \edag   \right)
\\   &   &   \\   &   &
\qquad   \quad   -  \;      \left(    a_{12  }'  \e  + a_{12}'' \edag \right)   \left(   a_{21}' \e  + a_{21}'' \edag   \right)
\\   &   &   \\    &   =   &
\displaystyle   \left(   a_{11}'  a_{22}'  - a_{12}'  a_{21}' \right) \, \e  \;  + \;    \left(   a_{11}'' a_{22}''  -  a_{12}''  a_{21}''  \right)  \, \edag
\\   &   &   \\   &   =  &
det  \mathscr A_1  \e   + det  \mathscr A_2  \edag \, ,
\end{array} $$

where  the  idempotent  decompositions  may  be  taken  with    coefficients  either  in  $ \C ( \i )$  or  in    $ \C ( \j )$. For the general case, one can simply use the Laplace theorem, that gives the general formula for the determinant of an $n \times n$ matrix in terms of the determinants of suitable $(n-1) \times (n-1)$ matrices.
\mbox{}  \qed  \mbox{}   \\

An immediate consequence of this result, is the Binet theorem for bicomplex matrices.

\begin{subCy}
Let $A$, and $B$ be two square bicomplex matrices. Then
$$
det (AB) = det A \cdot det B.
$$
\end{subCy}
{\sf Proof:}  \\
By using the previous result, we see that
$$\begin{array}{rcl}
det(AB)  & = &   \displaystyle    det \left(  \left(  \mathscr A_{ 1 , \i} \,  \e  + \mathscr A_{ 2 , \i}    \,  \edag \right)   \left(  \mathscr B_{ 1 , \i}  \,  \e  + \mathscr B_{ 2 , \i}    \,   \edag \right)  \right)
\\  &  &  \\   &  =  &
\displaystyle  det   \left(   \mathscr A_{ 1 , \i}  \mathscr B_{ 1 , \i}   \,   \e  + \mathscr A_{ 2 , \i}  \mathscr B_{ 2 , \i}   \,    \edag \right)
\\   &  &  \\   & = &
\displaystyle  det  \left(   \mathscr A_{ 1 , \i}  \mathscr B_{ 1 , \i}  \right)   \,     \e  +   det  \left(   \mathscr A_{ 2 , \i}  \mathscr B_{ 2 , \i}  \right)  \,      \edag
\\   &  &  \\   & = &
\displaystyle  det    \mathscr A_{ 1 , \i} \cdot  det  \mathscr B_{ 1 , \i}   \,   \e  +   det    \mathscr A_{ 2 , \i}  \cdot  det    \mathscr B_{ 2 , \i}     \,   \edag
\\   &  &  \\   & = &
\displaystyle   \left(   det    \mathscr A_{ 1 , \i}   \,    \e  +    det  \mathscr A_{ 2 , \i}   \,    \edag  \right)    \left(      det    \mathscr B_{ 1 , \i}    \,   \e   +   det    \mathscr B_{ 2 , \i}    \,     \edag  \right)
\\  &  &  \\   &  =  &
det A  \cdot  det B  \, ,
\end{array}  $$

where  one can take   $ \mathscr A_{ 1 ,  \, \j} $,    $ \mathscr A_{ 2 ,  \, \j} $,       $ \mathscr B_{ 1 ,  \, \j} $,       $ \mathscr B_{ 2 ,  \, \j} $       instead  of  their  $ \C ( \i)$--coun\-ter\-parts. This concludes the proof.
\mbox{}  \qed  \mbox{}  \\

Analogously, one can use the idempotent representation of a bicomplex matrix to determine its invertibility.

\begin{subPn}  Let $   A  =  \mathscr A_{ 1 , \i}    \,   \e  + \mathscr A_{ 2 , \i}  \,     \edag   =   \mathscr A_{ 1 ,  \, \j }   \,    \e  +     \mathscr A_{ 2 ,  \, \j  }    \,    \edag    \in \bc^{n \times n }$,       $ \mathscr A_{ 1 , \i} $,   $ \mathscr A_{ 2 , \i}  \in \C ^{n \times  n} ( \i )$,      $ \mathscr A_{ 1 , \, \j} $,   $ \mathscr A_{ 2 , \, \j}  \in \C ^{n \times n }(\j)$ be a bicomplex matrix. Then $A$ is  invertible  if  and  only  if   $ \mathscr A_{ 1 , \i} $   and      $ \mathscr A_{ 2 , \i}  $    are  invertible  in   $ \C  ^{n \times n }(  \i)$   and     $ \mathscr A_{ 1 , \,  \j} $,      $ \mathscr A_{ 2 ,  \,   \j}  $    are  invertible  in   $ \C^{n \times n } ( \,  \j)$.
\end{subPn}
{\sf Proof} \\
A  matrix  $A$  is  invertible  if  and  only  if  there  exists   $ B  =  \mathscr B_{ 1 , \i}    \,   \e  +  \mathscr B_{ 2 , \i}    \,   \edag  \in  \bc^{n \times n}$  such  that     $  AB =  BA  = I$.  This  is  equivalent  to
$$  \begin{array}{rcl}
I_{ \C^{n\times n} ( \i ) }   \, \e  +  I_{  \C^{n \times n}  ( \i ) }  \,  \edag  &  =   &  I_{  \bc^{ n \times n} }   =  \;  \mathscr A_{ 1 , \i} \,     \mathscr B_{ 1 , \i}  \,   \e  + \mathscr A_{ 2 , \i}    \,   \mathscr B_{ 2 , \i}    \,    \edag
\\   &   &   \\   &   =  &
\mathscr B_{ 1 , \i}  \,    \mathscr A_{1 , \i }   \,  \e  + \mathscr B_{ 2 , \i}  \,    \mathscr A_{ 2 , \i}   \, \edag
\end{array}  $$

which  is  equivalent   to     $\mathscr A_{ 1 , \i}   \,    \mathscr B_{ 1 , \i}  = I_{ \C ^{n \times n} ( \i )}  $  and   $\mathscr A_{ 2 , \i}     \,   \mathscr B_{ 2 , \i}  = I_{ \C ^{n \times n} ( \i ) }  $.  The  same  when  one  considers  the  invertibility  of    $ \mathscr A_{ 1 , \,  \j} $  and        $ \mathscr A_{ 2 ,  \,   \j}  $.   \mbox{}\qed\mbox{}\\

The following result is immediate.

\begin{subCy}  A  matrix  $A \in \bc^{n \times n} $  is  invertible  if  and  only  if   $ det A \notin  \mathfrak S   \cup  \{ 0 \} $.
\end{subCy}

We can naturally introduce three conjugations on each bicomplex matrix  $\displaystyle    A  =   \left(  a_{ \ell \,  j  }     \right)  \in   \bc^{m \times n } $,  as follows:

$$  \displaystyle   A^\dagger :=    \left(  a^\dagger_{ \ell \,  j  }     \right) \, , \qquad       \overline{ A } :=    \left(  \overline{ a }_{ \ell \,  j  }     \right) \,   ,   \qquad     A^\ast :=    \left(  a^\ast_{ \ell \,  j  }     \right) \,   .$$

It  is  immediate to  prove  that   all  the conjugations are  multiplicative  ope\-ra\-tions,  that  is,  each  of  them  applied  to the  product  of  two  matrices  becomes  the  product  of  the  conjugate  matrices.  Obviously  these  conjugations  are   additive  operations  also.

\medskip

As  usual,   $A^t$  denotes  the  transposed  matrix:
$$  \displaystyle  A^t  =    \left(  a_{   j   \ell  }     \right)   ,    \quad   1 \leq j \leq n , \quad    1 \leq \ell \leq m  ,$$

and we correspondingly have three adjoint  matrices:
$$ \displaystyle  A^{ t \,  \dagger}  :=   \left(  A^t \right) ^\dagger  =  \left(   A^\dagger \right) ^t  =    \left(  a^\dagger_{   j   \ell  }     \right)  ;$$

$$ \displaystyle  A^{ t \,  bar}  :=   \overline{ \left(  A^t \right)  }  =  \left(  \,   \overline{ A }  \,  \right) ^t  =    \left(  \overline{ a }_{   j   \ell  }     \right)  ;$$

$$ \displaystyle  A^{ t \,  \ast}  :=   \left(  A^t \right) ^\ast  =  \left(   A^\ast  \right) ^t  =    \left(  a^\ast_{   j   \ell  }     \right)  .$$

Since   for any  two compatible  matrices one has  that   $ (AB)^t = B^t A^t$,  then    $(AB)^{t \, \dagger}  =   B^{ t \, \dagger}     A^{ t \, \dagger}   $,
 $(AB)^{t \,bar}  =   B^{ t \, bar}     A^{ t \, bar}   $,    $(AB)^{t \, \ast}  =   B^{ t \, \ast}     A^{ t \, \ast}   $.

 \medskip

 The  idempotent representations
$$  A  =  \mathscr A_{1 , \i }    \,  \e  + \mathscr A_{2 , \i }  \,   \edag   =  \mathscr A_{1 , \,  \j }  \,    \e  + \mathscr A_{2 , \,  \j }    \,    \edag    $$

give:
$$  A^\dagger  =     \mathscr A_{2 , \i }   \,   \e  + \mathscr A_{1 , \i }  \,  \edag   =  \mathscr A^\ast_{2 , \,  \j }   \,   \e  + \mathscr A^\ast_{1 , \,  \j }   \,   \edag  \, ,  $$

$$  \overline{ A }   =    \overline{  \mathscr A}_{2 , \i }   \,   \e  +   \overline{ \mathscr A}_{1 , \i }  \,  \edag   =  \mathscr A_{2 , \,  \j }   \,   \e  + \mathscr A_{1 , \,  \j }   \,   \edag  \, ,    $$

$$  A^\ast  =    \overline{  \mathscr A }_{1 , \i }   \,   \e  +   \overline{ \mathscr A}_{2 , \i }  \,  \edag   =  \mathscr A^\ast_{1 , \,  \j }   \,   \e  + \mathscr A^\ast_{2 , \,  \j }   \,   \edag  \, ,  $$

$$  A^t  =  \mathscr A^t_{1 , \i }   \,    \e  + \mathscr A^t_{2 , \i }  \,      \edag   =  \mathscr A^t_{1 , \,  \j }   \,    \e  + \mathscr A^t_{2 , \,  \j }    \,    \edag,    $$

and hence
$$  A^{t \, \dagger }  =     \mathscr A^t_{2 , \i }   \,   \e  + \mathscr A^t_{1 , \i }  \,  \edag   =  \mathscr A^{ t \, \ast}_{2 , \,  \j }   \,   \e  + \mathscr A^{  t \, \ast}_{1 , \,  \j }   \,   \edag  \, ,  $$

$$ A^{t \, bar }   =    \overline{  \mathscr A}^{ \; t }_{2 , \i }   \,   \e  +   \overline{ \mathscr A}^{  \;  t }_{1 , \i }  \,  \edag   =  \mathscr A^t_{2 , \,  \j }   \,   \e  + \mathscr A^t_{1 , \,  \j }   \,   \edag  \,  ,   $$

$$  A^{ t \, \ast }  =    \overline{  \mathscr A }^{ \; t }_{1 , \i }   \,   \e  +   \overline{ \mathscr A}^{ \; t }_{2 , \i }  \,  \edag   =  \mathscr A^{ t \, \ast}_{1 , \,  \j }   \,   \e  + \mathscr A^{ t \, \ast}_{2 , \,  \j }   \,   \edag  \, .  $$

Similarly, we can define the notion of self-adjoint matrix, with respect to each one of the conjugations introduced above, specifically we will say that a matrix $A$ is self adjoint if, respectively, it satisfies one of the following equalities.
$$  A =   A^{t \, \dagger } ,  \quad   A  =   A^{ t \, bar }  ,  \quad   A =  A^{ t \, \ast }  . $$

We can express  these  definitions    in  terms  of  idempotent  components.

\medskip

\begin{enumerate}

\item[(a)]  The matrix $A$  is  $\dagger$--self--adjoint,   or   $ \dagger$--Hermitian,    if  and only if
$$ \mathscr A_{1 , \i }  =  \mathscr A^t_{2 , \i },   \quad    {\rm  if  \;   and   \;    only   \;    if }   \quad      \mathscr A_{1 ,  \, \j }  =  \mathscr A^{ t \, \ast }_{2 , \,  \j },   $$

\noindent
which  is  true  if  and  only  if
$$A^t = A^\dagger;  $$

\noindent
 thus  all  $\dagger$--self--adjoint  matrices  are  of  the  form
$$  A  =  \mathscr A_{1 , \i }    \,  \e  + \mathscr A^t_{1 , \i }  \,   \edag   =  \mathscr A_{1 , \,  \j }  \,    \e  + \mathscr A^{ t \, \ast }_{1 , \,  \j }    \,    \edag    ,  $$

\noindent
with   $ \mathscr A_{ 1 , \i } $  an  arbitrary  matrix  in  $ \C^{n \times n } ( \i ) $,  and   $ \mathscr A_{ 1 , \,  \j }  $  any  matrix  in  $ \C^{n \times n } ( \j ) $.  Note  that   $ \mathscr A^{ t \, \ast }_{1 , \,  \j }  $  is   the  usual  $ \C ( \j ) $  adjoint  of   $ \mathscr A_{1 , \,  \j }  $.

\medskip

\item[(b)]  The matrix $A$  is  bar--self--adjoint,   or  bar--Hermitian,    if  and  only  if
$$ \mathscr A_{1 , \i }  =  \overline{ \mathscr A}^{ \; t }_{2 , \i },   \quad    {\rm  if  \;   and   \;    only   \;    if }   \quad      \mathscr A_{1 ,  \, \j }  =  \mathscr A^{ t }_{2 , \,  \j },   $$

\noindent
which  is  true  if  and  only  if
$$A^t = \overline{ A };  $$

\noindent
 thus  all  bar--self--adjoint  matrices  are  of  the  form
$$  A  =  \mathscr A_{1 , \i }    \,  \e  +  \overline{ \mathscr A }^{ \;  t }_{1 , \i }  \,   \edag   =  \mathscr A_{1 , \,  \j }  \,    \e  + \mathscr A^{ t  }_{1 , \,  \j }    \,    \edag    ,  $$

\noindent
with   $ \mathscr A_{ 1 , \i } $  an  arbitrary  matrix  in  $ \C^{n \times n } ( \i ) $,  and   $ \mathscr A_{ 1 , \,  \j }  $  any  matrix  in  $ \C^{n \times n } ( \j ) $.  Note  that   $ \overline{ \mathscr A}^{ \;  t  }_{1 ,   \i }  $  is   the  usual  $ \C ( \i ) $  adjoint  of   $ \mathscr A_{1 ,   \i }  $.

\medskip

\item[(c)]  The matrix $A$  is  $\ast$--self--adjoint,   or   $\ast$--Hermitian,    if  and  only  if
$$ \mathscr A_{1 , \i }  =  \overline{ \mathscr A}^{ \; t }_{1 , \i },   \quad     \mathscr A_{2 , \i }  =  \overline{ \mathscr A}^{ \; t }_{2 , \i }  ,  $$

\noindent
 if  and  only  if
$$ \mathscr A_{1 , \,  \j }  =  \mathscr A^{  t \, \ast  }_{1 , \,   \j },   \quad     \mathscr A_{2 , \,  \j }  =  \mathscr A^{ t \, \ast }_{2 , \i };  $$

\noindent
that is,  both     $  \mathscr A_{1 , \i }  $   and  $  \mathscr A_{2 , \i }  $   are  usual   $ \C ( \i ) $  self--adjoint matrices,   and   both     $  \mathscr A_{1 ,   \,  \j }  $,   $  \mathscr A_{2 , \,  \j }  $   are  usual   $ \C ( \j ) $  self--adjoint matrices.

\end{enumerate}

\medskip

\begin{subParr} In  what   follows  we  are   interested   in   $\ast$--self--adjointness  since this property implies  a form of hyperbolic   ``positiveness"  of   bicomplex   matrices  which we will find quite useful.
\end{subParr}

 \begin{subDn}
 A   $\ast$--self-adjoint   matrix $A\in\bc^{n\times n}$ is called  {\it   hyperbolic     positive},  if   for  every   column   $c\in\bc ^n$,
 \begin{equation}
 \label{c*Ac}
 c^{\ast \, t } \cdot  A   \cdot  c   \in \mathbb D^+.
 \end{equation}
 \end{subDn}

 In  this  case  we  write   $A \gtrdot 0$.  Given  two  bicomplex  matrices   $A$,  $B$,  we  say  that   $A \gtrdot  B$  if  and  only  if   $A - B \gtrdot 0$.

 \begin{subPn}
 \label{pnmat}
 Let
 \begin{equation}
 \label{eq:repA} A=A_1+\mathbf jA_2=\mathscr A_1\mathbf e+\mathscr
 A_2\mathbf e^\dag,
 \end{equation}
 be an element of $\bc ^n$, with $A_1,A_2,\mathscr A_1$ and $\mathscr A_2$ in
 $\mathbb C^{n\times n}  (\i) $. Then, the following are equivalent:

 \medskip

 \begin{enumerate}

 \item[ ({\it a})]    $A\gtrdot 0$.

\medskip

\item[({\it b})]    Both $\mathscr A_1$ and $\mathscr A_2$ are complex  positive
 matrices.

\medskip

\item[({\it c})]   $A_1\ge 0$, the matrix $A_2$ is skew self adjoint, that is,
 $A_2+   \overline{ A}_2^{\, t}   =0$,
 and
 \begin{equation}
 \label{ineq}
 -A_1\le \i A_2\le A_1.
 \end{equation}

 \end{enumerate}
\end{subPn}
{\sf Proof:} \\     Assume $(a)$. Take   a  column  $c\in\bc^n$ with representations
 \[
 c=c_1+\mathbf jc_2=\mathscr \zeta_1\mathbf e+\mathscr
 \zeta_2\mathbf e^\dag,
 \]
 where the various   columns    are in $\mathbb C^n (\mathbf i)$. Then,
 \begin{equation}
 \label{cAc} c ^{ \ast \, t } A  c   =     \overline{  \zeta}_1^{\, t}  \,   \mathscr A_1  \,  \zeta_1\mathbf e   +
  \overline{  \zeta}_2 ^{\, t }    \,   \mathscr A_2   \,   \zeta_2\mathbf e^\dag.
 \end{equation}

 Thus, by definition, $A\gtrdot   0$ implies that both $\mathscr A_1$ and
 $\mathscr A_2$ are   $\C ( \i ) $   positive, and $(b)$ holds.

 \medskip

 Assume ({\it b}). Since
 \[
 \mathscr A_1=A_1-\i A_2\quad{\rm and}\quad \mathscr A_2=A_1+\i
 A_2,
 \]

 then
\[
 A_1=\frac{1}{2}\left(\mathscr A_1+\mathscr A_2\right)\ge 0;
 \]

 moreover    $\i A_2=A_1-\mathscr A_1$ is self-adjoint, and hence $A_2$ is
 skew self-adjoint. Furthermore, still in view of $(b)$,
 \[
 A_1-\i A_2\ge 0\quad{\rm and}\quad A_1+\i A_2\ge 0,
 \]
 and thus we obtain \eqref{ineq} and $(c)$ holds. Finally when
 $(c)$ holds, both the matrices $\mathscr A_1$ and $\mathscr A_2$ are
 positive, and thus $(a)$ holds as well in view of \eqref{cAc}.
 \mbox{}\qed\mbox{}\\

 \begin{subPn}
 Let $A\in\bc^{n\times n}$. The following are equivalent:

 \begin{enumerate}

 \item[(1)]   $A$ is hyperbolic  positive.

 \medskip

 \item[(2)]    $A=B^{\ast \, t }   \cdot  B$ where $B\in\bc^{m\times n}$ for some $m\in\mathbb N$.

 \medskip

 \item[(3)] $A=C^2$ where the matrix $C$ is hyperbolic positive.

 \end{enumerate}
 \end{subPn}
 {\sf Proof:}  \\     Let $A\in\bc^{n\times n}$ be represented as in
 \eqref{eq:repA}, and assume that $(1)$ holds. Then, by the preceding
 theorem, we have $\mathscr A_1\ge 0$ and $\mathscr A_2\ge0$, and thus
 we can write
 \begin{equation}
 \label{eqA1A2}
 \mathscr A_1= \overline{ U}^{\, t }U\quad{\rm  and}\quad \mathscr A_2=   \overline{V}^{\, t }    V,
 \end{equation}

 where $U$ and $V$ are matrices in $\mathbb C^{n\times n} ( \i ) $. Thus
 $A=B^{\ast \, t }   B$ with $B :  =U\mathbf e+V\mathbf e^\dag$, so that (2) holds
 with $m=n$. Assume now that (2) holds with $B\in\bc^{m   \times   n}$     for some $m\in\mathbb N$. Writing $B=U\mathbf e+V\mathbf
 e^\dag$, where now $U$ and $V$ both belong to $\mathbb C^{m\times
 n}  ( \i )  $ we have
 \[
 A=\overline{ U }^{\, t }    U\mathbf e    +     \overline{V}^{\, t}    V   \mathbf e^\dag,
 \]
 and so $(1)$ holds. The equivalence with $(3)$ stems from the fact
 that $U$ and $V$ in \eqref{eqA1A2} can be chosen positive.
 \mbox{}\qed\mbox{}\\

In the case of a complex matrix, it is equivalent to say that a
 positive
 matrix $A$ is Hermitian and to say that its eigenvalues are positive.
 We now give the corresponding result in the setting of $\bc$. We note that
 the existence of   zero   divisors   creates problems  in  an  attempt   to classify
 eigenvalues in general. For instance,     take two non-zero elements $a$
 and $b$ such that $ab=0$. Then every $c\in\bc$ such that $bc=0$
 is an eigenvalue of the matrix
\[
 \begin{pmatrix}a&a\\a&a\end{pmatrix}\quad
 \mbox{\rm
 with eigenvector}\quad
 \begin{pmatrix}b\\b\end{pmatrix},
 \]
 since
 \[
 \begin{pmatrix}a&a\\a&a\end{pmatrix}\begin{pmatrix}b\\b\end{pmatrix}=
 c\begin{pmatrix}b\\b\end{pmatrix}=
 \begin{pmatrix}0\\0\end{pmatrix}.
 \]

There  is  however     a  relation  between  the  eigenvalues  and  eigenvectors  of  a  bicomplex  matrix $A$  and  those   of  its   idempotent  components  $\mathscr A_1$ and  $\mathscr A_2$.  Indeed,  let  $ \lambda = \gamma_1 \, \e  + \gamma_2 \, \edag    \in  \bc  \setminus  \{  0 \} $    be  an  eigenvalue  of  $A$  with  a  corresponding  eigenvector   $u = v_1 \, \e  + v_2 \, \edag$,  then
$$ A \, u = \lambda \, u ,$$

which is equivalent to
$$ \displaystyle \left\{  \begin{array}{rcl}
\mathscr A_1 \, v_1 & = &  \gamma_1 \, v_1  ,
\\  &  &   \\
\mathscr A_2 \, v_2 & = &  \gamma_2 \, v_2  .
\end{array}  \right.  $$

If    $ \lambda   $  is  not  a  zero  divisor   and  $v_1 \neq 0$,   $v_2 \neq 0$    then  $\lambda$  is  an  eigenvalue of   $A$  if  and  only  if  $\gamma_1$  is  an  eigenvalue  of  $\mathscr A_1$ and   $\gamma_2$  is  an  eigenvalue  of  $\mathscr A_2$.

\begin{subtheorem}
A matrix  $A \in \bc^{n \times n}$ is hyperbolic    positive  if and only if

\begin{enumerate}

\item  $A$  is  $\ast$--Hermitian;

\item  none of its eigenvalue is a zero divisor  in  $\mathbb D^+$.
\end{enumerate}

\end{subtheorem}
 {\sf Proof:}  \\      That  $A$ is   hyperbolic   positive  is equivalent to   the  complex  matrices     $ \mathscr A_1$,  $\mathscr A_2$ being    positive,  which is in turn equivalent to state that  $\mathscr A_1$ and $\mathscr A_2$  are   complex   Hermitian  and  any of  their   eigenvalues  is  a  positive  number.  Finally  this is   equivalent to  say  that  $A$  is  $\ast$--Hermitian  and  that none   of its  eigenvalues   is   a zero divisor  in  $\mathbb D^+$.
 \mbox{}\qed\mbox{}\\

\begin{subCy}
A matrix $A \in \bc^{n \times n}$ is hyperbolic    positive  if and only if:

\begin{enumerate}

\item  $A$  is  $\ast$--Hermitian;

\item  if $\lambda =  \lambda_1 + \j \, \lambda_2$  is an eigenvalue for $A$, then $\lambda_1 >0$,   $ \i  \, \lambda_2 \in \mathbb R$  and
$$   - \lambda_1  <  \i \, \lambda_2 <  \lambda_1 \,   ;  $$

\end{enumerate}

\end{subCy}
 {\sf Proof:}  \\      This is  because   $ \lambda  =   \lambda_1  + \j \, \lambda_2 =  \gamma_1 \, \e  +  \gamma_2 \, \edag $  with   $\gamma_1 =  \lambda_1  - \i \,\lambda_2 $ and   $\gamma_2 =  \lambda_1  + \i \,\lambda_2 $.
 \mbox{}\qed\mbox{}\\

\begin{subremark}
The  last  inequality  in the statement of this   corollary is  equivalent to the  system
$$  \displaystyle     \left\{  \begin{array}{l}
\lambda_1 - \i \, \lambda_2 > 0  ,
\\   \\
\lambda_1 + \i \, \lambda_2 > 0 \, .
\end{array}  \right.  $$

\end{subremark}

\medskip

It  is  known that  in  the  case  of  complex   matrices  every    eigenvector  corresponds  to  only  one  eigenvalue. This is not the case for bicomplex matrices.    Specifically,     a  bicomplex  eigenvector  can  correspond  to  an infinite  family  of   bicomplex   eigenvalues.

\medskip

We  restrict   our considerations  to  hyperbolic     positive  matrices.   Let   $A$   be  such  a   matrix,  and let   $\lambda   =  \gamma_1 \, \e + \gamma_2 \, \edag$  be  one of its eigenvalues (in particular $\lambda$  is  a  non  zero  divisor  in  $\mathbb D^+$).  First   of  all    let us  show
that  any  such  egenvalue  has  an  eigenvector of the  form  $u =  v_2 \, \edag$  with   $ v_2 \in  \C^n (\i)$.   Since   $\gamma_2$  is  an  eigenvalue  of   $\mathscr A_2$,      let   $v_2$  be   a   corresponding    eigenvector:   $\mathscr A_2 \,v_2 = \gamma_2 \, v_2$.  Consider  $u : =  v_2 \, \edag$.  Let  us  show  that  it  is  an  eigenvector  of  $A$  corresponding  to  the  above   $\lambda$.    Indeed,  $A \, u  =   \mathscr A_2 \, v_2 \, \edag   =   \gamma_2 \, v_2  \, \edag     $  and   $ \lambda \, u =  \gamma_2 \, v_2 \, \edag$,  thus   $A \, u =  \lambda \, u$.

\medskip

Now we are in a    position    to  show  that  this  eigenvector   corresponds  to  an  infinite  family  of   eigenvalues.   For  any  $r>0 $  set    $ \lambda_r  :=  r \, \e  +  \gamma_2 \, \edag$,  then   $ \lambda_r \, u   =   \gamma_2 \, v_2 \, \edag =  A \, u$.  Hence  the  whole  family   $ \{  \,  \lambda_r \mid  r   > 0 \, \}   $    consists   of  the   eigenvalues   of  $A$    with  the  same  eigenvector    $u  \in   \bc_\edag^n$.   \\  

  We now study  bicomplex  $\ast$--unitary matrices, that is,  matrices $U\in\bc^{n\times n}$
 such that $UU^{\ast \, t }=U^{\ast \, t }  U=I_n$.

 \begin{subPn}
 \label{pn:real} Let $U=U_1+\mathbf j U_2    =  \mathscr U_1 \, e + \mathscr U_2 \, \edag     \in\bc^{n\times n}$. Then
 $U$ is unitary if and only if   its  idempotent  components  are  complex  unitary  matrices   or,   equivalently,  its  cartesian  components satisfy
$$  \displaystyle    U_1 \, \overline{ U}_1^{\, t}  +  U_2 \, \overline{ U}_2^{\, t}    =       \overline{ U}_1^{\, t}  \,   U_1  +  \overline{ U}_2^{\, t}   \,  U_2   =   I_n  $$

  and
$$    U_2   \,  \overline{ U}_1^{\, t}  =     U_1   \overline{ U}_2^{\, t}    ,   \qquad     \overline{ U}_1^{\, t}  \,   U_2   =   \overline{ U}_2^{\, t}  \,    U_1  .   $$

 \end{subPn}
 {\sf Proof:}   \\       Since   $U^{\ast \, t }   =    \overline{U}^{\, t}_1   -         \mathbf j \,   \overline{U}^{\, t}_2    = \overline{ \mathscr U}^{\, t}_1 \, e +   \overline{ \mathscr U}^{\, t}_2 \, \edag  $,  then  for   the   idempotent  representation  we  have:
$$ \mathscr U_1   \,   \overline{ \mathscr U}_1^{\, t}   \, \e  +    \mathscr U_2 \overline{ \mathscr U}_2^{\, t} \, \edag  =   \e \, I_n  + \edag \, I_n $$

 and
$$ \overline{ \mathscr U}_1^{\, t}    \,    \mathscr U_1  \, \e  +   \overline{ \mathscr U}_2^{\, t} \,   \mathscr U_2  \, \edag  =   \e \, I_n  + \edag \, I_n  .$$

This  means   that  $  \mathscr U_1 $  and     $ \mathscr U_2  $  are  complex   unitary   matrices.

 \medskip

 Similarly  for the cartesian representation we have:
$$\displaystyle    U_1 \, \overline{ U}_1^{\, t}  +  U_2 \, \overline{ U}_2^{\, t}    +   \j  \,  \left( U_2   \,     \overline{ U}_1^{\, t}      -   U_1   \,    \overline{ U}_2^{\, t}       \right)    = I_n  ,$$

$$\displaystyle    \overline{ U}_1^{\, t}  \,   U_1  +  \overline{ U}_2^{\, t}   \,  U_2    +   \j  \,  \left(   \overline{ U}_1^{\, t}  \,   U_2   -   \overline{ U}_2^{\, t}  \,    U_1  \right)    = I_n  ,$$

 which means  that
 $$  \displaystyle    U_1 \, \overline{ U}_1^{\, t}  +  U_2 \, \overline{ U}_2^{\, t}    =       \overline{ U}_1^{\, t}  \,   U_1  +  \overline{ U}_2^{\, t}   \,  U_2   =   I_n  $$
 and
 $$    U_2   \,  \overline{ U}_1^{\, t}  =     U_1   \overline{ U}_2^{\, t}    ,   \qquad     \overline{ U}_1^{\, t}  \,   U_2   =   \overline{ U}_2^{\, t}  \,    U_1  .   $$

  The result follows.
 \mbox{}\qed\mbox{}\\

\bigskip


\chapter{$\bc$--modules }


\section{$\mathbb B\mathbb C$--modules  and  involutions on  them }\label{BC modules}

\begin{subParr}  Let  $X$  be  a    $\mathbb B  \mathbb C$-module.     It  turns  out  that   some  structural  peculiarities   of   the  set   $\bc$  are immediately manifested  in  any  bicomplex   module, in contrast  to  the  cases  of   real,   complex, or  even  quaternionic  linear  spaces  where the structure of linear space does not imply anything immediate about the space itself.

\medskip

Indeed, consider  the  sets
$$  X_{ \mathbf e}:=  \mathbf e \cdot X  \qquad  {\rm and}  \qquad    X_{ \mathbf e^\dagger}:=  \mathbf e^\dagger \cdot X . $$
\end{subParr}

Since
$$  X_{ \mathbf e}   \cap   X_{ \mathbf e^\dagger}  =  \left\{ 0 \right\} , $$

and
\begin{equation}\label{descomposicion idempotente modulo bicomplejo}
 X  =     X_{ \mathbf e}  +  X_{ \mathbf e^\dagger}   \, ,
\end{equation}

one can define two  mappings
$$  P:  X \to  X   \, , \quad  Q: X \to X  $$

by
$$  P(x)  :=  \mathbf e \, x  \, \quad  \quad     Q(x)  :=  \mathbf e^\dagger \, x
\, . $$

Since
$$  P + Q = Id_X,  \quad       P \circ Q  =  Q  \circ  P  =  0,  \quad  P^2  =P  \quad    { \rm  and}   \quad  Q^2 =Q   , $$

then  the  operators  $P$  and  $Q$  are  mutually  complementary  projectors.   Formula    (\ref{descomposicion idempotente modulo bicomplejo})     is  called  the  idempotent  decomposition  of  $X$, and it  plays  an  extremely  important  role  in what  follows.   In  particular,  it  allows  to  realize  component--wise  the  operations  on  $X$:    if   $x =  \e x + \edag x$,  $y=  \e y + \edag y$  and if  $\lambda = \lambda_1 \e + \lambda_2 \edag$  then  $x+ y =  (\e x + \e y ) + ( \edag x + \edag y)$,  $\lambda x =  \lambda_1 x \e + \lambda_2 x \edag$.

\medskip

In  what  follows  we  will   write  $X_{ \mathbb C (\mathbf i ) } $  or   $X_{ \mathbb C (\mathbf j ) } $  whenever   $X$  is  considered  as  a   $\mathbb C ( \mathbf i )$  or   $\mathbb C ( \mathbf j )$   linear  space  respectively.

\medskip

Since     $X_{\mathbf e}$  and  $  X_{\mathbf e^\dagger}$  are  $\mathbb R$-,  $\mathbb C(\mathbf i)$--   and   $\mathbb C(\mathbf j)$--linear  spaces
as well as    $\mathbb B  \mathbb C$-modules, we have that $ X  =  X_{\mathbf e}
\oplus   X_{\mathbf e^\dagger} $  where   the  direct  sum    $\oplus$ can be  understood  in  the  sense
of   $\mathbb R$-,  $\mathbb C(\mathbf i)$--   or   $\mathbb C(\mathbf j)$--linear  spaces, as well as  $\mathbb B  \mathbb
C$-modules.

\medskip

Note  that   given  $ w \in  X_{\mathbf e}$    there  exists  $ x \in X$
such  that  $ w =  P(x)  =  \mathbf e \, x $  which  implies  $
\mathbf e \,  w   =  \mathbf e^2 \, x  =     \mathbf e \, x  = w
$,  i.e.,  any    $w \in  X_{\mathbf e} $  is  such  that    $w =  \mathbf
e \, w$.   Similarly  if $t \in    X_{\mathbf e^\dagger} $, then $  t =  \mathbf
e^\dagger \, t$.  As  a   matter  of  fact these identities characterize the  elements  in  $X_\e$  and  $X_\edag$:   $w \in  X_\e$  if  and  only  if   $w  =  \e \cdot w$  and    $t \in X_\edag$  if  and  only  if  \,  $t=  \edag \cdot  t$.  \\

In  order to be able  to say more  about  $ X_{\mathbf e}   $  and  $X_{\mathbf e^\dagger} $  we  need  additional  assumptions  about  $X$.  First of all, we are interested   in    the  cartesian--like  decompositions of   $\bc$--modules connected with the  analogs of the  bicomplex conjugations.

\begin{subDn} An involution  on a $\bc $-module  $X$ is   a  map
$$ \mathcal In : X \to X$$

such that

\begin{enumerate}

\item[{\it 1})] $ \mathcal In (x + y ) =  \mathcal In (x)  +  \mathcal In  (y)  $  for  any  $x, \, y \in X$;

\medskip

\item[{\it 2})]     $ \mathcal In^2  =  I_X \, $,  the  identity map;

\medskip

\noindent
and  one of the following conditions  holds  for all  $ Z  \in   \mathbb B
\mathbb C$   and al $x  \in  X$:

\medskip

\item[{\it 3a})]  $ \mathcal In (Z \, x) =\overline Z \, \mathcal In_{bar} (x)$; in this case we will say this is a bar-involution, and we will indicate it by the symbol $\mathcal In_{bar}$; in particular we note that such an involution is  $\mathbb C ( \mathbf j ) $-linear  while  it  is   $\mathbb C ( \mathbf i ) $-anti-linear and  $\mathbb D$-anti-linear;

\medskip

\item[{\it 3b})]   $ \mathcal In(Z \, x) =   Z^\dagger \, \mathcal In_\dagger (x)$; in this case we will say this is a $\dagger$-involution, and we will indicate it by the symbol $ \mathcal In_\dagger$; we note that such an involution is  $\mathbb C ( \mathbf i ) $-linear  while  it  is   $\mathbb C ( \mathbf j ) $--    and  $\mathbb D$-anti-linear;

\medskip

\item[{\it 3c})]     $ \mathcal In (Z \, x )     =      Z^\ast \, \mathcal In_\ast (x)$; in this case we will say this is a $\ast$-involution, and we will indicate it by the symbol $\mathcal In_\ast$; we note that such an involution is  $\mathbb D  $-linear  while  it  is   $\mathbb C ( \mathbf j ) $--   and  $\mathbb C  ( \mathbf i )$-anti-linear.

\end{enumerate}
\end{subDn}
\medskip

Each of these involutions induces a corresponding decomposition in $X$. Let us briefly describe them, beginning with $\mathcal In_\dagger.$
  Define
$$  \displaystyle    X_{1 , \dagger}  :=  \left\{  x \in X \mid  \mathcal In_\dagger \,
x  = x \, \right\}  \quad  {\rm and}  \quad     X_{2 , \dagger}   :=   \left\{
x \in X \mid  \mathcal In_\dagger  \, x  = -  x \, \right\}  \, .$$
We then have the following result:

\begin{subPn}
With the   above   notations. it   holds:
\begin{equation}\label{desc inv dagger}
X  =  X_{1 , \dagger}   + X_{2 , \dagger}  =  X_{1 , \dagger}   + \mathbf j \, X_{1 , \dagger}  \, .
\end{equation}

\end{subPn}

{\sf Proof}  \\  We first show that $\mathbf j X_{1, \dagger} = X_{2, \dagger}.$ Indeed if
$ x \in  X_{1 , \dagger} $,  then      $  \mathcal In_\dagger  ( \mathbf j \, x )  =
- \mathbf j \,  \mathcal In_\dagger (x)  =   - \mathbf j \, x$, that is,
$ \mathbf j \, x \in  X_{2 , \dagger} $, which means   that     $ \mathbf j \, X_{1 , \dagger}
\subset X_{2 , \dagger} $.   Conversely if  $ w \in X_{2 , \dagger} $,   then       $  \mathcal In_\dagger  ( -
\mathbf j \, w )  =   \mathbf j \,  \mathcal In_\dagger (w)  =   -
\mathbf j \, w$. This implies  $- \mathbf j \,  w = x \in X_{1 , \dagger} $  or,
equivalently,     $ w = \mathbf j \, x $,  thus
$$  \mathbf j \,  X_{1 , \dagger}
= X_{2 , \dagger}   . $$

Moreover, given  $x \in X$,  define
$$ \displaystyle  x_1  :=  \frac{1}{2}
\left( x + \mathcal I n_\dagger \, x \right)  \quad  {\rm and} \quad
x_2  :=  \frac{1}{2} \left( x - \mathcal I n_\dagger \, x \right)  . $$

 Clearly
$x_1  \in X_{1 , \dagger} $,  $x_2 \in X_{2 , \dagger} $,  and  $x =  x_1 + x_2  =  x_1  +
\mathbf j \, ( - \mathbf j \, x_2)$,  proving  the claim.
\mbox{}  \qed \mbox{}  \\

 We will say that $X_{1 , \dagger} $ is the  set  of  $\mathbb C ( \mathbf i )$--elements  in  $X$  since  it  is  of  course   a $\mathbb C ( \mathbf i )$--linear  subspace  of  $X$.    The representation   (\ref{desc inv dagger}) is a generalization  of the cartesian  representation  of     $\mathbb B
\mathbb C$.

\medskip

Note that   any $\dagger$--involution on a $\bc$-module $X$ is  consistent  with its
idempotent  re\-pre\-sen\-ta\-tion. Specifically, if $x \in X$,  $  x =  x_1
+ \mathbf j \, y_1$,    $x_1, \, y_1  \in X_{1 , \dagger} $    we
have  that
\begin{equation}\label{ecuacion BB}
x =  \mathbf e \, ( x_1 - \mathbf i \, y_1 )   +   \mathbf e^\dagger \, ( x_1 + \mathbf i \, y_1 )  \, .
\end{equation}

As we know the idempotent representation of a $\mathbb B
\mathbb C$--module exists  always,  but  if the  $\mathbb B
\mathbb C$--module  has also  a    $\dagger$--involution   which  generates a cartesian  decomposition,    then  both representations are related in the same way as it happens in $\mathbb B
\mathbb C$.     Note  that  both  elements   $  x_1 - \mathbf i \, y_1 $   and  $ x_1 + \mathbf i \, y_1 $  are   $\mathbb C (\mathbf i )$--elements,  and  in  general we  cannot  say  more.

\medskip

Exactly the same analysis can be made for $\mathcal In_{bar}$, and $\mathcal In_{\ast}.$ If  we  set
$$ X_{1, bar} := \left\{   \, x \in X \mid  \mathcal In_{bar} (x) = x \,    \right\}  $$

(the set of   $\mathbb C ( \mathbf j )$--complex  elements),  and
$$ X_{2, bar} := \left\{   \, x \in X \mid  \mathcal In_{bar} (x) =  - x \,    \right\}  \, ,  $$

then,  again,  $X_{2, bar}  = \mathbf i \, X_{1, bar}$  and
\begin{equation}\label{ecuacion AA}
X =  X_{1, bar}  + X_{2, bar}  = X_{1, bar}  + \mathbf i \, X_{1, bar} \, .
\end{equation}

This  decomposition  is  also  consistent  with  the  idempotent  one,    since  for     $x =  u_1  + \mathbf i \, u_2 \in X$  with  $u_1, \, u_2 \in     X_{1, bar}$  there  holds:
\begin{equation}\label{ecuacion CC}
x =  \mathbf e \, ( u_1 - \mathbf j \, u_2 )   +   \mathbf e^\dagger \, ( u_1 + \mathbf j \, u_2 )  \, .
\end{equation}

If  $X$  has  the  two  involutions   $ \mathcal In_{\dagger} $  and  $ \mathcal In_{bar}$,  then  both  decompositions   (\ref{ecuacion BB})  and  (\ref{ecuacion CC})  hold  for  any  $x \in X$;  moreover   
$$ \mathbf e \, ( u_1 - \mathbf j \, u_2 )   =   \mathbf e \, ( x_1 - \mathbf i \, y_1 ) $$

and
$$    \mathbf e^\dagger \, ( u_1 + \mathbf j \, u_2 )   =     \mathbf e^\dagger \,    ( x_1 + \mathbf i \, y_1 ) $$

although  not  necessarily      $  u_1 - \mathbf j \, u_2  =   x_1 - \mathbf i \, y_1 $   and      $  u_1 + \mathbf j \, u_2  =   x_1 + \mathbf i \, y_1 $.

\medskip

Finally,  for  $\mathcal In_\ast $ we have:
$$ X_{1, \ast} := \left\{   \, x \in X \mid  \mathcal In_\ast (x) = x \,    \right\}  $$

(the set of   hyperbolic elements),  and
$$ X_{2, \ast}  := \left\{   \, x \in X \mid  \mathcal In_\ast (x) =  - x \,    \right\}  \, ,  $$

with    $X_{2, \ast }  = \mathbf i \, \mathbf j  \, X_{1, \ast}  $  and
$$  X =  X_{1, \ast}  + X_{2 , \ast}  = X_{1, \ast}  + \mathbf i     \, \mathbf j    \, X_{1, \ast} \, .$$

The  existence  of  the  idempotent
representation  does  not necessarily imply  the  existence  of  an
involution.  Note  first  that  a  $\dagger$--involution  exists  if  and
only if  the  module  $X$  has  a  decomposition  (\ref{desc inv dagger}),  hence,  having   only   the      hypothesis  of  the
idempotent  representation  and  trying  to  look  for  components
$x_1 , \, x_2$ such  that  $x=  x_1 + \mathbf j \, x_2$,  one  is  led    to the  system
\begin{equation}\label{tres estrellas}
\begin{array}{rcl}
2 \, \mathbf e \, x & = &  ( 1 + \mathbf i \, \mathbf j ) \, x_1
+    ( \, \mathbf j - \mathbf i \, ) \, x_2 \, ;  \\  &  &  \\
2 \, \mathbf e^\dagger \, x & = &  ( 1 -  \mathbf i \, \mathbf j )
\, x_1  +    ( \, \mathbf j + \mathbf i \, ) \, x_2 \,  .
\end{array}
\end{equation}

Since  $\mathbf e $ and $\mathbf e^\dagger$  are  zero  divisors
there  is  no  guarantee  that  the  system  (\ref{tres estrellas})
has  a  solution.

\bigskip


\section{Constructing a a  $\bc$--module
from two  complex  linear  spaces.}
\label{construction}

Let   $X_1 , \, X_2$  be  $\C ( \i)$--linear  spaces.  We want  to  construct  a   $\bc$--module  $X$  such  that,
with suitable meaning of the symbols,
\begin{equation}\label{eqA}
X   =  \e X_1  + \edag X_2 .
\end{equation}

This means that we have to give a precise  meaning to  the  symbols  $ \e X_1  $  and   $  \edag X_2  $. To this purpose, consider  the  $\C ( \i ) $--linear  space   $ \e \,  \C ( \i)$  which  is  a   $ \C ( \i)$--linear  subspace  of  $\bc$, and define the tensor products
$$   \e X_1  :=  \e \, \C ( \i ) \otimes_{\C ( \i ) } X_1 \, ,$$

and
$$   \edag X_2  :=  \edag \, \C ( \i ) \otimes_{\C ( \i ) } X_2 \, .$$

It is clear that both  $   \e X_1 $   and   $ \edag  X_2 $  are  $\C (\i)$--linear  spaces.

\medskip

From  the  definition  of  tensor   product  any  elementary   tensor  in  $ \e X_1$    is  of  the   form     $ \e  \lambda  \otimes  x  $  with     $\lambda  \in     \C (\i)  $,  hence
$$    \e  \lambda  \otimes  x   =       \e    \otimes  \lambda  x   =      \e   \otimes   x_1  .  $$

That is, any element in  $\e X_1 $   is of the form  $ \e \otimes x_1$,  with    $x_1   \in X_1$, which we will write   as    $\e x_1$.  Similarly  with    $ \edag X_2$.

\medskip

Consider  now  the  cartesian  product    $\e X_1  \times  \edag X_2$  where,  as  usual,    $ \e X_1$  is  seen  as   $ \e X_1  \times  \{ 0 \} $  and   $ \edag  X_2   $  is  seen  as  $ \{ 0 \}  \times  \edag X_2$.    Since   $ \e X_1  \times   \edag  X_2$  is   an  additive  abelian   group  we   have    endowed     the  sum    $ \e X_1  +  \edag  X_2$  with  a  meaning:     for  any   $x_1 \in X_1$,   $x_2 \in  X_2$
$$    \e x_1  + \edag x_2   :=   (  \e x_1   ,  \edag  x_2  )   \in    \e X_1  \times   \edag  X_2    .$$

Hence  the  right--hand  side  in   (\ref{eqA})  is   defined  already  but  for  the  moment  as  a   $\C ( \i )$--linear  space  only.   Let's  endow   it    with  the  structure  of  $\bc$--module.      Given   $ \lambda  =   \beta_1 \e  + \beta_2  \edag   \in \bc$  and
$  \e x_1  + \edag  x_2    \in    \e  X_1  +   \edag  X_2$,  we  set:
$$  \lambda  (   \e x_1  + \edag x_2  )   =   ( \beta_1 \e  + \beta_2  \edag  )   \cdot  (     \e x_1  + \edag x_2   )   :=   \e  (  \beta_1 x_1 ) + \edag  ( \beta_2 x_2 ) .$$

It  is  immediate to  check  that  this is a  well  defined  multiplication  of  the  elements   of       $  \e  X_1  +   \edag  X_2$  by  bicomplex  scalars.  Thus  formula  (\ref{eqA})  defines  indeed  a   $\bc$--module.

\medskip

Of  course, in  complete  accordance  with    Section    \ref{BC modules}, one has
$$  X_\e :=  \e X  =  \e  X_1  ;$$
$$  X_ \edag   :=   \edag X  =   \edag X_2 , $$

While $X_1$  and  $X_2$  are  initially   $\C ( \i)$--linear  spaces  only, the new sets  $ \e X_1$  and   $\edag X_2$  are $\C ( \j)$--linear  spaces as well, and therefore  they  become   $\bc$--modules. This is a consequence of the fact that since  $ \j \e  =   - \i \e $  and   $ \j \edag  =   \i \edag $,  we can set

$$  \j  (\e  x_1 ):=  \e  (  - \i  x_1  )  \in  \e  X_1  \, , \quad  \forall  x_1 \in   X_1 \, ; $$

$$  \j  (\edag  x_2 ):=  \edag  (  \i  x_2  )  \in  \edag  X_2  \, , \quad  \forall  x_2 \in   X_2 \, . $$

  Notice  that $ X_1 $   and     $X_2$ are not, in general the  ``cartesian  components"  of  a   $\bc$--module.  Indeed,  as  before  we  can  define  $ \j X_2 :=  \j \C ( \i )  \otimes_{ \C ( \i ) } X_2$,  and   thus we may consider now the  additive   abelian  group   $X_1  + \j X_2$.  Thus  we  are  led  to  define  the  multiplication  by  bicomplex  scalars  by  the     formula
$$ \begin{array}{rcl}
 (z_1 + z_2 \j ) \cdot (x_1 + \j x_2)     &   :=   &     (z_1 x_1 - z_2 x_2 ) + \j ( z_1 x_2 + z_2 x_1 )
 \\   &   &   \\   &   =  &
   ( z_1 x_1 - z_2 x_2 \, , \, \j ( z_1 x_2 + z_2 x_1 ) )
\end{array}   $$

but  the symbols  $  z_1 x_1 - z_2 x_2   $   and  $  z_1 x_2 + z_2 x_1  $  are  meaningless.

\medskip

We could have obviously repeated the exact same procedure beginning with two  $\C (\j)$--linear  spaces, and construcing a new  $\bc$--module from them. This process, however, would not work if one starts with  two  hyperbolic  modules.


\chapter{Norms and  inner  products  on  $\bc$--modules}


%

\section{Real--valued norms on bicomplex modules}

In this section we consider a norm  on  a  bicomplex  module which  extends the usual  properties of the
Euclidean  norm  on  $\bc$.  Another  approach  which  generalizes the notion  of  $\D$-valued  norm  on  $\bc$  will be  considered    in  the  following  section.

\begin{subDn}\label{Def bicomplex norm}
 {\it  Let $X$  be  a   $\mathbb B \mathbb C$-module and   let  $ \| \cdot \|$  be  a  norm  on  $X$  seen  as  a  real  linear  space. We say that  $ \| \cdot \|$ is a (real--valued)    norm    on the   $\mathbb B  \mathbb C$--module  $X$ if  for any $ \mu \in  \bc$ it is }
$$  \| \mu \, x \|  \leq  \sqrt{ 2} \, | \mu | \cdot \| x \| \, .$$
\end{subDn}

\medskip

This definition  extends  obviously    property    (\ref{desigualdad norma euclideana})  but it  does  not  enjoy    necessarily  the  other  properties  of the  Euclidean  norm.   This   leads us to more  definitions.

\begin{subDn}
{\it  A     norm   $ \| \cdot \|$   on  a  $\mathbb B  \mathbb C$-module  $X$  is  called  a  norm  of  $ \mathbb C ( \mathbf i )$--type   if }
$$  \| \lambda \, x \| = | \lambda | \cdot \| x \| \qquad  {\rm for  \, all}  \;    \lambda \in \mathbb C ( \mathbf i ) \, ,  \; \;  {\rm for \, all}  \;  x \in X \, .$$
\end{subDn}

\begin{subDn}
{\it  A  norm   $ \| \cdot \|$   on  a  $\mathbb B  \mathbb C$-module  $X$  is  called  a  norm  of   $ \mathbb C ( \mathbf j )$--type   if }
$$  \| \gamma \, x \| = | \gamma | \cdot \| x \| \qquad  {\rm for  \, all}  \;    \gamma \in \mathbb C ( \mathbf j ) \, ,  \; \;  {\rm  for \, all}  \;   x \in X \, .$$
\end{subDn}

\begin{subDn}
{\it  A  norm   $ \| \cdot \|$   on  a  $\mathbb B  \mathbb C$-module  $X$  is  called  a  norm  of   complex   type  if  it is  simultaneously  of  $ \mathbb C ( \mathbf i )$--type  and   of  $ \mathbb C ( \mathbf j )$--type. }
\end{subDn}

\begin{subEx} {\rm   The  Euclidean  norm on    $\mathbb B  \mathbb C$  is  a      norm  of   complex    type.  Other  obvious  examples  are  the usual norms on $\mathbb B  \mathbb C^n$  and on $L_p ( \Omega, \, \mathbb B  \mathbb C) $  with  $p >1$.   }
\end{subEx}

\medskip


\begin{subParr} \label{norm construction}  Let us take up again the situation of Section  \ref{construction},  that  is, we  are   given  two  $\C ( \i )$--linear  spaces  $X_1$, $X_2$,  and   $X =   \e X_1  + \edag X_2$.   Now  assume  additionally  that  $X_1$  and  $X_2$   are normed  spaces  with  respective   norms    $\| \cdot \|_1$,     $\| \cdot \|_2$.  For  any   $x = \e x_1 + \edag x_2 \in X$, set
\end{subParr}
\begin{equation}\label{eqN}
\displaystyle   \| x \|_X :=   \frac{1}{ \sqrt{2} }    \sqrt{    \| x_1 \|_1^2   +   \| x_2 \|_2^2  \,   . }
\end{equation}

We can  show  that     (\ref{eqN})   defines  a   norm  on  $X$  in  the  sense  of  Definition   \ref{Def bicomplex norm}.   We  will  refer  to it  as  the  Euclidean--type       norm  on   $X =  \e X_1  + \edag  X_2$.    It  is  well--known  that   (\ref{eqN})   defines  a  real  norm  on  the  real  space  $X$.  But  besides we have, for any   $ \lambda = \lambda_1 \e   + \lambda_2 \edag \in \bc$,  $ \lambda_1 , \, \lambda_2 \in  \C ( \i ) $  and  any  $ x= \e x_1 + \edag x_2 \in X$:
$$
\begin{array}{rcl}
\| \lambda x \|_X & = &  \| \e ( \lambda_1 x_1 )  + \edag ( \lambda_2 x_2 ) \|_X  = \displaystyle \frac{1}{\sqrt{2}}  \sqrt{  \| \lambda_1 x_1 \|_1^2 + \| \lambda_2 x_2 \|_2^2 }
\\  &  &  \\  &  = &
\displaystyle  \frac{1}{\sqrt{2}}   \sqrt{   | \lambda_1 |^2 \cdot \| x_1 \|_1^2  + | \lambda_2 |^2 \cdot \| x_2 \|_2^2 \, }
\\  &  &   \\  & \leq  &
\displaystyle    \frac{1}{ \sqrt{2}}   \sqrt{    2 \cdot \| x \|_X^2  \left(  | \lambda_1 |^2  + | \lambda_2 |^2  \right) \, }
\\   &  &  \\   & =   &
\displaystyle    \| x \|_X \cdot   \sqrt{   | \lambda_1 |^2   + |  \lambda_2 |^2 \, }   =  \sqrt{2} \,  | \lambda | \cdot \| x \|_X \, .
\end{array} $$

Since  obviously  for any     $\lambda  \in  \C ( \i )$  we  have  that
$$  \| \lambda x \|_X =  |  \lambda  | \cdot  \| x \|_X \, ,$$

then  (\ref{eqN})  defines  a   norm  of     $\C ( \i )$--type.  Paradoxically,     (\ref{eqN})
proves to be  not  only a  $\C ( \i )$--norm  but  a   $\C ( \j )$--norm   as  well.  Indeed, take  $ \lambda = a + \j b   =   ( a - \i b ) \e  + ( a + \i b ) \edag \in \C ( \j ) \subset \bc$,  $a, \, b \in \R$,  then
$$   \begin{array}{rcl}
\| \lambda x \|_X & = &  \displaystyle   \|  \e ( a - \i b) x_1  +  \edag  ( a + \i b )  x_2  \, \|_X
\\  &  &  \\  & = &
\displaystyle  \frac{1}{ \sqrt{2} }  \sqrt{   | a - \i b |^2 \cdot \| x_1 \|_1^2   + | a + \i b |^2 \cdot \| x_2 \|_2^2 \, }
\\  &  &  \\  & = &
| \lambda | \cdot \| x \|_X \,  .
\end{array}$$

In the same  way, a  pair  of  $\C ( \j )$--linear  normed  spaces  generates  a  bicomplex  module  with  a       norm   of   $\C ( \j)$--type,  which  again  is  of  $ \C ( \i )$--type  also.

\medskip

Thus  any  pair    of  complex  normed  spaces  generates  a  bicomplex module,  with  a   norm  of  complex  type.

\bigskip


\section{$\D$--valued  norm on  $\bc$--modules}

\begin{subDn}
Let  $X$  be  a  $\bc$--module. A map
$$  \| \cdot \|_\D : X \longrightarrow \D^+$$
is said to be a hyperbolic norm on $X$ if it satisfies the following properties:

\begin{enumerate}

\item[({\sc i})]  $ \| x \|_\D =0$ if and only if  $x=0$;

\medskip

\item[({\sc ii})]  $ \| \mu x \|_\D =   | \mu |_\k \cdot  \| x \|_\D  $   \;  $\forall  x \in X$,  $\forall  \mu \in \bc$;

\medskip

\item[({\sc iii})]  $ \|  x  +  y \|_\D \lessdot   \| x \|_\D   +    \| y \|_\D  $  \;  $ \forall  x, y \in X$.

\end{enumerate}
\end{subDn}

\medskip

\begin{subParr}
If  $X$  is of the form  $X=  \e \cdot X_1 + \edag \cdot X_2$, with  $X_1$ and $X_2$  being  two  arbitrary  $\C ( \i)$--linear  normed  spaces  with  norms   $\| \cdot  \|_1$  and  $\| \cdot  \|_2$,  then  $X$  can be endowed  canonically  with the hyperbolic  norm  by the formula
\end{subParr}
\begin{equation}\label{def hyperbolic norm BC module}
\| x \|_\D = \|  \e x_1 + \edag x_2 \|_\D  :=  \| x_1 \|_1 \cdot \e  +  \|  x_2 \|_2 \cdot \edag  .
\end{equation}

That  this  is  a  hyperbolic  norm   is  easily verified as  follows:
$$  \begin{array}{rcl}
\|  \mu x \|_\D  & = &  \| ( \mu_1 \e  + \mu_2 \edag )  \cdot  ( \e \cdot x_1  + \edag \cdot x_2 ) \|_\D
\\   &  &  \\  & = &
\|  \e ( \mu_1 x_1)  + \edag ( \mu_2 x_2) \|_\D
\\  &  &  \\  & = &
\| \mu_1  x_1 \|_1 \cdot \e  +  \| \mu_2 x_2 \|_2 \cdot \edag
\\   &  &  \\   & = &
| \mu_1 | \cdot \| x_1 \|_1 \cdot \e  + |  \mu_2 | \cdot \| x_2 \|_2 \cdot \edag
\\  &  &  \\  &  = &
(  |\mu_1 |  \cdot \e+  |\mu_2 | \cdot \edag ) \cdot (  \| x_1 \|_1 \cdot \e  +  \| x_2 \|_2 \cdot \edag  )
\\   &  &  \\  &  =  &    |\mu|_\k  \cdot \|  x\|_\D  ;
\end{array}  $$

$$ \begin{array}{rcl}
\| u+v\|_\D & = &  \|  (\e \cdot u_1  + \edag \cdot u_2)    +   (\e \cdot w_1  + \edag \cdot w_2)   \|_\D
\\  &  &  \\  & = &
\| \e ( u_1 + w_1) + \edag ( u_2 + w_2 ) \|_\D
\\   &  &  \\  & = &
\|u_1 + w_1 \|_1 \cdot \e  +  \| u_2 + w_2 \|_2 \cdot \edag
\\  &  &  \\  &  \lessdot  &
(\| u_1 \|_1 + \|w_1\|_1 )\cdot \e + (\| u_2 \|_2  + \| w_2 \|_2 ) \cdot \edag
\\  &  &  \\  & = &
\|u\|_\D  +  \| v \|_\D .
\end{array}  $$

We can compare (\ref{def hyperbolic norm BC module}) with the   (real--valued)   norm    given  by formula  (\ref{eqN}), and we easily see that
\begin{equation}\label{lunes 5 nov}
\begin{array}{rcl}
\displaystyle  \left| \| x \|_\D \right|  &    = &
|  \e \cdot  \| x_1 \|_1  +  \edag \cdot  \|  x_2  \|_2  |
\\  &  &  \\  & = &
\displaystyle  \frac{1}{  \sqrt{2} }  \cdot  \sqrt{   \,  \| x_1 \|_1^2  +  \| x_2   \|_2^2  \,  }
\\   &  &  \\   & = &       \| x \|_X   \, ,
\end{array}
\end{equation}

exactly as it happens when  $X = \bc$.

\medskip

Let us now consider a $\bc$--module $X$ endowed with a real--valued   norm $\| \cdot \|$  of  complex  type.  Since  $X = \e \cdot X_\e   +  \edag  \cdot X_\edag $  where   $X_\e$  and  $X_\edag$  are  $\C (\i)$--  and  $\C (\, \j)$--linear  spaces, the  restrictions  of  the  norm  on  $X$  onto   $X_\e$  and  $X_\edag$  respectively:
$$  \|  \cdot \|_{X_\e}  :=  \|  \cdot  \|  \mid_{X_\e}   \qquad  {\rm and}  \qquad    \|  \cdot \|_{X_\edag}  :=        \|  \cdot  \|  \mid_{X_\edag}$$

 induce norms  on the  complex  spaces  $X_\e$  and  $X_\edag$    respectively. These norms, in accordance with the process outlined before, allow us to endow  the  bicomplex  module  $X$  with  a  new hyperbolic norm, given  by  the  formula
\begin{equation}\label{lunes 5 nov2}
 \|  x \|_\D  :=  \|  \e x \|_{X_\e}  \cdot  \e  +  \| \edag x \|_{X_\edag} \cdot \edag .
 \end{equation}

\bigskip

Conversely,  if the  $\bc$--module  $X$  has  a  hyperbolic  norm  $\| \cdot   \|_\D$,  then the formula
\begin{equation}\label{norm from hyperbolic norm in bc module}
\|  x \|_X :=  \left|  \| x \|_\D  \right|
\end{equation}

defines a  norm  on $X$.  We now show that in fact   (\ref{norm from hyperbolic norm in bc module})  defines  a  usual norm in the real  space  $X$. To do so we observe that:

\begin{enumerate}

\item[(a)]  $\| x \|_X = 0$  if and only if  $\|x\|_\D =0$    which  holds   if and only if  $x=0$,  since   $\| \cdot \|_\D$ is  a  hyperbolic  norm;

\medskip

\item[(b)]   given  any  $\mu \in \R$ and  $x \in X$,  then
$$\begin{array}{rcl}
\| \mu x \|_X  & =  &    |   \| \mu x \|_\D | = |  |\mu|_\k \cdot \| x \|_\D |
\\   &  &   \\   &  =   &
|  |\mu| \cdot \| x \|_\D |  =     |\mu|   \cdot   | \| x \|_\D |
\\   &   &   \\   &    =   &      |\mu|   \cdot    \| x \|_X;
\end{array}    $$

\medskip

\item[(c)]  for  any  $x, \, y \in X$ there holds:
$$ \begin{array}{rcl}
 \| x + y \|_X & = &    | \| x + y \|_\D |
 \\   &   &   \\   &  \leq &
 | \| x \|_\D  + \| y \|_\D |   \;  \leq \;  |  \| x \|_\D |  +  |    \|  y  \|_\D  |
 \\    &  &   \\    & =  &
 \| x \|_X  + \| y  \|_X  .
 \end{array}  $$
Note  that  we  have  used  (\ref{desig hyperb positiv implies desig Eucl}) in this reasoning.
\end{enumerate}

Finally,  we have to prove that this real norm satisfies  the   required  inequality  in Definition  \ref{Def bicomplex norm}. Given $W \in \bc$ and $x \in X$   it is
$$\begin{array}{rcl}
\| W \cdot x \|_X   &  =   &      | \|  W \cdot x \|_\D |  =    |   | W |_\k \cdot  \| x \|_\D |
\\  &   &    \\   &   =  &
|  | W |_\k | \cdot  | \| x \|_\D |  \
\\    &   &   \\  &   \leq   &   \; \sqrt{2} \cdot  |W | \cdot  \|  x \|_X ;
\end{array}$$

here we have used (1) and (2)  from  Remark \ref{very usseful remark}.

\medskip

We can summarize our arguments as follows:

\begin{subtheorem}
A  $\bc$--module  $X$  has  a    real--valued     norm  if  and  only  if   it  has  a  hyperbolic  norm,   and these norms  are connected  by  equality   (\ref{norm from hyperbolic norm in bc module}).
\end{subtheorem}

\bigskip

As we already saw in the case of $\bc$, the two norms on a  $\bc$--module  $X$ (the real valued norm  $\| \cdot \|_X$  and the hyperbolic  norm    $\| \cdot \|_\D$) give rise to equivalent notions of convergence.

\begin{subDn}
A  sequence  $\{ x_n \}_{n \in \N} $  in  $X$  converges  to  $x_0  \in X$  with  respect  to  the    norm  $\| \cdot \|_X$  if   for every $\epsilon >0$  there  exists  $N \in \N$ such that for any    $  n \geq N$  it is:
$$   \|  x_n - x_0 \|_X <  \epsilon . $$

\end{subDn}

Similarly  for the convergence with  respect to the hyperbolic  norm.

\begin{subDn}
A  sequence  $\{ x_n \}_{n \in \N} $  in  $X$  converges  to  $x_0  \in X$  with  respect  to  the hyperbolic    norm  $\| \cdot \|_\D $  if   $ \forall  \epsilon >0$  there  exists  $N \in \N$ such that $\forall  n \geq N$  there  holds:
$$   \|  x_n - x_0 \|_X \lessdot   \epsilon . $$

\end{subDn}

From  equalities  (\ref{lunes 5 nov}), (\ref{lunes 5 nov2}) and  (\ref{norm from hyperbolic norm in bc module})  it  follows  that  a  sequence  $\{ x_n \}_{n \in \N} $   converges  to   $x_0 \in X$  with  respect to the  norm  $\| \cdot \|_X$   if  and only  if  it  converges  to  $x_0$  with  respect to the  hyperbolic  norm  $\| \cdot \|_\D$.

\bigskip


\section{Bicomplex  modules with  inner  product.}

\begin{subDn}\label{def inner product}   (Inner  product).  Let $X$ be a  $\mathbb B
\mathbb C$-module. A  mapping
$$  \displaystyle  \langle \cdot  \, , \, \cdot \rangle  :  X \times X
\to  \mathbb B \mathbb C  $$

is said to be a  $\bc$--inner,  or  $\bc$--scalar,  product on $X$ if it satisfies
the following  properties:

\begin{enumerate}

\item  $\langle x , y+z \rangle  =   \langle x , y \rangle  +
\langle x , z \rangle  \;  $   for  all   $ x,y,z \in X $;

\medskip

\item  $\langle \mu \,   x ,  y \rangle = \mu \, \langle x,y
\rangle \;  $   for  all   $  \mu \in \mathbb B \mathbb C , \;  $   for  all   $ x,y  \in X$;

\medskip

\item  $\langle x,y \rangle = \langle y , x \rangle^\ast     \;  $   for  all   $  x,y \in X$;

\medskip

\item   $  \langle x ,x \rangle \in \mathbb D^+$,  and   $  \langle x,x \rangle =0 $  if  and  only  if    $x =0$.

\end{enumerate}

\end{subDn}

Note that property  (3)  implies  the  fo\-llo\-wing:
$$  a + b \, \mathbf j  = \langle  x,x   \rangle  =  \langle x,x
\rangle^\ast  =  \overline a  -  \overline b \, \mathbf j \, \quad
{\rm  with}  \quad   a,b  \in  \mathbb C (\i ) ,$$

\medskip

hence     $a \in \mathbb R$,  $b=  \mathbf i \, b_1$, with
$b_1 \in \mathbb R$,  that is,  $\langle x ,x \rangle  \in \mathbb
D$.   For    reasons that will be clear  later,  we are strengthening this property, by requiring additionally  that     $\langle x ,x \rangle $   will be not only  in   $ \mathbb
D$  but in   $ \mathbb
D^+$,  the  set  of  ``positive"  hyperbolic  numbers.

\medskip

Also  from  Properties  (2)  and (3),  we see that  for  any  $\mu \in  \mathbb B \mathbb C $ and any $x, y \in X$   one has:
$$  \langle  x  , \mu \,   y \rangle  =  \mu^\ast \, \langle x,y \rangle \, .$$

\medskip

\begin{subDn} A $ \mathbb B  \mathbb C$--module $X$ endowed with a  bicomplex  inner  pro\-duct   $\langle \cdot , \cdot  \rangle$ is said to be a   $ \mathbb B  \mathbb C$--inner  product   module.
\end{subDn}


\begin{subParr}\label{inner product on the sum of two complex spaces}  We consider now, as in Sections  \ref{construction}  and  \ref{norm construction},  two  $\C ( \i)$--linear  spaces   $X_1$  and  $X_2$. In addition we will assume that both  are  inner  product  spaces,  with  inner  products  $ \langle \cdot , \cdot \rangle_1$  and   $ \langle \cdot , \cdot \rangle_2$, and corresponding norms  $\| \cdot \|_1$  and  $\| \cdot \|_2$. We can then show  that  the  formula
\begin{equation}\label{avestruz}
 \begin{array}{rcl}
\langle x , w \rangle_X & = &   \langle \e x_1  + \edag x_2  , \e  w_1  + \edag w_2  \rangle_X
\\  &  &  \\   &  :=  &
\e \langle x_1 , w_1 \rangle_1   + \edag  \langle  x_2  , w_2  \rangle_2
\end{array}
\end{equation}
\end{subParr}
defines  a  bicomplex  inner  product on  the  bicomplex  module   $X =  \e X_1 + \edag X_2$.   We  begin  by proving  distributivity.

\begin{enumerate}

\item   $$ \begin{array}{rcl}
\langle x , y + w \rangle_X  &  =  &   \langle  \e x_1 + \edag x_2 ,  ( \e  y_1  + \edag y_2 ) +
\\  &  &  \\  &  &
\qquad  \qquad  \qquad     + \;    ( \e w_1 + \edag w_2 ) \rangle_X
\\  &  &  \\   &  = &
\langle \e x_1 + \edag x_2  ,  \e  (y_1  + w_1) + \edag ( y_2 + w_2 ) \rangle_X
\\   &  &  \\   & = &
\e \langle x_1 , y_1 + w_1 \rangle_1  + \edag \langle x_2 , y_2 + w_2 \rangle_2
\\  &  &  \\   & = &
\displaystyle \left(  \e \langle x_1 , y_1 \rangle_1  + \edag \langle x_2 , y_2 \rangle_2 \right)  \: +
\\  &  &  \\   &  &  \qquad   \quad
\displaystyle   + \;     \left(  \e \langle x_1 , w_1 \rangle_1  + \edag \langle x_2 , w_2 \rangle_2 \right)  \; =
\\  &  &  \\  &  = &
\langle x, y \rangle_X +    \langle x, w \rangle_X \, .
\end{array} $$

\medskip

\item    Next,  we  deal  with  the  linearity.    Given  $\mu = \mu_1 \e + \mu_2 \edag \in \bc$,   then
$$
\begin{array}{rcl}
\langle \mu \, x ,  y \rangle_X & = &  \langle    ( \mu_1 \e  + \mu_2 \edag ) (      \e x_1 + \edag x_2 )  ,  \e y_1 + \edag y_2  \rangle_X
\\  &  &  \\   & = &
\langle      \e   (\mu_1  x_1)   + \edag (   \mu_2 x_2)    , \e  y_1  + \edag    y_2     \rangle_X
\\  &  &  \\  & = &
\e  \langle    \mu_1    x_1 , y_1 \rangle_1  + \edag \langle  \mu_2   x_2 , y_2 \rangle_2
\\  &  &  \\  & = &
\e \mu_1\langle x_1 , y_1 \rangle_1  + \edag   \mu_2  \langle  x_2 , y_2 \rangle_2
\\  &  &  \\  & = &
( \e \mu_1 + \edag \mu_2 )( \e \langle x_1 , y_1 \rangle_1  + \edag \langle x_2 , y_2 \rangle_2 )
\\  &  &  \\   & = &
\mu \langle  x, y \rangle_X \, .
\end{array}  $$

\medskip

\item    Next  step  is  to  prove  what  is the analogue  of  being   ``Her\-mi\-tian  for  an  inner  product"   in  our  situation.
$$ \begin{array}{rcl}
\langle y , x \rangle_X^\ast  & = &  \langle \e y_1 + \edag y_2 , \e x_1 + \edag x_2 \rangle_X^\ast
\\  &  &  \\  & = &
\displaystyle \left(  \e \langle y_1 , x_1 \rangle_1 + \edag \langle y_2 , x_2 \rangle_2  \right)^\ast
\\  &  &  \\  & = &
\e \overline{ \langle y_1 , x_1 \rangle_1 }  + \edag  \overline{  \langle y_2 , x_2 \rangle_2 }
\\  &  &  \\  & = &
\e \langle x_1 , y_1 \rangle_1  + \edag \langle x_2 , y_2 \rangle_2 = \langle x , y \rangle_X \, .
\end{array} $$

\medskip

\noindent
Note that the  ``inner   product  square"
\begin{equation}\label{dos estrellas}
\langle  x ,x \rangle_X  = \e \| x_1 \|_1^2  + \edag \| x_2 \|_2^2\in  \mathbb D^+ \, .
\end{equation}
is  always a positive hyperbolic number, without a need for any  additional  requirements  on  $ \langle \cdot , \cdot \rangle_1$  and    $ \langle \cdot , \cdot \rangle_2$.
 The  marvelous  fact  here  is  that  with  the  above  described  process,  we  constructed  a class  of  bicomplex  modules  with  inner  product  in  such  a  way  that    property   (\ref{dos estrellas})    verifies  always.

\medskip

\item   Finally,  we  prove  the  non--degeneracy.
$$ \langle x,x \rangle_X =0    \;  \Longleftrightarrow \;   \e \| x_1 \|_1^2  + \edag \| x_2 \|_2^2 =0 \:   \Longleftrightarrow$$
$$ \,  \quad     x_1=0 , \; x_2=0 \; \;      \Longleftrightarrow \; \;    x=0.  $$

\end{enumerate}

As  we already noticed,  the  ``inner  product square"  is  a  hyperbolic  positive  number, and this suggests     the possibility of introducing  a   hyperbolic  norm  on  an  inner  product  $\bc$--module  consistent  with  the  bicomplex    inner  product.  Let  us  show    that  this is,  indeed,  possible.  Set
\begin{equation}\label{hyperbolic norm from inner prod}
\| x \|_\D = \| \e \cdot x_1  + \edag \cdot x_2 \|_\D :=  \langle x ,x \rangle^{1/2} ,
\end{equation}

where  $ \langle x ,x \rangle^{1/2} $ is the hyperbolic  number given by
$$ \begin{array}{rcl}
\langle x ,x \rangle^{1/2}  & =    & \left(   \e \langle x_1 ,x_1  \rangle_1  + \edag  \langle x_2  ,x_2   \rangle_2 \right)^{1/2}
\\   &  &  \\   &  = &
 \e \langle x_1 ,x_1  \rangle_1^{1/2}    + \edag  \langle x_2  ,x_2   \rangle_2^{1/2}
 \\   &  &  \\   & = &
\e \cdot \| x_1 \|_1  + \edag \| x_2 \|_2  \in \D^+ .
\end{array}  $$

It is easy to verify that the hyperbolic  norm satisfies the properties one would expect.

\begin{enumerate}

\item[({\sc i})]    $\| x \|_\D = \|  \e x_1 + \edag x_2 \|_\D  =0$ if and only if $x_1=0$ and  $x_2=0$.

\medskip

\item[({\sc ii})]    Given   $\mu = \mu_1 \e + \mu_2 \edag \in \bc$,  one has:
$$  \begin{array}{rcl}
\|  \mu x \|_\D  & = &   \langle   \mu  x       ,   \mu  x  \rangle^{1/2}
\\  &  &  \\  & = &
(\mu \mu^\ast)^{1/2} \cdot  \langle x ,x  \rangle^{1/2}
\\  &  &  \\  & = &
\left( | \mu_1 | \e  + | \mu_2 | \edag \right) \cdot \| x \|_\D
\\  &  &  \\  & =   &
| \mu |_\k  \cdot   \| x \|_\D \, .
\end{array}   $$

\medskip

\item[({\sc iii})]    For  any   $x , y \in X$,   it    follows:
$$   \begin{array}{rcl}
\|  x+ y \|_\D  & = &   \langle     x  + y        ,     x  + y  \rangle^{1/2}
\\  &  &  \\   & = &
\left(   \langle     x_1   + y_1         ,     x_1  + y_1   \rangle_1 \cdot \e  +   \langle     x_2   + y_2         ,     x_2  + y_2   \rangle_2 \cdot \edag  \right)^{1/2}
\\  &  &  \\   & = &
\left(   \|     x_1   + y_1  \|_1^2      \cdot \e  +   \|     x_2   + y_2   \|_2^2     \cdot \edag  \right)^{1/2}
\\  &  &  \\   &   \lessdot   &
\left(   \left(    \|     x_1 \|_1   +   \|  y_1  \|_1  \right)^2      \cdot \e  +   \left(    \|     x_2  \|_2    +   \|   y_2   \|_2  \right)^2      \cdot \edag  \right)^{1/2}
\\  &  &   \\   & = &
\| x \|_\D  + \| y \|_\D .
\end{array}  $$

\end{enumerate}


\begin{subParr}\label{bicomplex norm from inner product} On the other hand, it is also possible to endow  the   bicomplex  module   $X = \e X_1 + \edag X_2$  with  a real--valued norm naturally related to the bicomplex inner  product discussed above. To this purpose, we simply set
\begin{equation}\label{puntito}
\begin{array}{rcl}
\| x \|_X^2  & :=  &  \displaystyle  \frac{1}{2}  \left(   \langle  x_1 , x_1  \rangle_1  +  \langle x_2 , x_2 \rangle_2 \right)
\\  &  &  \\  & = &
\displaystyle  \frac{1}{2}  \left(  \| x_1 \|_1^2  + \| x_2  \|_2^2  \right)  \, ,
\end{array}
\end{equation}
\end{subParr}

and we compare this definition with  formulas   (\ref{eqN})  and   (\ref{dos estrellas}); we immediately see that (\ref{puntito})    defines   a  Euclidean--type     norm  on  $X$  which  is  related to the norms  on $X_1$  and  $X_2$  exactly in the same way as the Euclidean  norm on $\bc$ is related to the modules   of the coefficients  of its idempotent  representation.  The  relation between  (\ref{puntito})  and  (\ref{hyperbolic norm from inner prod})  is as one would expect:
$$ \begin{array}{rcl}
\displaystyle | \| x \|_\D  | & = & \displaystyle | \, \e \cdot \| x_1 \|_1  \, + \, \edag \cdot \| x_2 \|_2 \, |
\\  &  &  \\  & = &
\displaystyle  \frac{1}{\sqrt{2}}  \cdot \sqrt{  \, \| x_1 \|_1^2 + \| x_2 \|_2^2 \, }  \; = \;   \| x \|_X \, .
\end{array}  $$

\begin{subDn}  A $\bc$--inner  product  module $X $     is said to be a  bicomplex  Hilbert  module  if  it is  complete  with  respect  to  the  metric  induced  by  its  Euclidean--type    norm  generated  by the  inner  producet; this is equivalent to say that  $X$  is  complete  with  respect  to  the  hyperbolic  norm  generated   by  the  inner  product  square.
\end{subDn}

\begin{subremark}
We can restate this definition by saying that a  bicomplex  Hilbert  module is  a  4--tuple   $( X,  \| \cdot \|_X , \| \cdot \|_\D , \langle \cdot , \cdot \rangle )$  where  any  Cauchy  sequence is  convergent.
\end{subremark}

\medskip

It  follows  from  (\ref{puntito})  that   if    $X   = \e X_1 + \edag X_2 $,   where  $X_1$,  $X_2$ are  $\C ( \i)$--linear  inner  product  spaces,  then  $X$      is  a  bicomplex  Hilbert  space  if  and  only if  $X_1$  and  $X_2$  are  complex  Hilbert  spaces.

\medskip


\begin{subEx}\label{Ejemplo caso bicomplejo}  A bicomplex inner product in $ \mathbb B \mathbb C$ can be given through the formula:
\begin{equation}\label{inner product on BC}
\displaystyle  \left\langle Z , W \right\rangle :=  Z \, W ^\ast\, .
\end{equation}
\end{subEx}

Properties  (1),   (2)     and  (3)     of  Definition  \ref{def inner product}  are  obvious.  For  the
fourth  property  recall  that if $Z=  z_1  +  \mathbf j \, z_2  =  \beta_1 \cdot \e + \beta_2 \cdot \edag  $, then
\begin{equation}\label{HN}
\begin{array}{rcl}
\langle Z ,Z \rangle  &  =  &     Z \, Z ^\ast
\\   &   &   \\    &   =   &
(\beta_1 \cdot \e + \beta_2 \cdot \edag ) \cdot   ( \overline{  \beta}_1 \cdot \e +   \overline{ \beta}_2 \cdot \edag ) \;
\\   & & \\ & =  &     | \beta_1 | \cdot \e  + | \beta_2 | \cdot \edag
\\  &  &   \\   & = &
  \left(  | z_1 |^2
+ | z_2 |^2  \right)   -  2  Im  (z_1 \, \overline z_2 ) \,
\mathbf i \, \mathbf j  \, \in \D^+ .
\end{array}
\end{equation}

Thus   $ \langle Z ,Z \rangle =0$  if  and  only  if  $z_1 =0=z_2$  or equivalently
$Z=0$.  Hence  the  general  definition    of  the  inner  product
is simply a generalization of this definition for the bicomplex--valued inner product in  $ \mathbb B
\mathbb C$. Note in particular that    non--degeneracy    would be  lost  if  we
had used  any  of  the   two  other  conjugations  on  $ \mathbb B
\mathbb C$.  Moreover,  the  inner   product   square  of  any  element  is
a  hyperbolic  number  in  $\mathbb D^+$.  This  is the reason why, in    \cite{LMR2010}, we required   that
 the  inner     product   square  of  any  element   (in  an  arbitrary  bicomplex  inner  product  module) had to be
a  positive     hyperbolic   number.   Note  that  in  this  example  the  real   part  of  the  hyperbolic  number   (\ref{HN}) is  already  a  good  norm.  \\

Let us now compare these   two  approaches    with Subsection   \ref{inner product on the sum of two complex spaces}.
The  idempotent   decomposition  of  $\bc$:
$$  \bc  =  \C ( \i ) \e  + \C ( \i ) \edag \, ,    $$

can be reformulated,     by saying that $X_1  =    \C ( \i ) \e$,    $ X_2 =   \C ( \i ) \edag  $
and  with  the   $ \C ( \i )$--valued  inner  products
$$  \displaystyle  \langle \beta_1 \e , \beta_1' \e \rangle_1  :=   \beta_1  \, \overline{   \beta_1' }  \quad  \forall  \beta_1 \e , \;   \beta_1' \e  \in X_1 ,  \;  \beta_1 , \, \beta_1' \in \C ( \i ) ,$$

$$  \displaystyle  \langle \beta_2 \edag , \beta_2' \edag \rangle_2  :=   \beta_2  \, \overline{ \beta_2' }  \quad  \forall  \beta_2 \edag , \;   \beta_2' \edag  \in X_2 ,  \;  \beta_2 , \, \beta_2' \in \C ( \i ) .$$

\medskip

Now, in accordance with (\ref{avestruz}), if we  take   $Z= \beta_1 \e  + \beta_2 \edag $, and $  W=  \beta_1' \e    +    \beta_2' \edag$, we obtain   that
$$  \displaystyle    \langle Z ,W \rangle_X  :=   \langle  \beta_1  \e , \beta_1' \e \rangle_1 \, \e  + \langle \beta_2 \edag , \beta_2' \edag \rangle_2 \, \edag \, ,$$

and therefore
$$ \begin{array}{rcl}
 \displaystyle    \langle Z ,W \rangle_X   &  =  &      \beta_1 \, \overline{  \beta_1' }  \e  +  \beta_2 \,  \overline{ \beta_2' }   \edag  =   \left(    \beta_1 \e  +  \beta_2 \edag \right)   \left(   \overline{ \beta_1' }  \e + \overline{ \beta_2' }   \edag \right)
 \\   &  &   \\   &  =  &    \displaystyle     ZW^\ast  = \langle Z ,W \rangle .
 \end{array} $$

\medskip

This  means  that    the  inner   product   (\ref{inner product on BC})    on  $\bc$   coincides  with  the  one  generated  by  its   idempotent  decomposition. In addition, the Euclidean norm on $\bc$  coincides with the  norm   on  $\bc$  given  by formula   (\ref{puntito}).  In  other  words   the  Euclidean  norm of  bicomplex  numbers  is  the  Euclidean--type      norm  on  the  bicomplex  inner  pro\-duct  module  $\bc$.

\medskip

Finally, if  one  takes,  in   (\ref{HN}),  $Z$  in the idempotent  form  $ Z =  \beta_1 \e + \beta_2 \edag $  then  the  inner  product  square  is
$$  \langle  Z ,Z  \rangle_X  = Z \cdot Z^\ast  = | Z |_\k^2  =   | \beta_1 |^2  \cdot \e  +   | \beta_2 |^2  \cdot \edag  $$

giving    us  a  hyperbolic  norm on the $\bc$--module  $\bc$.

\bigskip


\section{Inner  products  and  cartesian  decompositions}
\label{cartesian decomp}
\begin{subParr}
\label{Sub231}
We will now assume  that      a   bicomplex   inner   product   module  $X$  has  a   bar-involution, and we denote by   $X_{1, bar}$  the
$\mathbb C (\mathbf j)$--linear  space  of  its    $\mathbb C (\mathbf j)$--elements.
Thus  $X =   X_{1, bar}  + \mathbf i \, X_{1, bar}$  and   $X_{ \mathbb C (\mathbf j) } =   X_{1, bar}  \oplus_{ \mathbb C (\mathbf j) }   \mathbf i \, X_{1, bar}$,   when  we restrict  the  scalar  multiplication  to  $\mathbb C (\mathbf j)$. Nevertheless, the inner  product  takes  values in  $\mathbb B \mathbb C $.  If we now define an inner product
$$ \langle \cdot  , \cdot  \rangle_{bar}:  X_{ \mathbb C (\mathbf j) }  \times  X_{ \mathbb C (\mathbf j) }  \to  \mathbb C (\mathbf j) $$
\end{subParr}

by
$$ \langle x  ,y \rangle_{bar}:=  \Pi_{1, \mathbf j}  \left(   \langle x  ,y \rangle  \right) \in    \mathbb C (\mathbf j)  \, ,$$

we see  that  $ \langle \cdot  , \cdot  \rangle_{bar}$  is  an  inner  product  on the   $    \mathbb C (\mathbf j)$--linear  space  $X_{ \mathbb C (\mathbf j) }$,  and that it  induces  an  inner  product  on the  $   \mathbb C (\mathbf j)$--linear  space  $X_{1,bar}$  as  follows:
$$ \langle x_1   ,y_1 \rangle_{  \mathbb C (\mathbf j) }   :=        \langle x_1  ,y_1 \rangle_{bar}  \quad  \forall  x_1, \, y_1    \in X_{1,bar}   \, .$$

The relationship between the bicomplex the $\C ( \, \j)$--valued inner products,  is established by the following formula:
if   $ x=  x_1 +  \mathbf i \, x_2$,  $y =  y_1 + \mathbf i \, y_2$  are  in  $X$,     then
\begin{equation}\label{prod int inducido en X desde X1bar}
\begin{array}{l}
\langle  x , y \rangle  =   \langle x_1  +  \mathbf i \, x_2  ,   y_1 + \mathbf i \, y_2  \rangle  =
\\  \\   \qquad  \  =  \;  \langle  x_1 , y_1 \rangle_{\mathbb C (\mathbf j) }   \,   +  \langle  x_2 , y_2 \rangle_{\mathbb C (\mathbf j) }  \,  +
\\  \\   \qquad  \qquad   \ + \;    \mathbf i \, \left(   \langle  x_2 , y_1 \rangle_{\mathbb C (\mathbf j) }     -    \langle  x_1 , y_2 \rangle_{\mathbb C (\mathbf j) }     \right)   \, .
\end{array}
\end{equation}

This   formula  provides  even  more  information.  Indeed,  if $X$  is  just  a  bicomplex  module,  without  any  bicomplex  inner  product, and if    $ X_{1, bar}$  has a  $ \mathbb C (\mathbf j) $--valued   inner  product,  then  formula    (\ref{prod int inducido en X desde X1bar})  tells  us  how  to  endow  $X$  with  a  bicomplex  product   which   extends     the one on $X_{1, bar}$.


\begin{subParr}\label{Sub232}   We can treat in exactly the same way the  case in which $X$ admits a  $\dagger$--involution.  In  this  case    $X =   X_{1, \dagger}  + \mathbf j \, X_{1, \dagger}  $,  and    $X_{  \mathbb C (\mathbf i)}   =   X_{1, \dagger}  \oplus_{\mathbb C (\mathbf i) }    \mathbf j \, X_{1, \dagger}  $.    The  mapping
$$ \langle \cdot  , \cdot  \rangle_{\dagger}:  X_{ \mathbb C (\mathbf i) }   \times  X_{ \mathbb C (\mathbf i) }  \to  \mathbb C (\mathbf i) $$
\end{subParr}
given  by
$$ \langle x  ,y \rangle_{\dagger}:=  \Pi_{1, \mathbf i}  \left(   \langle x  ,y \rangle  \right) \in    \mathbb C (\mathbf i)  \, $$

is  an  inner  product  on  the   $  \mathbb C (\mathbf i) $--linear  space  $X_{  \mathbb C (\mathbf i) } $.  The  induced  inner  pro\-duct  on  the  $ \mathbb C (\mathbf i) $--linear  space    $ X_{1, \dagger}  $  for  any   $s_1 , \,  t_1  \in   X_{1, \dagger}  $  is:
$$ \langle s_1   ,t_1 \rangle_{  \mathbb C (\mathbf i) }   :=        \langle s_1  ,t_1 \rangle_{\dagger}   \, .$$

\medskip

Now, if  $ x=  s_1 +  \mathbf j \, s_2$,  and    $y =  t_1 + \mathbf j \, t_2$  are  in  $X$  then
\begin{equation}\label{prod int inducido en X desde X1dagger}
\begin{array}{l}
\langle  x , y \rangle  =   \langle s_1  +  \mathbf j \, s_2  ,   t_1 + \mathbf j \, t_2  \rangle  =
\\  \\   \qquad  \quad    =  \;  \langle  s_1 , t_1 \rangle_{\mathbb C (\mathbf i) }     +    \langle  s_2 , t_2 \rangle_{\mathbb C (\mathbf i) }   + \mathbf j \, \left(   \langle  s_2 , t_1 \rangle_{\mathbb C (\mathbf i) }    -     \langle  s_1 , t_2 \rangle_{\mathbb C (\mathbf i) }      \right)   \, .
\end{array}
\end{equation}

\medskip

Again,  formula (\ref{prod int inducido en X desde X1dagger})  provides  an  extension of a  $\mathbb C (\mathbf i ) $--valued  inner  product  defined  on  $X_{ 1, \dagger }$  up  to  a  bicomplex   inner  product  on  $X$.

\bigskip


\begin{subEx}  \label{Ex 2a parte}   We  continue  now the  analysis   iniciated  in  Example    \ref{Ejemplo caso bicomplejo}.  
      The space $ X =  \mathbb B \mathbb C$  has  both  involutions,  so that  $X_{1,bar} = \mathbb C (\mathbf j)$  and     $X_{1,\dagger} = \mathbb C (\mathbf i)$.  Given  $ Z, \, W  \in    \mathbb B \mathbb C$,  with  $Z  = z_1 + \mathbf j \, z_2 = s_1 + \mathbf i \, s_2$  and   $W =  w_1 + \mathbf j \, w_2 = t_1 + \mathbf i \, t_2$,  $z_1$,  $z_2$, $w_1$,  $w_2 \in \mathbb C ( \mathbf i)$,   $s_1, \, s_2 , \, t_1, \, t_2 \in   \mathbb C (\mathbf j)$,  then
\end{subEx}
$$ \begin{array}{rcl}
\displaystyle \left\langle Z,W \right\rangle  &  =  &   \displaystyle    Z W ^\ast =   \left(   z_1 + \mathbf j \, z_2   \right)  \cdot  \left(  w_1 + \mathbf j \, w_2  \right)^\ast
\\  &  &  \\   & = &
 \displaystyle      \left(   z_1 +  \mathbf j \,   z_2   \right) \cdot  \left(   \overline w_1 - \mathbf j \,  \overline w_2  \right)
 \\   &  &   \\   &   =   &
\displaystyle     \left(   z_1 \, \overline w_1 +  z_2 \, \overline w_2 \right)  + \mathbf j \, \left(   z_2 \, \overline w_1   - z_1 \,  \overline w_2    \right)
\\  &  &  \\  &  =  &   \displaystyle     \left(   s_1 + \mathbf i \, s_2   \right)  \cdot  \left(  t_1 + \mathbf i \, t_2  \right) ^\ast
\\  &  &  \\   & = &
 \displaystyle      \left(   s_1 +  \mathbf i \,  s_2   \right) \cdot  \left(  t_1 ^\ast   - \mathbf i \, t_2^\ast  \right)  =  \left( s_1 \, t_1^\ast + s_2 \, t_2 ^\ast   \right)  + \mathbf i \, \left(          s_2 \, t_1^\ast      -    s_1 \, t_2 ^\ast      \right) \, ,
\end{array}  $$

which  implies  that  on  the   $\C (\j)$--linear  space   $\bc_{ \C ( \j ) }$  the  $\C ( \j )$--valued  inner  product  $ \langle \cdot , \cdot  \rangle_{bar} $  is  given  by
$$  \displaystyle  \left\langle Z , W \right\rangle_{bar}  = s_1 \, t_1 ^\ast+ s_2 t_2^\ast     \in \mathbb C (\mathbf j )$$

\medskip

and  on  the  $\C ( \i )$--linear  space   $\bc_{ \C ( \i )}$  the $ \C ( \i )$--valued  inner  product    $ \langle \cdot , \cdot  \rangle_\dagger $   is
$$  \displaystyle  \left\langle Z, W \right\rangle_\dagger =    z_1 \,  \overline w_1 +    z_2 \,  \overline w_2 \in \mathbb C ( \mathbf i ) \, ;$$

finally it is
$$  \langle s_1 , t_1 \rangle_{\mathbb C (\mathbf j ) }  :=  s_1\, t_1^\ast  \, ,  \qquad    \langle z_1 , w_1 \rangle_{\mathbb C (\mathbf i ) } :=   z_1 \, \overline w_1  \, ,   $$

\medskip

which  are  the usual inner  products  on    $\C ( \j )$   and  $\C ( \i )$.

\medskip


\section{Inner  products  and  idempotent  decompositions.}
\label{restrictions of inner product}


\begin{subParr}    Consider  a  bicomplex  module  $X$  with  an inner  product
$\langle \cdot , \cdot  \rangle$.  For any  $x , y \in X$, the inner product  $\langle
x , y \rangle$ is given by
\begin{equation}\label{eq A}
\begin{array}{rcl}
\langle
x , y \rangle  &  =   &     \beta_1 (x,y) \, \mathbf e  + \beta_2 (x,y) \,
\mathbf e^\dagger \, \\  &  &  \\   &  =  &
\gamma_1  (x,y)  \, \mathbf e   +  \gamma_2  (x,y) \, \mathbf e^\dagger \, ,
\end{array}
\end{equation}
\end{subParr}

where  $ \beta_\ell  (x,y) = \pi_{\ell , \mathbf i}    (  \langle x , y \rangle ) \in
\mathbb C  (   \mathbf i) $,    $ \gamma_\ell (x,y) = \pi_{\ell , \mathbf j} (  \langle x , y \rangle )
\in \mathbb C ( \mathbf j )$,  for  $\ell  = 1,2$.    Define
$$ \mathcal B_1 , \; \mathcal B_2 : X \times  X \to \mathbb C  ( \mathbf i)   \qquad   {\rm and}  \qquad    \mathcal F_1 , \, \mathcal F_2 : X \times X \to \mathbb C (\mathbf j  )  $$

by
\begin{equation}\label{def Bl}
 \mathcal B_\ell \, (x,y) := \pi_{\ell  , \mathbf i }   \, ( \langle x,y \rangle )  =
\beta_\ell\,  (x,y)     \in  \mathbb C ( \mathbf i )  \,
\end{equation}

and
\begin{equation}\label{def Fl}
    \mathcal F_\ell \, (x,y) := \pi_{\ell  , \mathbf j }   \, ( \langle x,y \rangle )  =
\gamma_\ell\,  (x,y)  \in  \mathbb C ( \mathbf j )   \, .
\end{equation}

We  will  prove  that  the  restrictions
\begin{equation}\label{prod int e i}
\langle \cdot , \cdot \rangle_{ X_{\mathbf e ,  \mathbf i } }  :=     \mathcal B_1  \mid_{  X_{\mathbf e}  \times X_{\mathbf e}  }  : X_{\mathbf e}  \times  X_{\mathbf e}   \to \mathbb C  ( \mathbf i )  \,  ;
\end{equation}

\begin{equation}\label{prod int e j}
\langle \cdot , \cdot \rangle_{ X_{\mathbf e , \,   \mathbf j } }  :=     \mathcal F_1  \mid_{  X_{\mathbf e}  \times X_{\mathbf e}  }  : X_{\mathbf e}  \times  X_{\mathbf e}   \to \mathbb C  ( \mathbf j )  \,  ;
\end{equation}

\begin{equation}\label{prod int e dagger i}
\langle \cdot , \cdot \rangle_{ X_{\mathbf e^\dagger ,  \mathbf i } }  :=     \mathcal B_2  \mid_{  X_{\mathbf e^\dagger}  \times X_{\mathbf e^\dagger}  }  : X_{\mathbf e^\dagger}  \times  X_{\mathbf e^\dagger}   \to \mathbb C  ( \mathbf i )  \,  ;
\end{equation}

\begin{equation}\label{prod int e dagger j}
\langle \cdot , \cdot \rangle_{ X_{\mathbf e^\dagger ,  \,  \mathbf j } }  :=     \mathcal F_2  \mid_{  X_{\mathbf e^\dagger}  \times X_{\mathbf e^\dagger}  }  : X_{\mathbf e^\dagger}  \times  X_{\mathbf e^\dagger}   \to \mathbb C  ( \mathbf j )  \,  ,
\end{equation}

are the usual  inner   products on the    $\mathbb C( \mathbf i )$--  and    $\mathbb C( \mathbf j )$--linear spaces
$ X_{\mathbf e}  $  and  $ X_{\mathbf e^\dagger} $  respectively.   Indeed,  for  any    $x,y,z \in  X_{\mathbf e} $     and  $\lambda \in \mathbb C  ( \mathbf i ) $,     we have:

\begin{enumerate}

\item[(a)]
$$ \begin{array}{rcl}
\langle x ,   y + z   \rangle_{ X_{\mathbf e ,  \mathbf i } }  &   =    &
\pi_{ 1  , \mathbf i }   \, ( \langle x , y+z \rangle ) = \pi_{1  ,  \mathbf i }   \, ( \langle
x,y \rangle + \langle x,z \rangle )  \\  &  &  \\ & = & \pi_{ 1  , \mathbf i }
\, ( \langle x,y \rangle ) +  \pi_{1  ,  \mathbf i }    \, ( \langle x, z \rangle )  \\  &  &   \\   &
=  &      \langle x, y \rangle_{ X_{\mathbf e ,  \mathbf i } }    +    \langle x, z \rangle_{ X_{\mathbf e ,  \mathbf i } }  \, .
\end{array} $$

\medskip

\item[(b)]
Note that  given  $\mathbf e \, x , \; \mathbf e \, y \in  X_{\mathbf e}    $  then
$$  \langle \mathbf e \, x , \mathbf e \, y \rangle = \mathbf e \,
\mathbf e^\ast \, \langle x,y \rangle = \mathbf e \, \langle x,y
\rangle \, ,$$

\noindent
which implies
\begin{equation}\label{eqU}
 \langle  \mathbf e \, x , \mathbf e \, y   \rangle_{ X_{\mathbf e ,  \mathbf i } }     = \pi_{1 , \mathbf i }  (
\mathbf e \, \langle x,y \rangle )  = \pi_{1 , \mathbf i }  ( \langle x,y \rangle
)  =  \beta_1 (x , \, y ) .
\end{equation}

\medskip

\noindent      Similarly we have       that
\begin{equation}\label{eqV}
 \langle  \edag  x , \edag  y   \rangle_{ X_{\mathbf e^\dagger} ,   \mathbf i } =      \pi_{2 ,  \,   \mathbf i }  (  \langle x,y \rangle )   =  \beta_2(x,y) \, ,
 \end{equation}

\noindent and
$$  \langle  \e \, x , \e \, y   \rangle_{ X_{\mathbf e ,  \,  \mathbf j } }      =       \pi_{1 , \,  \mathbf j }  (  \langle x,y \rangle )   ,   \; \;     \langle  \edag  x , \edag  y   \rangle_{ X_{\mathbf e^\dagger , \,   \mathbf j } } =      \pi_{2 , \,  \mathbf j }  (  \langle x,y \rangle )   . $$

\noindent
Now,  given  $ \lambda  \in  \mathbb C ( \mathbf i ) $  and      $\mathbf e \, x , \; \mathbf
e \, y \in   X_{\mathbf e}       $     then
$$\begin{array}{rcl}
\langle   \lambda \, (  \mathbf e \, x )  ,\mathbf e \, y   \rangle_{ X_{\mathbf e ,  \mathbf i } }     &  =  &
\pi_{1  , \mathbf i }   (
\langle   \lambda \, \mathbf e \, x ,  \mathbf e \, y  \rangle )
\\  &  &   \\  &   = &
\pi_{1 , \mathbf i }  (  \lambda \,  \langle  \mathbf e \, x , \mathbf e \, y
\rangle )  =  \lambda \, \pi_{1 , \mathbf i }  ( \langle \mathbf e \, x , \mathbf e
\, y \rangle )
\\  &  &   \\   &  =  &     \lambda \, \langle   \mathbf e \, x ,
\mathbf e \, y  \rangle_{ X_{\mathbf e  ,  \mathbf i } }    \, .
\end{array} $$

\medskip

\item[(c)]  One has:
$$  \langle  y , x  \rangle_{ X_{\mathbf e ,  \mathbf i } }      =   \pi_{1 , \mathbf i } \, ( \langle y ,
x \rangle )  =  \pi_{ 1 , \mathbf i }  \, ( \langle  x,y \rangle^\ast ) \, ;$$

\noindent
if we write
$$   \langle  x,y \rangle    =   \beta_1
(x,y)   \, \mathbf e  +   \beta_2 (x,y)  \, \mathbf
e^\dagger  , $$

\noindent   with    $  \beta_1 (x,y)$  , $    \beta_2 (x,y)   \in  \mathbb C ( \mathbf i ) $,  then
$$   \langle  x,y \rangle^\ast    = \overline{  \beta_1
(x,y)}   \, \mathbf e  +  \overline{  \beta_2 (x,y) }  \, \mathbf
e^\dagger   ,  $$

\noindent  and
$$\begin{array}{rcl}
 \pi_{1, \mathbf i }  ( \langle  x,y \rangle^\ast )  &      =     &    \langle  x,y    \rangle^\ast_{ X_{\mathbf e ,  \mathbf i } }   =        \overline{  \beta_1 \,
(x,y) }
\\   &  &   \\   &     =   &
\overline{ \pi_{1 , \mathbf i} \,  ( \langle x,y \rangle )  }  =    \overline{ \langle  x , y  \rangle}_{ X_{\mathbf e ,  \mathbf i } }  \, .
\end{array}$$

\medskip

\item[(d)]  First  note  that  if  $ x, y \in  X_{\mathbf e}$, then   $ x =  \mathbf e \, x$,  $ y =  \mathbf e \, y$,  and  thus
\begin{equation}\label{prod int de eltos en X e}
\begin{array}{l}
\mathbb B \mathbb C \ni \beta_1 (x,y) \, \mathbf e +  \beta_2 (x,y) \, \mathbf e^\dagger  \;   =  \;  \langle x , y \rangle  =   \langle   \mathbf e \, x ,  \mathbf e \, y \rangle  \; =
\\     \\  \qquad \quad  \qquad   = \;
 \mathbf e \, \langle  x , y \rangle =   \mathbf e \, \left(  \beta_1(x,y) \mathbf e +  \beta_2(x,y) \mathbf e^\dagger \right)  \; =
 \\    \\   \qquad  \quad   \qquad    =   \;    \beta_1 (x,y) \, \mathbf e \, ,
\end{array}
\end{equation}

\noindent
hence  $\beta_2 (x,y)   =0$.  In  particular,  taking  $y = x    \in  X_{\mathbf e} $  one  has    $  \langle x,x \rangle  =    \beta_1 (x,x) \, \mathbf e $.   Recall  that   $  \langle x,x \rangle  =0$  if  and  only  if  $x =0$,  that  is,    $x =0$  if  and  only  if   $  \beta_1(x,x) =0$.  On  the  other  hand
$$    \beta_1(x,x)   =  \pi_{1 , \mathbf i }  \,  \left(  \langle  x,x \rangle  \right)  =     \langle  x,x \rangle_{    X_{\mathbf e ,  \mathbf i } }$$

\medskip

\noindent
and  thus  the  last   property  of  the  complex  inner  product  is  proved.

\end{enumerate}

\medskip

Equation (\ref{prod int de eltos en X e})  gives   one   more  reason  to  call  the  $\mathbb C (\mathbf i)$--elements  the  elements  of    $ X_{\mathbf e}$.  \\

Quite  analogously we can show that  formulas    (\ref{prod int e j}),     (\ref{prod int e dagger i}),   (\ref{prod int e dagger j})  define  (complex)     inner  products  on  the  corresponding  linear  spaces.  Thus  we  have  just  proved  the  following

\begin{subtheorem} \label{Teorema sobre restricciones}
Let $X$ be a  $ \mathbb B \mathbb C$-inner  product  space.  Then:

\begin{enumerate}

\item[(a)]  $ \left(  X_{\mathbf e}, \langle \cdot , \cdot  \rangle_{ X_{\mathbf e , \mathbf i } }  \right)  $  and   $ \left(  X_{\mathbf e^\dagger}  , \langle \cdot , \cdot  \rangle_{ X_{\mathbf e^\dagger , \mathbf i } } \right)  $
   are  $\mathbb C( \mathbf i )$-inner  product  spa\-ces.

\item[(b)]  $ \left(  X_{\mathbf e}, \langle \cdot , \cdot  \rangle_{ X_{\mathbf e , \,   \mathbf j }  } \right)  $  and   $ \left(  X_{\mathbf e^\dagger}  , \langle \cdot , \cdot  \rangle_{ X_{\mathbf e^\dagger , \,  \mathbf j } }  \right)  $
   are  $\mathbb C( \mathbf j )$-inner  product  \ spa\-ces.

\end{enumerate}

\end{subtheorem}

\medskip


\begin{subParr}\label{bicomplex inner producet and norm}  We will now   show that the results from  Section    \ref{restrictions of inner product} , are fully consistent with those in Subsections  \ref{inner product on the sum of two complex spaces}    and  \ref{bicomplex norm from inner product}.  Indeed, since
$$  X  =  \e  X  + \edag X  =  X_\e  + X_\edag    =  \e X_\e  + \edag X_\edag  =: \e X_1 + \edag X_2  \, , $$
\end{subParr}

\noindent with  $X_1  = \e X = X_\e$,   $X_2 = \edag X = X_\edag$, we can  consider the bicomplex module  $X$  as  generated by  the   $ \C ( \i )$--linear  spaces   $X_1 $  and  $X_2$.  Moreover  for  $x , y  \in X$  one has:
$$  x  =   \e x +  \edag  x  =  \e  ( \e  x )  + \edag (\edag x)  =:  \e x_1  + \edag x_2 \in  \e X_1  + \edag X_2  , $$
 $$  y  =   \e y +  \edag  y  =  \e  ( \e  y )  + \edag (\edag y)  =:  \e y_1  + \edag y_2 \in  \e X_1  + \edag X_2 $$

where  $x_1 , y_1  \in X_1$  and  $x_2 , y_2 \in X_2$.

\medskip

By Theorem \ref{Teorema sobre restricciones},    we have that    $ \displaystyle  \left( X_1 , \langle \cdot , \cdot \rangle_1 \right)$  and   $ \displaystyle  \left( X_2 , \langle \cdot , \cdot \rangle_2 \right)$,  with   $  \langle \cdot , \cdot \rangle_1 $     $  :=       \langle \cdot , \cdot \rangle_{  X_{\e  , \i } } $  and    $  \langle \cdot , \cdot \rangle_2   :=       \langle \cdot , \cdot \rangle_{  X_ {\edag  ,    \i}  } $,   are  $ \C ( \i )$--inner  product  spaces.    Thus,  the  results  in   Subsection  \ref{inner product on the sum of two complex spaces}  imply   that $X$ can be endowed with the bicomplex inner product
$$ \begin{array}{rcl}
\langle x , y \rangle_X & :=   & \e \langle x_1 , y_1 \rangle_1 + \edag \langle x_2 , y_2 \rangle_2
\\  &  &   \\  & =    &
\e  \langle \e x , \e y \rangle_{X_{ \e , \i } } +   \edag  \langle \edag   x , \edag y \rangle_{X_{ \edag , \i } }  \, .
\end{array}$$

It turns out that
$$ \langle x,y \rangle_X  =  \langle x , y \rangle \,  . $$

Indeed,  using   (\ref{eqU})  and (\ref{eqV})
$$   \langle x,y \rangle_X  =   \e  \beta_1 (x,y)  + \edag \beta_2  (x,y) =  \langle x , y \rangle \,  . $$

\medskip

We summarize this analysis in the next theorem.

\medskip

\begin{subtheorem}\label{T2521}   Given a bicomplex module $X$, then
\end{subtheorem}
\begin{enumerate}
\item[1)] $X$ has a bicomplex inner  product $ \langle \cdot , \cdot \rangle$  if   and  only  if  both   $X_{\e ,\i }  $    and  $X_{\edag, \i } $  have  complex  inner  products    $ \langle \cdot , \cdot \rangle_{X_{\e , \i} }  $   and     $ \langle \cdot , \cdot \rangle_{X_{\edag , \i} }  $;  if  this  holds  then  the  complex  inner  products  are  the  components   of  the   idempotent  representation  of  the  bicomplex  inner  product, i.e.,
\begin{equation}\label{prod int inducido por descomp idemp}
 \langle x , y \rangle  =  \e \langle \e x , \e y \rangle_{X_{\e , \i} }   +     \edag \langle \edag  x , \edag  y \rangle_{X_{\edag  , \i} }  \, .
 \end{equation}

\medskip

\item[2)]  If $X$  has  a  bicomplex  inner  product  then  it  has  a  real--valued   Euclidean--type      norm    and  a  hyperbolic  norm  which  are  related  with  the  norms  on   $ X_{\e , \i} $  and $X_{\edag , \i} $  by the formulas
\begin{equation}\label{norma bicompleja}
\displaystyle  \| x \|_X  :=   \frac{1 \, }{ \sqrt{ 2 } }  \sqrt{  \, \| \e x \|_{  X_{ \e , \i }}^2   +    \| \edag x \|_{  X_{ \edag , \i }}^2  \, }   \, .
\end{equation}

and
\begin{equation}
\| x  \|_\D  :=    \| \e x \|_{  X_{ \e , \i }}   \cdot \e  +    \| \edag x \|_{  X_{ \edag , \i }}    \cdot \edag \, .
\end{equation}

\medskip

\item[3)]  $ \displaystyle \left(  X , \langle \cdot , \cdot   \rangle   \right) $  \,   is    \, a   \,    bicomplex  \,    Hilbert  space  if  and  only  if   $ \displaystyle   \left(  X_{\e , \i }   , \langle \cdot , \cdot   \rangle_{ X_{ \e , \i } }    \right) $   and    $ \displaystyle   \left(  X_{\edag , \i }   , \langle \cdot , \cdot   \rangle_{ X_{ \edag , \i } }    \right) $  are  complex  Hilbert  spaces.

\end{enumerate}

\medskip

The  previous  Theorem  deals  with arbitrary  bicomplex  modules on which we have made no  additional assumptions; the results  explains the  relations  between the  inner   products which can be defined.   In  Section  \ref{cartesian decomp}, on the other hand,   we  developed  the theory of   bicomplex   inner  products  in the case in which a  bicomplex  module  $X$  has   a  dagger--    or   a   bar--involution.  We  will now  show  how  these situations   interact.

\bigskip

\begin{subtheorem}  Let  $X$ be a $\bc$--module  with  a  $\dagger$--involution.  Then  the following  three  statements are   equivalent:

\begin{enumerate}

\item[(1)]    $X$ has  a  bicomplex  inner  product.

\medskip

\item[(2)]   Each of  the  $\C (\i)$--linear  spaces  $X_\e$  and  $X_\edag$  has  a  $\C (\i)$--valued  inner  product,  $ \langle \cdot , \cdot \rangle_{X_\e , \i } $  and
$ \langle \cdot , \cdot \rangle_{X_\edag , \i } $  respectively.

\medskip

\item[(3)]  The $\C (\i)$--linear  space   $X_{1, \dagger}$  has  a   $\C (\i)$-valued  inner  product  $ \langle \cdot , \cdot \rangle_{X_1 , \dagger } $.

\end{enumerate}

\end{subtheorem}

{\sf Proof:}  \\       We already proved that   (1)  is equivalent  to   (2)  and   that  (1)  implies   (3).

\medskip

Assume now  that   $X_{1 , \dagger}$  has  a   $ \C ( \i ) $--valued  inner  product     $ \langle \cdot , \cdot   \rangle_{ \C ( \i ) } $.  Take  $x = s_1 + \j \, s_2$  and   $y =  t_1 + \j  \, t_2 $ in $X$. By (\ref{prod int inducido en X desde X1dagger}) we have that formula
$$ \begin{array}{rcl}
\langle x , y \rangle  &  :=  &   \displaystyle \left(  \langle s_1 + \i \, s_2 , t_1 \rangle_{ \C ( \i ) }   +    \langle s_2  -  \i \, s_1 , t_2 \rangle_{ \C ( \i ) }   \right) \, \e  \; +
\\   &   &   \\   &   &
\qquad  \quad
\displaystyle   + \;    \left(  \langle s_1 - \i \, s_2 , t_1 \rangle_{ \C ( \i ) }   +    \langle s_2  +  \i \, s_1 , t_2 \rangle_{ \C ( \i ) }   \right) \, \edag
\end{array} $$

gives  a  bicomplex   inner  product  on  $X$;  thus each  of  its  complex  components   defines  an  inner  product  on  the  corresponding  idempotent  components  of  $X$:
$$  \langle \e x , \e y \rangle_{X_{\e , \i } }  :=    \langle s_1 + \i \, s_2 , t_1 \rangle_{ \C ( \i ) }   +    \langle s_2  -  \i \, s_1 , t_2 \rangle_{ \C ( \i ) }  \, , $$
$$  \langle \edag x , \edag  y \rangle_{X_{\edag , \i } }  :=    \langle s_1 - \i \, s_2 , t_1 \rangle_{ \C ( \i ) }   +    \langle s_2  +  \i \, s_1 , t_2 \rangle_{ \C ( \i ) }  \, ;$$

this shows that  (3)  implies   (1)  and  (2). \mbox{} \qed \mbox{} \\

It  is  obvious  that when  $X$  has  a  bar--involution  then    conditions  (2)  and  (3)  in  the  above  Theorem will require $\C (\j)$--linearity  instead  of $\C (\i )$--linearity.

\bigskip


\section{Complex  inner  products  on  $X$  induced  by
 idempotent  decompositions.}      \label{complex inner product}

It is possible to give a different  approach to the notion  of  bicomplex Hilbert  space, by involving only complex structures. We briefly study this approach  below.

\medskip

We can consider $X$  as  $X_{\mathbb C (\mathbf i )}$  and    $X_{\mathbb C (\mathbf j )}$; from this point of view, $X$ is  endowed     with a specific  inner   product   if  this is true for  both  $ X_{\mathbf e , \mathbf i }  $  and    $ X_{\mathbf e^\dagger , \mathbf i}$,  as  well  as  both       $ X_{\mathbf e , \,  \mathbf j }  $  and    $ X_{\mathbf e^\dagger , \,    \mathbf j}$,  where  the  subindexes   $\mathbf i$  and $\mathbf j$  mean $ \mathbb C ( \mathbf i )$  or $   \mathbb C ( \mathbf j )$   linearities.      Indeed,  since
$$  X_{ \mathbb C ( \mathbf i ) }   =  X_{\mathbf e , \mathbf i}    \oplus_{ \mathbb C ( \mathbf i ) }   X_{\mathbf e^\dagger , \mathbf i}  $$

and
$$  X_{ \mathbb C ( \mathbf j ) }   =  X_{\mathbf e , \,  \mathbf j}    \oplus_{ \mathbb C ( \mathbf j ) }   X_{\mathbf e^\dagger ,  \,   \mathbf j}   \, ,  $$

\medskip

we  can  endow  $  X_{ \mathbb C ( \mathbf i ) }  $  and  $  X_{ \mathbb C ( \mathbf j ) } $   with  the  following  inner  products  respectively:  given   $ x, y \in   X_{ \mathbb C ( \mathbf i ) }  $  or  in  $   X_{ \mathbb C ( \mathbf j ) }  $  set
\begin{equation}\label{prod int complejo i en X}
 \begin{array}{rcl}
\langle x,y \rangle_{ \mathbb C ( \mathbf i ) }  &  :=  &   \langle   P (x) + Q (x) , P (y) + Q (y)  \rangle_{\mathbb C ( \mathbf i ) }
\\ &  &  \\   &  : = &
\langle   P(x) , P(y)  \rangle_{ X_{ {\mathbf e} , \mathbf i } }    +   \langle Q (x)  , Q (y)  \rangle_{ X_{ {\mathbf e^\dagger} , \mathbf i  } } \;    \in  \mathbb C ( \mathbf i )  \, ,
\end{array}
\end{equation}

and
\begin{equation}\label{prod int complejo j en X}
\begin{array}{rcl}
\langle x,y \rangle_{ \mathbb C ( \mathbf j ) }  &  :=  &   \langle   P (x) + Q (x) , P (y) + Q (y)  \rangle_{\mathbb C ( \mathbf j ) }
\\ &  &  \\   &  : = &
\langle   P(x) , P(y)  \rangle_{ X_{ {\mathbf e} , \,  \mathbf j }   } +   \langle Q (x)  , Q (y)  \rangle_{ X_{ {  \mathbf e^\dagger} ,   \, \mathbf j }    } \;  \in  \mathbb C (\mathbf j )   \, .
\end{array}
\end{equation}

Of  course,  if  both  $ (  X_{\mathbf e}, \langle \cdot , \cdot  \rangle_{ X_{  {\mathbf e} , \mathbf i } } )  $  and   $ (  X_{\mathbf e^\dagger}  , \langle \cdot , \cdot  \rangle_{ X_{ {\mathbf e^\dagger} , \mathbf i } } )  $  are  Hilbert  spaces,  then   $ (  X_{ \mathbb C (\mathbf i )}, \langle \cdot , \cdot  \rangle_{ \mathbb C (\mathbf i)}   )  $  is  also  a  Hilbert  space.  The  same is true in the    $\mathbb C ( \mathbf j )$ case.

\bigskip


Assume  now  that  the  bicomplex  module  $X$  has  a   bicomplex  inner  \,  product   \, $ \langle \cdot , \cdot \rangle $.   \,  By  Theorem     \ref{Teorema sobre restricciones},     $ (  X_{\mathbf e}, \langle \cdot , \cdot  \rangle_{ X_{ {\mathbf e} , \mathbf i } }  )  $  and   $ (  X_{\mathbf e^\dagger}  , \langle \cdot , \cdot  \rangle_{ X_{ {\mathbf e^\dagger} , \mathbf i }  } )  $
   are  $\mathbb C( \mathbf i )$-inner  product  spaces,   and     $(  X_{\mathbf e}, \langle \cdot , \cdot  \rangle_{ X_{ {\mathbf e} , \,   \mathbf j } }  )  $  and   $ (  X_{\mathbf e^\dagger}  , \langle \cdot , \cdot  \rangle_{ X_{  {\mathbf e^\dagger} , \,  \mathbf j } }  )  $
   are  $\mathbb C( \mathbf j )$-inner  product  spaces,  hence       $X_{\mathbb C ( \mathbf i )} $  and   $X_{\mathbb C ( \mathbf j ) }$
are  (complex)  linear  spaces  with  inner  products    $ \langle \cdot , \cdot  \rangle_{\mathbb C ( \mathbf i ) }$  and   $ \langle \cdot , \cdot  \rangle_{\mathbb C ( \mathbf j ) } $
determined  by  equations   (\ref{prod int complejo i en X})  and      (\ref{prod int complejo j en X});   moreover,       $X_{\mathbb C ( \mathbf i )} $  and   $X_{\mathbb C ( \mathbf j ) }$  become  (complex)  normed  spaces    with  the  norms  determined  in the usual  way  by  the  formulas:
\begin{equation} \label{norma i en X}
 \begin{array}{rcl}
\| x \|_{\mathbb C ( \mathbf i ) }  &  :=   &  \langle x,x \rangle_{\mathbb C ( \mathbf i ) }^{1/2}
\\  &  &  \\   & := &
\displaystyle    \left(   \langle  P(x) , P(x) \rangle_{X_{ {\mathbf e } , \mathbf i } }   +    \langle  Q(x) , Q(x) \rangle_{X_{ {\mathbf e^\dagger } , \mathbf i }  }  \right)^{1/2}
\\  &  &  \\   & = &
\displaystyle  \left(  \| P(x) \|_{X_{ { \mathbf e } , \mathbf i } }^2  +    \| Q(x) \|_{X_{ {\mathbf e^\dagger  , \mathbf i } } }^2  \right)^{1/2}
\end{array}
\end{equation}

and  similarly
\begin{equation}\label{norma j en X}
 \begin{array}{rcl}
\| x \|_{\mathbb C ( \mathbf j ) }  &  :=   &  \langle x,x \rangle_{\mathbb C ( \mathbf j ) }^{1/2}
\\  &  &  \\   & := &
\displaystyle    \left(   \langle  P(x) , P(x) \rangle_{X_{\mathbf e  ,  \,   \mathbf j } }   +    \langle  Q(x) , Q(x) \rangle_{X_{ \mathbf e^\dagger  ,  \,    \mathbf j } }  \right)^{1/2}
\\  &  &  \\   & = &
\displaystyle  \left(  \| P(x) \|_{X_{ \mathbf e  ,  \,  \mathbf j }  }^2  +    \| Q(x) \|_{X_{\mathbf e^\dagger  ,   \,    \mathbf j } }^2  \right)^{1/2}   \, .
\end{array}
\end{equation}

\medskip

Recalling  that
$$ \begin{array}{rcl}
\mathbb D^+ \ni \langle x,x \rangle & = &  \langle P(x) , P(x) \rangle_{X_{\mathbf e  , \mathbf i } }  \, \mathbf e  +     \langle Q(x) , Q(x) \rangle_{X_{\mathbf e^\dagger  , \mathbf i }}  \, \mathbf e^\dagger
\\  &  &  \\  & = &
\langle P(x) , P(x) \rangle_{X_{\mathbf e  , \, \mathbf j }  }  \, \mathbf e  +     \langle Q(x) , Q(x) \rangle_{X_{\mathbf e^\dagger  , \,  \mathbf j } }  \, \mathbf e^\dagger \, ,
\end{array} $$

one  sees  that
$$   \langle P(x) , P(x) \rangle_{X_{\mathbf e  , \mathbf i } }    =    \langle P(x) , P(x) \rangle_{X_{\mathbf e  ,  \, \mathbf j } }    \in  \mathbb R^+ \cup \{ 0 \}  $$

and
$$    \langle Q(x) , Q(x) \rangle_{X_{\mathbf e^\dagger  , \mathbf i }  }   =       \langle Q(x) , Q(x) \rangle_{X_{\mathbf e^\dagger  , \,   \mathbf j }   }     \in  \mathbb R^+ \cup \{ 0 \}  \, ,   $$

which implies
\begin{equation}\label{igualdad normas complejas}
\|  x \|_{ \mathbb C ( \mathbf i ) }   =    \|  x \|_{ \mathbb C ( \mathbf j ) }    \quad \forall x \in X  \, .
\end{equation}

\medskip

Formulas  (\ref{estrella}),   (\ref{eq A}),  and  the    fact  that  we  can  take  the  square  root  of  a  bicomplex  number  by  taking  the  square root  of  each  of its idempotent  components,    imply   that
$$  \sqrt{   \langle x ,x \rangle }  =   \|  P (x)  \|_{  X_{ \mathbf e , \mathbf i } }   \cdot \mathbf e  +    \|  Q (x)  \|_{  X_{ \mathbf e^\dagger , \mathbf i } }   \cdot \mathbf e^\dagger   $$

and  therefore
\begin{equation}\label{normas complejas y prod int bicomplejo}
\sqrt{2} \,    | \sqrt{   \langle x ,x \rangle \, } |   =     \|  x \|_{ \mathbb C ( \mathbf i ) }   =    \|  x \|_{ \mathbb C ( \mathbf j ) }     \quad \forall x \in X  \,
\end{equation}

which shows the  relation  between  both   norms  and  the original   bicomplex  inner  product.  \\

Equation   (\ref{igualdad normas complejas})  tells  us  that  although  we  have  on  $X$  only  one  map
$$ x \in X  \mapsto  \|  x \|_{ \mathbb C ( \mathbf i ) }   =    \|  x \|_{ \mathbb C ( \mathbf j ) }  $$

such map  determines  two  \underline{different}  norms on $X$  seen as  $X_{\mathbb C ( \mathbf i )} $  and   $X_{\mathbb C ( \mathbf j ) }$; in particular,  $X$ can be  endowed  with two  linear  topologies:  one  with scalars  in
$ \mathbb C ( \mathbf i ) $  and another one  with  scalars  in   $ \mathbb C ( \mathbf j ) $. When we do not consider the field in which we take the scalars, we have the same  underlying  topology, which we will denote  by  $\tau_{ \mathbb C }$; note that this topology is   related to the  initial  bicomplex   inner  product     via   (\ref{normas complejas y prod int bicomplejo}).

Note also that  $ \displaystyle  |  \sqrt{ \,  \langle x , x \rangle \, } |$  coincides     with the norm   (\ref{norma bicompleja}).     The  approach  given  in  Section   \ref{bicomplex inner producet and norm}    allows  us   to  deal  with  the  bicomplex  norm  on  the  bicomplex  module  $X$.  What  is  more, the  topology  $\tau_\C$  is  in  fact  a  bicomplex  linear  topology, that  is,  the  multiplication  by  bicomplex  scalars  (not  only  by   $\C ( \i)$  scalars)  is  a  continuous  operation.  \\

We conclude this section by illustrating  some  aspects  of  this     (complex)  approach  to  the  notion  of  inner  products  on  bicomplex  modules  using  again  $X = \bc$.  \\

Note   that  formulas   (\ref{repr_idempotent})  and  (\ref{repr_idempotent_Cj})  give  the  following  representations of a  given  bicomplex  number:
$$ \begin{array}{rcl}
Z & = & (x_1 + \mathbf i \, y_1 ) + (x_2 + \mathbf i \, y_2 ) \, \mathbf j
\\  &  &  \\  & = &
(x_1 + \mathbf j \, x_2 ) + (y_1 + \mathbf j \, y_2 ) \, \mathbf i
\:  =  \;  \beta_1 \, \mathbf e  +  \beta_2 \, \mathbf e^\dagger
\\   &  &  \\    &  :=  &
\left(   (x_1 + y_2 ) + \mathbf i \,  (y_1 - x_2) \right) \, \mathbf e  +   \left(  (x_1 - y_2 ) + \mathbf i \, (y_1 + x_2 ) \right) \, \mathbf e^\dagger
\\  &  &   \\   & = &
\gamma_1 \, \mathbf e  + \gamma_2 \, \mathbf e^\dagger
\\   &   &   \\  & :  = &
\left(   (x_1 + y_2 ) + \mathbf j \,  (x_2 - y_1) \right) \, \mathbf e  +   \left(  (x_1 - y_2 ) + \mathbf j \, (y_1 + x_2 ) \right) \, \mathbf e^\dagger \, .
\end{array} $$

Consider a second bicomplex number
$$ \begin{array}{rcl}
W & = &  (s_1 + \mathbf i \, t_1) + (s_2 + \mathbf i \, t_2 ) \, \mathbf j \\  &  &  \\  & = &
\xi_1 \, \mathbf e + \xi_2 \, \mathbf e^\dagger    \\  &  &  \\  & := &
\left(   (s_1 + t_2) + \mathbf i \, (t_1 - s_2) \right) \, \mathbf e  +  \left(  ( s_1 - t_2 )  + \mathbf i \, (t_1 + s_2 ) \right) \, \mathbf e^\dagger   \\  &  &  \\  & = &
\eta_1 \, \mathbf e + \eta_2 \, \mathbf e^\dagger   \\  &  &  \\  & :  = &
\left(   (s_1 + t_2) + \mathbf j \, ( s_2  - t_1) \right) \, \mathbf e  +  \left(  ( s_1 - t_2 )  + \mathbf j \, (t_1 + s_2 ) \right) \, \mathbf e^\dagger  \, .
\end{array}   $$

In  terms  of  the  previous  sections, one  has:  \\

$X = \mathbb B  \mathbb C$;  \;   $X_{\mathbf e , \mathbf i } =  \mathbb B  \mathbb C \cdot \mathbf e   =   \mathbb C ( \mathbf i )  \cdot \mathbf e    $;  \;  $   X_{\mathbf e^\dagger , \mathbf i } =  \mathbb B  \mathbb C \cdot \mathbf e^\dagger   =   \mathbb C ( \mathbf i )  \cdot \mathbf e^\dagger    $;

$$   \langle   \beta_1 \, \mathbf e  , \xi_1 \, \mathbf e \rangle_{  \mathbb C ( \mathbf i )  \cdot \mathbf e }  : =  \beta_1 \,   \overline \xi_1 \, ; $$

$$   \langle   \beta_2 \, \mathbf e^\dagger  , \xi_2 \, \mathbf e^\dagger \rangle_{  \mathbb C ( \mathbf i )  \cdot \mathbf e^\dagger }  : = \beta_2 \,  \overline \xi_2 \, ; $$

\medskip

 $X_{\mathbf e , \mathbf j } =  \mathbb B  \mathbb C \cdot \mathbf e   =   \mathbb C ( \mathbf j )  \cdot \mathbf e    $;  \;  $   X_{\mathbf e^\dagger , \mathbf j } =  \mathbb B  \mathbb C \cdot \mathbf e^\dagger   =   \mathbb C ( \mathbf j )  \cdot \mathbf e^\dagger    $;

$$   \langle   \gamma_1 \, \mathbf e  , \eta_1 \, \mathbf e \rangle_{  \mathbb C ( \mathbf j )  \cdot \mathbf e }  : = \gamma_1 \, \eta_1^\ast \, ; $$

$$   \langle   \gamma_2 \, \mathbf e^\dagger  , \eta_2 \, \mathbf e^\dagger \rangle_{  \mathbb C ( \mathbf j )  \cdot \mathbf e^\dagger }  : = \gamma_2 \, \eta_2 ^\ast   \, . $$

Then,  formulas   (\ref{prod int complejo i en X})   and  (\ref{prod int complejo j en X})    give:
$$  \begin{array}{l}
\langle Z , W \rangle_{ \mathbb C ( \mathbf i ) }   \;   :=
\\   \\
\quad = \;   \displaystyle   \left\langle    \left(   (x_1  + y_2 ) + \mathbf i \,  ( y_1 - x_2) \right) \, \mathbf e \, , \,    \left(   (s_1  + t_2 ) + \mathbf i \,  ( t_1 - s_2) \right) \, \mathbf e \,  \right\rangle_{ \mathbb C ( \mathbf i ) \, \mathbf e }  \;  +
\\    \\
\quad   + \;
\displaystyle   \left\langle    \left(   (x_1  - y_2 ) + \mathbf i \,  ( y_1 + x_2) \right) \, \mathbf e^\dagger \, , \,    \left(   (s_1  - t_2 ) + \mathbf i \,  ( t_1 + s_2) \right) \, \mathbf e^\dagger \,  \right\rangle_{ \mathbb C ( \mathbf i ) \, \mathbf e^\dagger }
\\   \\
\quad = \;   \displaystyle   \left(  ( x_1 + y_2 ) \, ( s_1 + t_2 )   + ( y_1 - x_2 ) \, ( t_1 - s_2 )    +  \right.
\\  \\
\qquad   \qquad  + \;     \displaystyle   \left.     (x_1 - y_2 ) \, (s_1 - t_2)   \; +      (y_1  + x_2 ) \, ( t_1 + s_2 )  \right)    \; + \;
\\  \\
\qquad  \qquad   \quad   + \; \displaystyle      \mathbf i \,   \left(    (y_1 - x_2 ) \, ( s_1 + t_2 )   \; -    (x_1 + y_2 ) \, (  t_1 - s_2 )   \; +   \right.
\\  \\
\qquad  \qquad \quad   \quad   \displaystyle   + \;  \left.   (y_1 + x_2 ) \, ( s_1  - t_2)     -   (x_1 - y_2 ) \, ( t_1    + s_2)    \right)  \in  \mathbb C ( \mathbf i ) \, ;
\end{array} $$

$$  \begin{array}{l}
\langle Z , W \rangle_{ \mathbb C ( \mathbf j ) }   \;   :=
\\   \\
\quad = \;   \displaystyle   \left\langle    \left(   (x_1  + y_2 ) + \mathbf j \,  ( x_2  -  y_1  ) \right) \, \mathbf e \, , \,    \left(   (s_1  + t_2 ) + \mathbf j \,  (  s_2  -  t_1  ) \right) \, \mathbf e \,  \right\rangle_{ \mathbb C ( \mathbf j ) \, \mathbf e }  \;  +
\\    \\
\quad   + \;
\displaystyle   \left\langle    \left(   (x_1  - y_2 ) + \mathbf j \,  ( y_1 + x_2) \right) \, \mathbf e^\dagger \, , \,    \left(   (s_1  - t_2 ) + \mathbf j \,  ( t_1 + s_2) \right) \, \mathbf e^\dagger \,  \right\rangle_{ \mathbb C ( \mathbf j ) \, \mathbf e^\dagger }
\\   \\
\quad = \;   \displaystyle   \left(  ( x_1 + y_2 ) \, ( s_1 + t_2 )   + (  x_2  - y_1   ) \, (  s_2  - t_1   )    +  \right.
\\  \\
\qquad   \qquad  + \;     \displaystyle   \left.     (x_1 - y_2 ) \, (s_1 - t_2)   \; +      (y_1  + x_2 ) \, ( t_1 + s_2 )  \right)    \; + \;
\\  \\
\qquad  \qquad   \quad   + \; \displaystyle      \mathbf j \,   \left(  ( x_2  -  y_1  ) \, ( s_1 + t_2 )   -   (x_1 + y_2 ) \, (   s_2  - t_1  )   \; +   \right.
\\  \\
\qquad  \qquad \quad   \quad   \displaystyle   + \;  \left.   (y_1 + x_2 ) \, ( s_1  - t_2)     -   (x_1 - y_2 ) \, ( t_1    + s_2)    \right)  \in  \mathbb C ( \mathbf j ) \, .
\end{array} $$

Of  course,    $\langle Z,W \rangle_{ \mathbb C ( \mathbf i ) } \in    \mathbb C ( \mathbf i ) $  and    $\langle Z,W \rangle_{ \mathbb C ( \mathbf j ) } \in    \mathbb C ( \mathbf j ) $   are  different  bicomplex  numbers;  but  their  ``real  parts"  coincide, while their  ``imaginary  parts"     are  in general different. They coincide, for example, when $Z=W.$, in which case
$$ \langle  Z ,Z \rangle_{ \mathbb C ( \mathbf i ) }  =      \langle  Z ,Z \rangle_{ \mathbb C ( \mathbf j ) }    \geq  0 \, . $$

\bigskip

\section{The bicomplex  module $\bc^n$.}


\begin{subParr}\label{inner_prod_in_bc_n}     We already know that   $\bc^n$  is  a  bicomplex  module.  A  bicomplex  inner  product  on it  can  be  introduced  in  a  canonical  way:  if  $Z=  ( Z_1 , \ldots , Z_n)$,  and $W = (W_1 , \ldots , $     $   W_n)  \in  \bc^n$   then  we set
$$ \langle Z \, , \,  W \rangle :=  Z_1 W_1^\ast  + \cdots  +  Z_n \, W_n^\ast .$$
\end{subParr}

It  is  obvious  that  properties  (1), (2), (3)  of  Definition  \ref{def inner product}  are  satisfied.  The  last  property  can be proved  similarly to  the  proof  in  Example \ref{Ejemplo caso bicomplejo},  but   we  prefer  to  give  an  alternative  proof  using  the  idempotent  representation.  Indeed,  write $Z_k = \beta_{k , \, 1} \, \e + \beta_{k , \, 2} \, \edag $, and correspondingly  $Z_k^\ast   = \overline{ \beta}_{k , \, 1} \, \e + \overline{ \beta}_{k , \, 2} \, \edag $. Then we have
$$ \begin{array}{rcl}
\langle Z , \, Z \rangle & = &  \displaystyle \sum_{k=1}^n Z_k \cdot  Z_k^\ast  =  \sum_{k=1}^n  \left( | \beta_{k, \, 1} |^2 \, \e  + | \beta_{k , \, 2 } |^2 \, \edag \right)
\\  &  &  \\  & = &
\displaystyle  \left(  \,  \sum_{k=1}^n  | \beta_{k , \, 1 } |^2 \right) \, \e   +     \left(  \,   \sum_{k=1}^n  | \beta_{k , \, 2 } |^2 \right) \, \edag  \in \mathbb D^+  \, .
\end{array}  $$

Note  that  the square of the inner product of some   $Z \in \bc^n$ can  be  a  zero  divisor:  this  happens
exactly  when  $Z \in \bc^n_\e = \bc^n \cdot \e $  or   $Z \in \bc^n_\edag = \bc^n \cdot \edag $.   Note  also  that  this  formula  defines  the  hyperbolic   norm  on  the  $\bc$--module  $\bc^n$.

\begin{subParr}  Similarly  to what happens in the  complex  case,  one can consider  more  ge\-ne\-ral  inner  products on  $\bc^n$. To do so, take  any  bicomplex  positive  matrix   $A = \mathscr A_1 \, \e + \mathscr A_2 \, \edag$  and define,  for   $Z, \; W \in \bc^n$,
$$  \langle Z , \, W \rangle_A :=  Z^t \cdot A \cdot W^{ \ast \,  } \, .$$
\end{subParr}

We can easily show  that this is  an  inner  product in the sense  of  Definition \ref{def inner product}:

\begin{enumerate}

\item[(1)]    $  \langle Z , \, W+V \rangle_A \;  = \;   Z^t \cdot A \cdot ( W+V )^{\ast  }  =  Z^t \cdot A \cdot  \left(  W^{ \ast }  +  V^{ \ast } \right)  $
\\  \\
$ \mbox{ }   \qquad  \qquad   = \;  Z^t   \cdot A \cdot   W^{ \ast }  + Z^t \cdot A \cdot     V^{ \ast }   $
\\  \\
$  \mbox{}  \qquad  \qquad  =  \;    \langle Z , \, W \rangle_A  + \langle Z , \, V \rangle_A .  $

\medskip

\item[(2)]   Given  $\mu \in \bc$,   then
$$ \begin{array}{rcl}
 \langle \mu \, Z , \, W \rangle_A  & = &  \displaystyle \left( \,  \mu \, Z \right)^t  \cdot A \cdot W^{\ast  }  =  \mu \, \left(  Z^t  \cdot A \cdot W^{\ast } \right)
 \\  &  &   \\   & = &  \mu \, \langle Z , \, W \rangle_A .
 \end{array} $$

 \medskip

 \item[(3)]  $   \displaystyle   \langle Z , \, W \rangle_A =  Z^t   \cdot A  \cdot   W^{\ast  }  =  \left(  \, Z^t   \cdot A \cdot   W^{\ast  }  \right)^t  $
 \\  \\
$ \mbox{}  \qquad \qquad  \displaystyle =  \;   \left(   \left(   W^{\ast \, t }   \cdot A^t \cdot Z \right)^\ast \, \right)^\ast  $
\\   \\
$ \mbox{}  \qquad \qquad  \displaystyle =  \;  \left(  W^t  \cdot A \cdot Z^{\ast  }  \, \right)^\ast  \; = \;  \langle W , \, Z \rangle_A^\ast \, , $

\medskip

\noindent   where  we took  into account  that  $ Z^t \cdot A \cdot   W^{\ast } $  is  a  $1 \times 1$    matrix.

\medskip

\item[(4)]   Since $A$  is  a  bicomplex  positive  matrix,  for  any  column     $Z \in \bc^n$ we have  that
$$ Z^t   \cdot A \cdot  Z^{ \ast  } \in \mathbb D^+ \, .$$

\end{enumerate}


    Let us now reinterpret  formula  (\ref{prod int inducido por descomp idemp}) in this setting.  The  simplest  way to do this is via a direct computation.  Take  $Z= \mathscr Z_1 \, \e  + \mathscr Z_2 \, \edag$,  $W= \mathscr W_1 \, \e  + \mathscr W_2 \, \edag  \in \bc^n$   and  let     $A = \mathscr A_1 \, \e  + \mathscr A_2 \, \edag$ be  a  
  bicomplex  positive  matrix.   Then
$$  \begin{array}{rcl}
\langle Z , \, W \rangle_A & = &  \displaystyle \left( \mathscr Z_1 \, \e  + \mathscr Z_2 \, \edag   \right)^t \cdot \left(  \mathscr A_1 \, \e  + \mathscr A_2 \, \edag \right) \cdot \left(
\mathscr W_1 \, \e  + \mathscr W_2 \, \edag \right)^{\ast  }
\\  &  &  \\  & =  &
\mathscr Z_1^t \cdot \mathscr A_1 \cdot  \overline{ \mathscr W}_1  \, \e   +  \mathscr Z_2^t   \cdot \mathscr A_2 \cdot  \overline{ \mathscr W}_2  \, \edag
\\   &  &  \\  & = :&
\langle \mathscr Z_1 , \, \mathscr W_1 \rangle_{ \mathscr A_1 } \cdot \e  +   \langle \mathscr Z_2 , \, \mathscr W_2 \rangle_{ \mathscr A_2 } \cdot \edag \, ,
\end{array}  $$

which should be compared  with  formula  (\ref{prod int inducido por descomp idemp}).  \\

The  complex  $( \C (\i) ) $  matrices  $\mathscr A_1$  and $\mathscr A_2$  are  positive  and  thus  $ \langle \cdot , \, \cdot \rangle_{\mathscr A_1 }$  and  $ \langle \cdot , \, \cdot \rangle_{\mathscr A_2 }$   are  $( \C (\i ) )$  complex  inner  products  on  $ \C^n (\i)$.  In particular,
$$ \begin{array}{rcl}
\langle Z , \, Z \rangle_A & = & \langle \mathscr Z_1 , \,  \mathscr Z_1   \rangle_{\mathscr A_1 }  \cdot \e  +      \langle \mathscr Z_2 , \,  \mathscr Z_2   \rangle_{\mathscr A_2 }  \cdot \edag
\\  &  &   \\  & = &
\displaystyle  \| \mathscr Z_1 \|_{ \C^n (\i ) , \, \mathscr A_1 } ^2 \cdot \e  +    \| \mathscr Z_2 \|_{ \C^n (\i ) , \, \mathscr A_2 } ^2 \cdot \edag   \, ,
\end{array} $$

which is
the  square  of  the  hyperbolic    norm  on $\bc^n$.
The  bicomplex  norm on  $\bc^n$  generated by the inner  product  $ \langle \cdot , \, \cdot \rangle_A$  is  given  by
$$  \displaystyle  \| Z \|_A := \frac{1}{ \sqrt{2 } }  \,  \sqrt{  \,    \| \mathscr Z_1 \|_{ \C^n (\i ) , \, \mathscr A_1 } ^2  \;   +  \;      \| \mathscr Z_2 \|_{ \C^n (\i ) , \, \mathscr A_2 } ^2  \, }  \, . $$

\bigskip


\section{The ring   $\H (\C)$ of biquaternions  as a $\bc$--module}


Let us show now how this theory can clarify some  algebraic  properties  of the ring of  biquaternions (sometimes  called   complex  quaternions).  First of all  recall  that  the  ring   $\H (\C)$   of  biquaternions    is
the set of  elements   of the form
$$     \mathscr Z = z_0 +  \i_1 z_1 + \i_2 z_2 + \i_2 z_2 $$

with   $z_0 , \, z_1 , \, z_2 , \,  z_3 \in \C (\i )$;   the three  quaternionic   imaginary units  $\i_1 , \, \i_2 , \, \i_3$ are  such that
$$ \i_1^2=\i_2^2=\i_3^2 = -1 ,   $$

the multiplication between them  is  anti--commutative:
$$ \i_2 \cdot \i_1 = -  \i_1 \cdot \i_2  = - \i_3 ; \quad   \i_3 \cdot \i_2 = - \i_2 \cdot \i_3=- \i_1 ; \quad  $$
$$ \i_1 \cdot \i_3 = - \i_3 \cdot \i_1 = - \i_2, $$

while their  products    with  the  complex    imaginary unit $\i$   are  commutative:
$$ \i \cdot \i_1 = \i_1 \cdot \i ; \quad   \i \cdot \i_2 = \i_2 \cdot \i ; \quad  \i \cdot \i_3 = \i_3 \cdot \i .$$

Since any of the products $ \i \cdot \i_\ell$   satisfies   $( \i \cdot \i_\ell)^2=1$,  there  are  many  hyperbolic  units inside   $\H (\C)$ and thus there  are many subsets  inside it  which  are   isomorphic to   $\bc$, hence there  are  many different ways of  ma\-king    $\H (\C)$ a $\bc$--module.  Let's fix     one of these  ways  and  write   each biquaternion   as
$$ \mathscr Z =  (z_0 +  \i_1 z_1 )  + ( z_2 + \i_1 z_2)  \i_2  =:Z_1 + Z_2 \i_2 , $$

where  $Z_1$,  $Z_2  \in \bc$,  (i.e.,   our   copy  of  $\bc$  has  as      imaginary  units    $\i$,     $\i_1$,   $\k =  \i \i_1$).  It is clear  that   $\H (\C)$   is  a  $\bc$--module,  with the multiplication  by  bicomplex  scalars as follow.  For any    $\Lambda  \in  \bc$  and  any  $\mathscr Z  = Z_1 +Z_2 \i_2   \in    \H (\C)$ one has:
$$  \Lambda \cdot \mathscr Z =  \Lambda \cdot Z_1   +  \Lambda \cdot Z_2 \i_2.   $$

 Given  $\mathscr W = W_1 + W_2 \i_2$ and taking into account   that   for   any $Z \in \bc$ it is  $\i_2 Z = Z^\dagger \i_2$, one has:
$$  \displaystyle    \mathscr Z \cdot \mathscr W   =     \left(   Z_1 W_1 - Z_2 W_2^\dagger    \right)    +    \left(   Z_1 W_2  + Z_2 W_1^\dagger \right) \i_2  .  $$

That is,   $\H (\C)$  becomes    an algebra--type  $\bc$--module.    \\

Because of the numerous imaginary units inside $\H (\C)$,  there  are  many  conjugations, and they can be  expressed using the conjugations  from  $\bc$.   Some of them   are   (there is   no   standard notation in the literature):

$$  \overline{ \mathscr Z} :=   ( \overline{z}_0  +  \overline{z}_1 \i_1 ) +  (  \overline{ z}_2 +   \overline{z}_3 \i_1 ) \i_2  =    \overline{Z}_1  + \overline{Z}_2 \i_2  \, ;$$

$$\mathscr Z^{\dagger_1} :=   (z_0  - z_1 \i_1 ) +  (z_2 - z_3 \i_1 ) \i_2 =  Z_1^\dagger  + Z_2^\dagger \i_2 \,  ;$$

$$  \mathscr Z^\star   :=    (\overline{z}_0  - \overline{z}_1 \i_1 ) +  (\overline{z}_2 - \overline{z}_3  \i_1 ) \i_2   =      Z_1^\ast + Z_2^\ast \i_2  \, ;  $$

$$   \mathscr Z^{\dagger_2}   :=     (z_0  +   z_1 \i_1 ) -    (z_2 +   z_3 \i_1 ) \i_2   =   Z_1  - Z_2 \i_2  \, ;  $$

$$  \mathscr Z^{\dagger_3}:=   (z_0  +  z_1 \i_1 ) +  (z_2 - z_3 \i_1 ) \i_2   =   Z_1  +  Z_2^\dagger \i_2 \, ; $$

$$  \left(  \mathscr Z^{\dagger_1} \right)^{\dagger_2}  =     \left(  \mathscr Z^{\dagger_2} \right)^{\dagger_1}  =  Z_1^\dagger  - Z_2^\dagger \i_2  \, ; $$

$$ \mathscr Z^\odot  :=    Z_1^\ast  - \overline{Z}_2 \i_2 \,   ;$$

$$ \mathscr Z^\diamond  :=    Z_1^\dagger  - Z_2 \i_2 \,   ;$$

etc.  \\

These conjugations interact with the product of two biquaternions as follows:

$$ \overline{ \mathscr Z \cdot \mathscr W }  =   \overline{ \mathscr Z} \cdot  \overline{ \mathscr W}  \, ;  \qquad   (  \mathscr Z \cdot \mathscr W )^{\dagger_1}   =   \mathscr Z^{\dagger_1}   \cdot \mathscr W^{\dagger_1}    \, ;  $$

$$ (  \mathscr Z \cdot \mathscr W )^\star   =   \mathscr Z^\star    \cdot \mathscr W^\star   \, ;  \qquad   (  \mathscr Z \cdot \mathscr W )^{\dagger_2}   =   \mathscr Z^{\dagger_2}   \cdot \mathscr W^{\dagger_2}    \, ;  $$

$$ (  \mathscr Z \cdot \mathscr W )^{\dagger_3}   =   \mathscr W^{\dagger_3}   \cdot \mathscr Z^{\dagger_3}    \, ;  \qquad   (  \mathscr Z \cdot \mathscr W )^\odot   =   \mathscr W^\odot    \cdot \mathscr Z^\odot   \, ;  $$

$$ (  \mathscr Z \cdot \mathscr W )^\diamond   =   \mathscr W^\diamond     \cdot \mathscr Z^\diamond    \, .  $$

One may ask  then   if  at  least   one of these conjugations is appropriate  in order to generate  an  $\H (\C)$--valued  inner  product  with  reasonably  good  properties.   The answer  is  positive  and  is  given  by

\begin{subtheorem}
The mapping
$$  \langle \cdot , \cdot \rangle_{ \H (\C) }  :   \H (\C)  \times   \H (\C)  \longrightarrow    \H (\C)  $$

given by
$$  \langle \mathscr Z  , \, \mathscr W \rangle_{ \H (\C) }    :=   Z \cdot  W^\odot ,$$
possesses  the  following  properties:

\begin{enumerate}

\item[(I)]  $  \langle \mathscr Z  , \, \mathscr W  +  \mathscr S   \rangle_{ \H (\C) }    =    \langle \mathscr Z  , \, \mathscr W     \rangle_{ \H (\C) }   +    \langle \mathscr Z  , \,   \mathscr S   \rangle_{ \H (\C) } ;$

\medskip

\item[(II)]   $   \langle     \Lambda \cdot    \mathscr Z  , \, \mathscr W     \rangle_{ \H (\C) }   =   \Lambda \cdot    \langle \mathscr Z  , \, \mathscr W     \rangle_{ \H (\C) } $;

\medskip

\item[(III)]    $\langle \mathscr Z  , \, \mathscr W    \rangle_{ \H (\C) }   =    \langle \mathscr W  , \, \mathscr Z     \rangle_{ \H (\C) }^\odot$;

\medskip

\item[(IV)]     $  \langle \mathscr Z  , \, \mathscr Z    \rangle_{ \H (\C) } =  \eta  + r\,  \i \,   \i_2 + s \,  \i  \,   \i_3$  with  $\eta \in \D^+$, $r , \, s \in \R$;  \\

\noindent   $  \langle \mathscr Z  , \, \mathscr Z    \rangle_{ \H (\C) }  =0 $  if  and  only  if  $Z=0$.

\end{enumerate}
\end{subtheorem}

{\sf Proof:}  \\      The  proof  of   (I),  (II)   and  (III)    follows  by      direct  computation  and  taking  into  account  that
\begin{equation}\label{ZWodot}
\begin{array}{rcl}
\mathscr Z \cdot \mathscr W^\odot    &   =  &   \displaystyle    \left(   Z_1 + Z_2 \i_2  \right)   \cdot   \left(   W_1^\ast  + \overline{W}_2 \i_2  \right)
\\  &  &  \\  & = &
\displaystyle    \left(   Z_1 W_1^\ast  +Z_2 W_2^\ast    \right)   +   \left(     Z_2  \overline{W}_1   - Z_1  \overline{W}_2 \right) \i_2 .
\end{array}
\end{equation}

Since
$$   \mathscr Z \cdot \mathscr Z^\odot       =       \left(   Z_1 Z_1^\ast  +Z_2 Z_2^\ast    \right)   +   \left(     Z_2  \overline{Z}_1   - Z_1  \overline{Z}_2 \right) \i_2 ,$$

we may conclude that       $  \langle \mathscr Z  , \, \mathscr Z    \rangle_{ \H (\C) }=0   $      if and  only  if   $Z=0$. In addition  we  know  that   $ Z_1 Z_1^\ast  +Z_2 Z_2^\ast  \in \D^+$  and   direct  computations  show  that  $ \left(  Z_2  \overline{Z}_1   - Z_1  \overline{Z}_2 \right)  \i_2 $  is  of  the  form   $  r \,  \i  \,   \i_2  \, +  \,   s  \,   \i  \,   \i_3$  with   $r , \, s \in \R$.
\mbox{} \qed \mbox{}  \\

\begin{subCy}  For  any  $\mathscr Z , \, \mathscr W \,$  and  $\Lambda $  in  $\H ( \C)$ it is:
$$  \displaystyle  \langle \mathscr Z  , \,  \Lambda  \cdot   \mathscr W    \rangle_{ \H (\C) }   =    \langle \mathscr Z  , \, \mathscr W     \rangle_{ \H (\C) }  \cdot  \Lambda^\odot .  $$

\end{subCy}

Thus,  the  biquaternionic  module  $\H ( \C)$  becomes  endowed   with  an  $\H ( \C)$--valued  inner    product  with  all  the  usual  properties.  Of  course the inner  product  square  is   neither  in  $\D^+$  nor  a  real  number but    it    has the peculiarity of  involving  the  three hyperbolic  imaginary  units.  We  believe  that   in    this  way  a  theory  of  biquaternionic  inner  product  spaces  can  be  constructed  along the same lines we have used to introduce the  theory  of  bicomplex  inner  product  modules.

\begin{subPn}  Let
  $ \H (\C)$  be  seen    as  a  $\bc$--module.  Then  the  first  bicomplex  component
$$    {\sf Z}_1 :=    Z_1 W_1^\ast  +Z_2 W_2^\ast  $$

of  the  biquaternionic  inner  product is    a       $\bc$--valued   inner  product.
\end{subPn}
{\sf Proof:}  \\     This  is an immediate consequence of the fact that we are considering  $ \H (\C)$  as  $\bc^2$ and thus  we can apply  the definition  and  properties  given in  paragraph  \ref{inner_prod_in_bc_n}.
\mbox{} \qed \mbox{}  \\

Some authors consider  linear sets  formed by  $ \H (\C)$--valued  functions  endowing  them  with  complex--valued  inner  products,   but  in  this  case   some     essential  structures    of  such  functions  are  lost.  The  Proposition  above  allows  to catch  more structures  endowing  such  sets  with  bicomplex--valued   inner  products  and consequenctly hyperbolic--valued  norms.


\chapter{Linear functionals  and  linear  operators  on  $\bc$--modules}

\section{Bicomplex  linear  functionals}

\begin{subParr}\label{primera seccion funcionales}     {\rm  Let $X$   be  a     $  \mathbb B  \mathbb C$-module, and let  $f: X \to   \mathbb B  \mathbb C$
be  a   $   \mathbb B  \mathbb C$--linear  functional. Then  for  any  $x \in X$  one  has:}

$$  \begin{array}{rcl}
f(x)  & =   &    f_{1 , \mathbf i }  (x)  \cdot \mathbf e  +f_{2 , \mathbf i }  (x)  \cdot \mathbf e^\dagger
=       f_{1 ,  \,  \mathbf j }  (x)  \cdot \mathbf e  +f_{2 ,  \,   \mathbf j }  (x)  \cdot \mathbf e^\dagger
\\  &   &   \\  & = &
F_1 (x)  + F_2 (x)\cdot  \mathbf j
=
G_1 (x)  + G_2 (x)\cdot  \mathbf i     \in   \mathbb B  \mathbb C
\end{array} $$

\end{subParr}

where   $ f_{\ell , \mathbf i } (x) \in \mathbb C ( \mathbf i ) $,   $   f_{\ell , \,  \mathbf j } (x)   \in  \mathbb C (\mathbf j) $,   $F_\ell (x)   \in \mathbb C ( \mathbf i ) $    and   $G_\ell  (x)    \in \mathbb C ( \mathbf j ) $,  for   $ \ell  =  1 , \, 2$.    This means that  $f$  induces  maps
$$ f_{1 , \mathbf i }   ,  \,f_{2 , \mathbf i } :  X    \to   \mathbb C ( \mathbf i ) , \qquad    f_{1 , \, \mathbf j }   ,  \,f_{2 ,  \,   \mathbf j } :  X    \to   \mathbb C ( \mathbf j ) ,  $$

$$F_1 , \, F_2  :   X    \to   \mathbb C ( \mathbf i ),     \qquad       G_1 , \, G_2  :   X    \to   \mathbb C ( \mathbf j ) ,$$

such that

$$f_{1 , \mathbf i }  =   F_1  -  \mathbf i \, F_2   ,  \qquad   f_{2 , \mathbf i }  =   F_1  +  \mathbf i \, F_2  , $$

and 

$$  f_{1 , \,  \mathbf j }  =   G_1  -  \mathbf j \, G_2  , \qquad    f_{2 ,  \,    \mathbf j }  =   G_1  +  \mathbf j \, G_2  .  $$

\medskip

In  particular we have  that
$$ \displaystyle     f_{ \ell , \i }  =    \pi_{ \ell , \i  } \circ f = \left(  \Pi_{ 1 , \i }  +  ( -1)^\ell  \, \i \, \Pi_{2 , \i }   \right) \circ f    ;   $$

$$ \displaystyle     f_{ \ell , \,  \j }  =    \pi_{ \ell , \,  \j  } \circ f = \left(  \Pi_{ 1 ,  \,  \j  }  +  ( -1)^\ell  \, \j \, \Pi_{2 ,  \,  \j }   \right) \circ f    ;   $$

$$ \displaystyle    F_\ell  =   \Pi_{ \ell , \, \j } \circ f =   \frac{  \,  ( \, \j  \, )^{ \ell  - 1  } }{ 2}  \left(  \pi_{ 1 , \, \j }  + (-1 )^{ \, \ell - 1 }  \, \pi_{2 , \, \j }  \right)   \circ f ; $$

$$ \displaystyle    G_\ell    =   \Pi_{ \ell ,  \i } \circ f =   \frac{  \,  ( \, \i   \,   )^{ \ell  - 1  } }{ 2}  \left(  \pi_{ 1 ,  \i }  + (-1 )^{ \, \ell - 1 }  \, \pi_{2 ,  \i }  \right)  \circ f  . $$

\medskip

Let  us  show   that     $ f_{1 , \mathbf i }   ,  \,f_{2 , \mathbf i } :  X_{  \mathbb C  (  \mathbf i ) }    \to   \mathbb C ( \mathbf i ) $  are  $  \mathbb C ( \mathbf i )$--linear  functionals  and  that     $ f_{1 , \, \mathbf j }   ,  \,f_{2 , \,  \mathbf j } :  X_{  \mathbb C  (  \mathbf j ) }  \to   \mathbb C ( \mathbf j ) $   are  $  \mathbb C ( \mathbf j )$--linear  functionals. This is equivalent  to requiring linearity for   $F_1 , \, F_2, \, G_1$  and  $G_2$.   Given  $\lambda  =  \lambda_1  +  \lambda_2 \, \mathbf j  \in \mathbb B \mathbb C$,  with   $ \lambda_1 , \,  \lambda_2   \in \mathbb C (\mathbf i) $,  and  given   $x, \, y  \in X$,  one has:
$$  \begin{array}{rcl}
f( \lambda \, x + y )  &  = &    f_{1, \mathbf i}  (  \lambda \, x + y )  \,   \mathbf e    +    f_{2, \mathbf i}  (  \lambda \, x + y )  \,   \mathbf e^\dagger
\\   &  &  \\   & = &
\lambda \, f(x)  +  f(y)
\\  &  &  \\  & = &
\displaystyle  \left(   ( \lambda_1  -  \mathbf i  \, \lambda_2 ) \, f_{1 , \mathbf i }  (x)   +   f_{1 , \mathbf i }  (y)   \right)  \, \mathbf e  \; +
\\  &  &  \\  &       &
\displaystyle    \qquad  \qquad   \quad  + \;     \left(   ( \lambda_1  +  \mathbf i  \, \lambda_2 ) \, f_{2 , \mathbf i }  (x)   +   f_{2 , \mathbf i }  (y)   \right)  \, \mathbf e^\dagger \, ,
\end{array} $$
hence
\begin{equation}\label{equation A}
 f_{1, \mathbf i}  (  \lambda \, x + y )   =   ( \lambda_1  -  \mathbf i  \, \lambda_2 ) \, f_{1 , \mathbf i }  (x)   +   f_{1 , \mathbf i }  (y)
 \end{equation}

and
\begin{equation}\label{equation B}
 f_{2, \mathbf i}  (  \lambda \, x + y )   =   ( \lambda_1  -  \mathbf i  \, \lambda_2 ) \, f_{2 , \mathbf i }  (x)   +   f_{2 , \mathbf i }  (y) \, .
\end{equation}

In  particular,  setting  $\lambda_2 =0$  we  have:
\begin{equation}\label{linealidad Ci}
 f_{1, \mathbf i}  (  \lambda_1 \, x + y )   =    \lambda_1   \, f_{1 , \mathbf i }  (x)   +   f_{1 , \mathbf i }  (y)
 \end{equation}

and
\begin{equation}\label{2a linealidad Ci}
 f_{2, \mathbf i}  (  \lambda_1 \, x + y )   =    \lambda_1  \, f_{2 , \mathbf i }  (x)   +   f_{2 , \mathbf i }  (y) \, .
\end{equation}

Thus,  the mappings  $f_{1 ,\i}$,   $f_{2 ,\i}$  are   $\C (\i)$   linear   functionals  on  $X_{ \C (\i)}$.  The case of   $  f_{1, \,  \mathbf j} $  and$  f_{2, \,  \mathbf j} $  follows in the same way.

\medskip

   In  particular,  equations  (\ref{equation A})   and  (\ref{equation B})  show us  the   ``type  of  homogeneity"   of  the   $\mathbb C ( \mathbf i)$--linear  functionals   $ f_{1, \mathbf i}$  and     $ f_{2, \mathbf i}$    with  respect   to  general  bicomplex   scalars.  Indeed,  for  $y=0$  equation  (\ref{equation A})  gives:
$$  f_{1, \mathbf i}  (  (  \lambda_1 +\mathbf j   \,       \lambda_2 ) x )  =   ( \lambda_1   -  \mathbf i \, \lambda_2 )  \, f_{1, \mathbf i} (x) \, ,$$

that is,
\begin{equation}\label{pi linealidad}
 f_{1, \mathbf i}  (    \lambda  \,  x )  =  \pi_{1, \mathbf i }  ( \lambda )  \cdot    f_{1, \mathbf i} (x) \,  .
\end{equation}

In  the  same  way  equation  (\ref{equation B})  gives:
\begin{equation}\label{2a pi linealidad}
f_{2, \mathbf i}  (    \lambda  \,  x )  =  \pi_{1, \mathbf i }  ( \lambda )  \cdot    f_{2, \mathbf i} (x) \,  .
\end{equation}

\medskip

To understand this phenomenon we note that if  $\lambda  =  \lambda_1 \e  + \lambda_2 \edag \in \bc$, and $x \in X$, then
$$ \begin{array}{rcl}
f ( \lambda x ) & = &  f_{1 , \i } ( \lambda x ) \e + f_{2 , \i } (\lambda x ) \edag
\\  &  &  \\  & = &
\pi_{1 , \i }  \circ f (\lambda \,x)  \, \e  +  \pi_{2 , \i }  \circ   f ( \lambda \,  x)   \, \edag
\\  &  &  \\  & = &
\pi_{1 , \i }  ( \lambda \, f (x) ) \, \e  +  \pi_{2 , \i }  ( \lambda \, f (x) ) \, \edag
\\  &  &  \\  & = &
\pi_{1 , \i }  ( \lambda ) \,   \pi_{1 , \i } ( f (x) ) \, \e  +  \pi_{2 , \i }  ( \lambda )  \,   \pi_{2 , \i } ( f (x) ) \, \edag
\\  &  &  \\  & = &
\pi_{1 , \i }  ( \lambda ) \,  f_{1 , \i }  (x)  \, \e  +  \pi_{2 , \i }  ( \lambda )  \,   f_{2 , \i }  (x)  \, \edag  .
\end{array} $$

Thus, if   $ \lambda = \lambda_1 \in \C ( \i )$, we obtain the  $ \C ( \i )$--linearity  of  $ f_{1 , \i }$  and $ f_{2 , \i } $.    In  general  when  $ \lambda \in \bc$  we  have   (\ref{pi linealidad})  and  (\ref{2a pi linealidad}).


\begin{subParr}\label{5.1.2}   As  in  the  classic  complex  case  the components  of  a   bicomplex  linear    functional   $f$  are  not  independent.  For example, one has that 
$$  \begin{array}{rcl}
f (\,    \mathbf j \, x )  & = &   F_1 (\,    \mathbf j \, x )   +   \mathbf j \, F_2 (\,    \mathbf j \, x )
\\  &  &  \\   & =  &
\mathbf j \, f(x)  =  -  F_2 ( x )   +   \mathbf j \, F_1 (  x ) \,
\end{array}  $$
\end{subParr}
which implies
$$   F_1 ( x )   =   F_2 (\,    \mathbf j \, x ) \, ,   $$
$$   F_2 ( x )   =   -     F_1 (\,    \mathbf j \, x ) \,  , $$

which  are of  course  mutually   reciprocal    relations. In  the  same  way
$$   G_1 ( x )   =   G_2 (  \mathbf i \, x ) \,   ,$$
$$   G_2 ( x )   =   -     G_1 (  \mathbf i \, x ) \,    .$$

This phenomenon is even more remarkable for the functionals $  f_{1,  \mathbf i} $,   $  f_{2,   \mathbf i} $,    $  f_{1, \,  \mathbf j} $  and    $  f_{2, \,  \mathbf j} $.  Indeed
$$  \begin{array}{rcl}
f_{1 , \mathbf i }    (\,    \mathbf j \, x )  & = &   F_1 (\,    \mathbf j \, x )   -   \mathbf i \, F_2 (\,    \mathbf j \, x )   =   - F_2(x)  - \mathbf i \, F_1 (x)
\\  &  &  \\   & =  &
-  \mathbf i \,    \left(    F_1 (x)    - \mathbf i \,   F_2 ( x )  \right)   =   - \mathbf i \,   f_{1 , \mathbf i }   (  x )
\end{array}  $$

and
$$  \begin{array}{rcl}
f_{2 , \mathbf i }    (\,    \mathbf j \, x )  & = &   F_1 (\,    \mathbf j \, x )  +   \mathbf i \, F_2 (\,    \mathbf j \, x )   =   - F_2(x)  + \mathbf i \, F_1 (x)
\\  &  &  \\   & =  &
 \mathbf i \,    \left(    F_1 (x)    + \mathbf i \,   F_2 ( x )  \right)   =    \mathbf i \,   f_{2 , \mathbf i }   (  x )  \, ,
\end{array}  $$

that is,
$$
f_{1 , \mathbf i }    (\,    \mathbf j \, x )   =  - \mathbf i \,   f_{1 , \mathbf i }   (  x )  ,  $$

$$
f_{2 , \mathbf i }    (\,    \mathbf j \, x )    =    \mathbf i \,   f_{2 , \mathbf i }   (  x )  \, . $$

In  the  same  way
$$
f_{1 , \, \mathbf j }    (   \mathbf i \, x )   =  - \mathbf j \,   f_{1 , \,  \mathbf j }   (  x )    \, ,$$

$$
f_{2 , \,  \mathbf j }    (  \mathbf i \, x )    =    \mathbf j \,   f_{2 , \,  \mathbf j }   (  x )  \, . $$


Consider now a bicomplex  linear  functional  $ f : X \to  \mathbb B \mathbb C$. One has that  $ f(x) = f ( \mathbf e \, x + \mathbf e^\dagger \, x ) =  f (x  \, \mathbf e ) \, \mathbf e + f (x  \, \mathbf e^\dagger) \, \mathbf e^\dagger   =  f_{1 , \mathbf i }  (x  \, \mathbf e ) \, \mathbf e  +     f_{2 , \mathbf i }   (x  \, \mathbf e^\dagger) \, \mathbf e^\dagger   =    f_{1 , \,  \mathbf j }  (x  \, \mathbf e ) \, \mathbf e  +     f_{2 ,  \,   \mathbf j }   (x  \, \mathbf e^\dagger) \, \mathbf e^\dagger $,  \,  i. e.,
\begin{equation}\label{decomposition functional}
\begin{array}{rcl}
f(x) & =   &   f_{1 , \mathbf i }  (x  \, \mathbf e ) \, \mathbf e  +     f_{2 , \mathbf i }   (x  \, \mathbf e^\dagger) \, \mathbf e^\dagger
\\   &  &  \\  &   =  &      f_{1 , \,  \mathbf j }  (x  \, \mathbf e ) \, \mathbf e  +     f_{2 ,  \,   \mathbf j }   (x  \, \mathbf e^\dagger) \, \mathbf e^\dagger \, .
\end{array}
\end{equation}

  Let  us  assume  now  that  $X$  has  a  real--valued  or  a  hyperbolic  norm.        Then     equalities   (\ref{decomposition functional}) imply   that  $f$  is  continuous if and  only  if   $f_{\ell , \mathbf i }$  ,  $f_{\ell, \, \mathbf j }$,  $\ell = 1,2$,  are  continuous  on  $X_{\mathbb C (\mathbf i ) }$  and  on   $X_{\mathbb C (\, \mathbf j ) }$ respectively.  In  particular  if $f$  is continuous,  then the restrictions     $f_{1 , \mathbf i }   \mid_{ X_{ \mathbf e , \mathbf i} } $  and     $f_{2,  \mathbf i } \mid_{X_{  \mathbf e^\dagger, \mathbf i } }  $  are   $\mathbb C (\mathbf i ) $--linear  continuous  functionals;  in the same manner   $f_{1 , \,  \mathbf j }   \mid_{ X_{ \mathbf e , \,  \mathbf j} } $  and     $f_{2,   \,  \mathbf j } \mid_{X_{  \mathbf e^\dagger,  \,   \mathbf j } }  $  are  $\mathbb C (\mathbf j ) $--linear  continuous  functionals.   \\

We are now ready to prove a bicomplex version of the celebrated Riesz theorem (technically, such theorem was already proved in \cite{LMR2010}, but 
we decided to include its proof here for the sake of completeness, and because we can frame the proof in the context of our more general treatment of bicomplex functionals).

\begin{subtheorem}  (Riesz).  Let $X$  be  a  bicomplex Hilbert     space  and  let  $f: X \to    \mathbb B  \mathbb C$  be  a  continuous  $ \mathbb B  \mathbb C$-linear  functional.  Then  there  exists  a  unique  $y \in X$  such  that for every   $x \in X$,  $f(x) =  \langle  x  ,  y   \rangle$.
\end{subtheorem}

{\sf Proof.}  \\   The  first  representation  in   (\ref{decomposition functional})  generates two   $\mathbb C (\mathbf i ) $--linear    continuous  functionals  on  $ X_{  \mathbf e, \mathbf i }  $  and
$X_{  \mathbf e^\dagger, \mathbf i } $ for  each of  which  the  complex  Riesz  representation   Theorem  gives a  unique element    $ u   =   \e \,    u   \in     X_{  \mathbf e, \mathbf i }  $  and  $ v   =  \edag   \,   v  \in     X_{  \mathbf e^\dagger, \mathbf i }  $  such  that
\begin{equation}\label{u en Xei}
f_{1 , \mathbf i }   \mid_{ X_{ \mathbf e , \mathbf i} }  ( \mathbf e \, x)  =  \langle   \e  x , \mathbf e \, u \rangle_{ X_{ \mathbf e , \mathbf i}  }
\end{equation}

and
\begin{equation}\label{v en Xe dagger i}
f_{2 , \mathbf i }   \mid_{ X_{ \mathbf e^\dagger , \mathbf i} }  ( \mathbf e^\dagger \, x)  =  \langle   \edag    x , \mathbf e^\dagger \,  v \rangle_{ X_{ \mathbf e^\dagger , \mathbf i}  }    \, .
\end{equation}

The  obvious   candidate    $y:=  u + v     =  \e u + \edag v  \in X$  allows us to conclude the proof since
$$ \begin{array}{rcl}
 \langle      y   ,  x  \rangle   &  =  &      \langle    \mathbf e \,  x  +  \mathbf e^\dagger x  ,      \mathbf e \, u + \mathbf e^\dagger \, v  \rangle
\\  &  &  \\  & = &     \e   \,   \langle      \mathbf e \,  x ,     \mathbf e \, u   \rangle_{X_{\e  , \i } }      +   \edag \,     \langle     \mathbf e^\dagger \,  x   \, , \,    \mathbf e^\dagger \, v  \rangle_{X_{\edag  , \i } }
\\  &  &  \\  & = &    \mathbf e \,  f_{1, \mathbf i }  (  \mathbf e \, x  )       +       \mathbf e^\dagger \,  f_{2 , \mathbf i }  (    \mathbf e^\dagger \, x  )  =   f(x) \, .
\end{array}
$$

 Since  the  bicomplex  inner  product  is  non  degenerate  then  the  element  $y$  is unique.  At the same  time  this  unicity  implies  that  if we repeat the above  process  using  the  $ \mathbb C ( \mathbf j )$--linear  spaces    $ X_{ \mathbf e ,  \,  \mathbf j}   $  and       $    X_{ \mathbf e^\dagger ,  \,  \mathbf j}  $  the  element  that  we would get will  be  the  same.      \mbox{}  \qed  \mbox{}  \\

\medskip

The next result can also be found in   \cite{LMR2010}.

\begin{subtheorem} (First  bicomplex Schwarz inequality).
Let  $X$   be    a  bicomplex   Hilbert  module, and    let   $x,y \in X$.  Then
$$  \displaystyle     |   \langle      x  ,  y  \rangle |  \leq   \,   \sqrt{ 2 }     \,    \| x \|  \,  \| y \| \, .$$
\label{Schwarz}
\end{subtheorem}
{\sf Proof:} \\     Recall   that   $X$  is  the direct  sum    of    $X_{ \mathbf e }$  and    $X_{ \mathbf e^\dagger }$  when  they  are  seen   as  $ \mathbb C ( \mathbf i)$--  or   $ \mathbb C ( \mathbf j)$--complex  Hilbert  spaces. For any  $x,y \in X$ one has
$$ \begin{array}{rcl}
|   \langle      x  ,  y  \rangle |  & = &   \displaystyle  \left|    \,       \langle      \mathbf e \,  x   +  \mathbf e^\dagger \, x ,     \mathbf e \, y   +    \mathbf e^\dagger \, y   \rangle     \right|
\\  &  &  \\  & = &
\displaystyle  \left|    \,  \mathbf e \,     \langle      \mathbf e \,  x ,     \mathbf e \, y   \rangle_{X_{ \mathbf e  , \mathbf i } }    +    \mathbf e^\dagger
\,        \langle     \mathbf e^\dagger \,  x , \,    \mathbf e^\dagger \, y  \rangle_{X_{ \mathbf e^\dagger  , \mathbf i } }   \right|
\end{array} $$

\qquad   \qquad  \qquad  \qquad     \qquad  \qquad     using  (\ref{estrella})
$$ \qquad  \qquad   \qquad  \quad   =  \;      \displaystyle   \frac{1}{ \sqrt{2} } \,   \left(  \,   \left|    \langle     \mathbf e \,       x ,     \mathbf e \, y     \rangle_{X_{ \mathbf e  , \mathbf i } }  \right|^2  +   \left|    \langle   \mathbf e^\dagger
\,  x , \,    \mathbf e^\dagger \, y    \rangle_{X_{ \mathbf e^\dagger  , \mathbf i } }    \right|^2 \right)^{1/2}   $$

\qquad    \qquad     \qquad  \qquad     using  the  complex  Schwarz  inequality
$$ \begin{array}{rcl}
\qquad  &  \leq  &
 \displaystyle   \frac{1}{ \sqrt{2} } \,   \left(  \,   \left\|        \mathbf e \,       x \right\|_{X_{ \mathbf e  , \mathbf i } }^2  \cdot  \left\|      \mathbf e \, y   \right\|_{X_{ \mathbf e  , \mathbf i } }^2     +     \left\|     \mathbf e^\dagger
\,  x \right\|^2_{X_{ \mathbf e^\dagger  , \mathbf i } } \cdot  \left\|      \mathbf e^\dagger \, y   \right\|^2_{X_{ \mathbf e^\dagger  , \mathbf i } }       \right)^{1/2}
\\  &  &  \\  & \leq  &
 \displaystyle   \frac{1}{ \sqrt{2} } \,   \left(  \,    2 \,     \left\|       x \right\|^2   \left(    \left\|      \mathbf e \, y   \right\|_{X_{ \mathbf e  , \mathbf i } }^2     +     \left\|      \mathbf e^\dagger \, y   \right\|^2_{X_{ \mathbf e^\dagger  , \mathbf i } }   \, \right)         \right)^{1/2}
\\  &  &  \\  & =     &
 \displaystyle    \,   \left(  \,   2  \,    \left\|       x \right\|^2  \cdot     \left\|       y   \right\|^2      \right)^{1/2}   =   \,   \sqrt { 2 }   \,   \,  \| x \|  \cdot  \|  y \| \, ,
\end{array}$$

which concludes the proof.        \mbox{} \qed \mbox{} \\

Since  we  have  a  hyperbolic  norm  which  takes  values  in  the  partially  ordered  set  $\D^+$,  we  are able to prove  a  ``hyperbolic  version"  of  the  Schwarz  inequality.

\begin{subtheorem} (Second  bicomplex Schwarz inequality).
Let  $X$   be    a  bicomplex   Hilbert  module,   and  let   $x,y \in X$.  Then
$$  \displaystyle     |   \langle      x  ,  y  \rangle |_\k  \lessdot   \,      \| x \|_\D  \cdot   \| y \|_\D \, .$$
\label{Schwarz2}
\end{subtheorem}
{\sf Proof:}  \\    As  before,
$$    \langle      x  ,  y  \rangle  =      \,  \mathbf e \,     \langle      \mathbf e  x ,     \mathbf e  y   \rangle_{X_\e }    +    \mathbf e^\dagger
\,        \langle     \mathbf e^\dagger \,  x , \,    \mathbf e^\dagger \, y  \rangle_{X_\edag}  \, ,$$

hence  
$$ \begin{array}{rcl}
\displaystyle     |   \langle      x  ,  y  \rangle |_\k  & = &
  \mathbf e \,   |    \,    \langle      \mathbf e  x ,     \mathbf e  y   \rangle_{X_\e } |      +    \edag \cdot  |
        \langle     \mathbf e^\dagger   x , \,    \mathbf e^\dagger  y  \rangle_{X_\edag }   |
\end{array} $$

\qquad   \qquad  \qquad  \qquad     \qquad      using  the  complex  Schwarz  inequality
$$\begin{array}{l}
     \qquad  \qquad     \lessdot    \;      \displaystyle    \e  \,  \|  \e x \|_{X_\e}  \cdot    \|  \e y \|_{X_\e}  +    \edag  \,  \|  \edag x \|_{X_\edag}  \cdot    \|  \edag y \|_{X_\edag}     \\  \\
 \qquad     \qquad     = \left(  \e \cdot  \|  \e x \|_{X_\e}  +   \edag  \,  \|  \edag x \|_{X_\edag}  \right)  \cdot   \left(  \e \cdot  \|  \e y \|_{X_\e}  +   \edag  \,  \|  \edag y \|_{X_\edag}  \right)  \\  \\       \qquad  \qquad  = \;     \| x \|_\D \cdot \| y \|_\D .
 \end{array}$$
 \mbox{} \qed  \mbox{}  \\

Notice the  absence of  the coefficient  $\sqrt{ 2 }$   in  the  second  inequality;  this is  due to  the compatibility  between  the  hyperbolic  modulus  of  the  inner  product  (which  is  a  bicomplex  number)  and  the  hyperbolic  norm  on  $X$.


\begin{subEx}

We  are interested  again in the situation in which  $  X  =   \mathbb B \mathbb C.$ Take then  $x=Z =  \beta_1 \, \mathbf e + \beta_2  \,  \mathbf e^\dagger$, and  $ y = W =    \xi_1 \, \mathbf e + \xi_2 \, \mathbf e^\dagger$,    $\beta_1 $,   $\beta_2$, $\xi_1$,  $\xi_2    \in \mathbb C ( \mathbf i ) $  so that
$$  \displaystyle    \| Z \|_X^2  =   \frac{1}{2}   \,  \left(    |  \beta_1 |^2  + | \beta_2 |^2     \right)   =   \, | Z |^2    \,   ,   \quad    \| W \|_X^2  =   \frac{1}{2}  \, \left(  |  \xi_1 |^2  + | \xi_2 |^2   \right)    =   \, | W |^2   .  $$
\end{subEx}

One has directly:
$$ \begin{array}{rcl}
| \langle Z ,W \rangle |  & = &  | Z \cdot W^\ast |  =  |  (  \beta_1 \, \mathbf e +   \beta_2 \, \mathbf e^\dagger ) (  \overline{ \xi}_1 \, \mathbf e + \overline{ \xi}_2 \, \mathbf e^\dagger ) |
\\  &  &  \\  & = &
| \beta_1  \,  \overline{ \xi}_1 \, \mathbf e +  \beta_2 \, \overline{\xi}_2 \, \mathbf e^\dagger |
\\  &  &  \\  & = &
\displaystyle  \frac{1}{ \sqrt{ 2 } }  \, \sqrt{     | \beta _1 |^2 \cdot   |  \overline{ \xi}_1 |^2    \,   +    \, | \beta _2 |^2 \cdot   |  \overline{ \xi}_2 |^2  \,   }
\\  &  &  \\  & \leq &
\displaystyle   \frac{1}{ \sqrt{2} }   \sqrt{   \, 2 \,    \|Z \|^2_{ \mathbb C }  \left(   | \xi_1 |^2  + | \xi_2 |^2  \right)  }
\\   &  &  \\   &  =  &
\displaystyle   \sqrt{   \, 2 \,       \|Z \|^2_X   \cdot     \| W \|^2_X  }
 = \;         \sqrt{2}   \,   \|   Z  \|_X     \cdot   \|  W \|_X  \, ,
\end{array} $$

in accordance with  Theorem   \ref{Schwarz}.

\medskip

Similarly,  for the hyperbolic  modulus  of  the  inner  product  and  for    the  corresponding  hyperbolic    norm  one  has:
$$ \begin{array}{rcl}
|  \langle  Z , W \rangle |_\k & = &   |  Z \cdot W^\ast |_\k  =  | Z |_\k  \cdot  |  W^\ast |_\k
\\   &  &   \\  & = &
 | Z |_\k  \cdot  |  W |_\k   \; = \;     \| Z \|_\D  \cdot  \|  W  \|_\D \, ;
 \end{array}$$

This means that,  in this case, the  Second  Schwarz  inequality  becomes  an equality, exactly as   it    happens with the Schwarz  inequality in $\C$.

\medskip


\section{Polarization identities}

Let  $X$  be  a  real  linear  space  with  a  bilinear  symmetric  form    $\mathcal B_\R ( \cdot , \cdot ) $, and  denote  by   $\mathcal Q_\R ( \cdot )$  the  corresponding  quadratic  form  generated  by $\mathcal B_\R$, i.e., for  any  $x   \in X  $,
$$  \mathcal Q_\R (x)  :  = \mathcal B_\R (x,x) . $$

Then the following  relation  holds  for  any  $x , y \in X$:
\begin{equation}\label{polar real}
\displaystyle  \mathcal B_\R ( x,y)   = \frac{1}{4}  \left(   \mathcal Q_\R ( x+y)  - \mathcal Q_\R (x-y) \right) .
\end{equation}

This relationship usually goes under the name of  Polarization  Identity  for  real  linear  spaces.  In  particular,  a  real  inner product  and 
its associated norm satisfy  (\ref{polar real}).

\medskip

For  a  complex  linear  space  $X$    the  form  $\mathcal B_\C$  is  assumed  to  be  Her\-mi\-tian  and  sesquilinear,  that  is,    $\mathcal B_\C (x,y)   =  \overline{ \mathcal B_\C (y,x) } $  and   $\mathcal B_\C ( \lambda x , \mu y )  =  \lambda \overline{ \mu } \mathcal B_\C (x,y)$  \, for any  $x,y \in X $, and for all $\lambda , \, \mu \in \mathbb C$.  let $\mathcal Q_\C (x) :  =   \mathcal B_\C (x,x)$ denote  the  quadratic  form.  The  polarization  identity  in  this  case  becomes
\begin{equation}\label{complex polar}
\begin{array}{l}
\displaystyle  \mathcal B_\C ( x,y)   = \frac{1}{4}  \left(   \mathcal Q_\C ( x+y)  - \mathcal Q_\C (x-y)    +  \right.
\\  \\
\displaystyle  \qquad   \qquad  \qquad  \quad    \left.   + \;     \i   \, \left(    \mathcal Q_\C ( x+ \i \, y)  - \mathcal Q_\C (x-   \i \, y)   \right)    \right) .
\end{array}
\end{equation}

\medskip

Again  any  complex   inner  product  and  its  corresponding  norm  satisfy  (\ref{complex polar}).

Let   now  $X$  be  a  bicomplex  module.   We   will   consider  maps
$$ \mathcal B  :  X \times  X  \to  \bc $$

such  that   $ \mathcal B  ( x , y_1 + y_2)  = \mathcal B (x,y_1)  + \mathcal B (x ,y_2)$  and    $ \mathcal B  ( x_1 + x_2 , y ) = \mathcal B ( x_1 , y) + \mathcal B (x_2 , y )$  \, for all $ x_1, \, x_2 , \,  x , \, y_1 , \, y_2  , \, y \in X$.   For    such a map  $\mathcal B$,  to be  ``Hermitian"  may mean one of the  three  properties:

\medskip

\begin{enumerate}

\item[({\sc i})]   If  $ \mathcal B (x,y)  =  \overline{ \mathcal B (y,x) }$, (in this case $\mathcal B$ is called  bar--Hermitian),

\medskip

\item[({\sc ii})]   If  $ \mathcal B (x,y)  =   ( \mathcal B (y,x) )^\dagger$,  (in this case $\mathcal B$ is called  $\dagger$--Hermitian),

\medskip

\item[({\sc iii})]   If  $ \mathcal B (x,y)  =   ( \mathcal B (y,x) )^\ast  $, (in this case $\mathcal B$ is called  $\ast$--Hermitian).

\end{enumerate}

\medskip

Similarly  we have the  three types of sesquilinearity. If    $ \nu , \, \mu \in \bc$  then

\medskip

\begin{enumerate}

\item[({\sc i'})]      $\mathcal B$ is   bar--sesquilinear  if
$$ \mathcal B ( \nu  x ,  \mu  y )  =   \nu \, \overline{ \mu } \, \mathcal B (x,y) . $$

\medskip

\item[({\sc ii'})]   $\mathcal B$ is called  $\dagger$--sesquilinear  if
$$ \mathcal B ( \nu  x ,  \mu  y )  =   \nu \,  \mu^\dagger  \, \mathcal B (x,y) . $$

\medskip

\item[({\sc iii'})]  $\mathcal B$ is  $\ast$--sesquilinear  if
$$ \mathcal B ( \nu  x ,  \mu  y )  =   \nu \,  \mu^\ast  \, \mathcal B (x,y) . $$

\end{enumerate}

\medskip

For reasons that will be readily apparent,  we  are  going  to  consider  $\ast$--Hermitian  and  $\ast$--ses\-qui\-li\-near  maps  $\mathcal B$.  Note  that  since  $\mathcal Q (x) : = \mathcal B (x,x) = \mathcal B (x,x)^\ast$  then  $ \mathcal Q (x)$  takes  values  in  $\mathbb D$.  Hence  the  inner  products that we have defined fit  into  this  notion:  if  $\mathcal B (x,y) = \langle x , y \rangle_X$  then  $\mathcal B$ is a $\ast$--Hermitian and $\ast$--sesquilinear  form, and  $\mathcal Q(x) = \mathcal B (x,x)$  is  the  square of    a  hyperbolic norm on $X$.

\medskip

It  turns  out  that   given  a  $\ast$--Hermitian  and  $\ast$--sesquilinear  bicomplex  form  $\mathcal B$, there are several polarization--type
 identities that connect it to its  quadratic  form  $\mathcal Q$. Indeed,  it  is  a  simple   and  direct  computation  to  show  that    for  any   $x, \, y \in X$  one has
\begin{equation}\label{1i}
\begin{array}{l}
\displaystyle  \mathcal B (x,y)   = \frac{1}{4} \left(  \mathcal Q (x+y)  -    \mathcal Q (x  - y)  \; + \right.
\\  \\
\displaystyle   \qquad \qquad  \quad   \left.   + \;   \i   \left(    \mathcal Q (x +  \i  y)   -     \mathcal Q (x    - \i  y)  \right)  \right)  ;
\end{array}
\end{equation}
\begin{equation}\label{jk}
\begin{array}{l}
\displaystyle  \mathcal B (x,y)   = \frac{1}{4} \left(  \j  \left(    \mathcal Q (x+  \j y)  -    \mathcal Q (x  -   \j  y)   \right)   \;   +  \right.
\\  \\
\displaystyle  \qquad \qquad \quad   \left.  + \;  \k   \left(    \mathcal Q (x +  \k  y)   -     \mathcal Q (x    - \k  y)  \right)  \right)  ;
\end{array}
\end{equation}
\begin{equation}\label{1j}
\begin{array}{l}
\displaystyle  \mathcal B (x,y)   = \frac{1}{4} \left(    \mathcal Q (x+   y)  -    \mathcal Q (x  -   y)    \;    +   \right.
\\  \\
\displaystyle \qquad \qquad \quad  \left.   + \;   \j   \left(    \mathcal Q (x +  \j  y)   -     \mathcal Q (x    - \j  y)  \right)  \right)  ;
\end{array}
\end{equation}
\begin{equation}\label{ik}
\begin{array}{l}
\displaystyle  \mathcal B (x,y)   = \frac{1}{4} \left(  \i  \left(    \mathcal Q (x+  \i y)  -    \mathcal Q (x  -   \i  y)   \right)    \;  +  \right.
\\  \\
\displaystyle    \qquad  \qquad  \quad  \left.   + \;     \k   \left(    \mathcal Q (x +  \k  y)   -     \mathcal Q (x    - \k  y)  \right)  \right)  .
\end{array}
\end{equation}
The validity of formulas   (\ref{1i}) and (\ref{1j}) can be concluded  from   (\ref{complex polar})  but   there are deep   bicomplex   reasons
  for the validity of (\ref{1i})--(\ref{ik}).  Indeed,  $\mathcal B (x,y)$  can be written as
\begin{equation}\label{idempotent repr of a form}
 \begin{array}{rcl}
\mathcal B (x,y) & = & \mathcal B ( x \e  + x \edag , y \e  + y \edag )
\\  &  &  \\  & = &
\e  \,  \mathcal B ( x \e , y \e  )  + \edag  \mathcal B  ( x \edag  ,  y  \edag)
\\  &  &  \\   & = &
\e   \,    \pi_{1 , \i }  (  \mathcal B  ( x \e , y \e ) )  +   \edag  \pi_{2 , \i }   (   \mathcal B  ( x \edag  ,  y  \edag) )
\\  &  &  \\   & = &
\e   \,    \pi_{1 , \,  \j }  (  \mathcal B  ( x \e , y \e ) )  +   \edag  \pi_{2 ,  \,  \j }   (   \mathcal B  ( x \edag  ,  y  \edag) )  ,
\end{array}
\end{equation}
so that
$$ \begin{array}{rcl}
\mathcal Q (x)   & = &  \mathcal B (x,x)
\\  &  &  \\  & = &
\e \,  \pi_{1 , \i }  (    \mathcal B ( x \e , y \e  )  )  + \edag   \pi_{ 2 , \i }  (   \mathcal B  ( x \edag  ,  y  \edag)   )
\\  &  &  \\   & = &
\e  \,   \pi_{1 , \i }  (    \mathcal Q ( x \e  )  )  + \edag   \pi_{ 2 , \i }  (   \mathcal Q  ( x \edag  )   )
\\  &  &  \\   & = &
\e \,     \pi_{1 , \,  \j }  (    \mathcal B ( x \e , y \e  )  )  + \edag   \pi_{ 2 ,  \,  \j }  (   \mathcal B  ( x \edag  ,  y  \edag)   )
\\  &  &  \\   & = &
\e     \,     \pi_{1 ,  \,   \j }  (    \mathcal Q ( x \e  )  )  + \edag   \pi_{ 2 ,  \,  \j }  (   \mathcal Q  ( x \edag  )   )  .
\end{array}  $$

By applying this to the  right--hand side  of  (\ref{jk}), we see that for the first  summand one has
$$  \begin{array}{l}
\displaystyle  \j  \,   \left(   \mathcal Q  ( x + \j y )  -     \mathcal Q  ( x - \j y )  \right)      =
\\    \\
\qquad  =  \;      \displaystyle   \j \, \left(   \e \, \pi_{1 , \i }   \left(       \mathcal Q ( ( x + \j y )   \e )  \right)    +   \edag  \pi_{2 , \i }  \left(      \mathcal Q   (   ( x + \j y )   \edag )  \right)      \right.
\\    \\
\qquad   \quad     \left.  \displaystyle   - \;        \e \, \pi_{1 , \i }   \left(       \mathcal Q ( ( x - \j y )   \e )  \right)    -   \edag  \pi_{2 , \i }  \left(      \mathcal Q   (   ( x - \j y )   \edag )  \right)    \right)
\\    \\
\qquad  =  \;      \displaystyle    -  \i    \e \, \pi_{1 , \i }   \left(       \mathcal Q ( ( x -    \i y )   \e )  \right)    + \i   \edag  \pi_{2 , \i }  \left(      \mathcal Q   (   ( x + \i y )   \edag )  \right)
\\    \\
\qquad   \quad      \displaystyle   + \;      \i  \e \, \pi_{1 , \i }   \left(       \mathcal Q ( ( x + \i y )   \e )  \right)    -  \i     \edag  \pi_{2 , \i }  \left(      \mathcal Q   (   ( x - \i y )   \edag )  \right)
\\    \\
\qquad  =  \;      \displaystyle      \i   \,  \left(       \e \, \pi_{1 , \i }   \left(       \mathcal Q ( ( x +   \i y )   \e )  \right)    +    \edag  \pi_{2 , \i }  \left(      \mathcal Q   (   ( x + \i y )   \edag )  \right)   \right.
\\    \\
\qquad   \quad   \left.      \displaystyle   - \;      \left(     \e \, \pi_{1 , \i }   \left(       \mathcal Q ( ( x - \i y )   \e )  \right)    +      \edag  \pi_{2 , \i }  \left(      \mathcal Q   (   ( x - \i y )   \edag )  \right)    \right)    \right)
\\   \\
\qquad  =   \displaystyle  \;   \i \, \left(    \mathcal Q  (  x +  \i y )  -   \mathcal Q  ( x - \i y )  \right)   ;
\end{array}  $$
as to the  second  summand  we  have:
$$  \begin{array}{l}
\displaystyle  \k  \,   \left(   \mathcal Q  ( x + \k y )  -     \mathcal Q  ( x - \k y )  \right)      =
\\    \\
\qquad  =  \;      \displaystyle   \k \, \left(   \e \, \pi_{1 ,  \i }   \left(       \mathcal Q ( ( x + \k y )   \e )  \right)    +   \edag  \pi_{2 , \i }  \left(      \mathcal Q   (   ( x + \k y )   \edag )  \right)      \right.
\\    \\
\qquad   \quad     \left.  \displaystyle   - \;        \e \, \pi_{1 , \i }   \left(       \mathcal Q ( ( x - \k y )   \e )  \right)    -   \edag  \pi_{2 , \i }  \left(      \mathcal Q   (   ( x - \k y )   \edag )  \right)    \right)
\\    \\
\qquad  =  \;      \displaystyle        \e \, \pi_{1 , \i }   \left(       \mathcal Q ( ( x + y )   \e )  \right)    -   \edag  \pi_{2 , \i }  \left(      \mathcal Q   (   ( x -  y )   \edag )  \right)
\\    \\
\qquad   \quad      \displaystyle   - \;       \e \, \pi_{1 , \i }   \left(       \mathcal Q ( ( x -   y )   \e )  \right)    +     \edag  \pi_{2 , \i }  \left(      \mathcal Q   (   ( x + y )   \edag )  \right)
\\    \\
\qquad  =   \displaystyle  \;    \mathcal Q  (  x +   y )  -   \mathcal Q  ( x -  y )  .
\end{array}  $$

One can now compare the  right--hand  sides  in    (\ref{1i})-- (\ref{ik}) to conclude that  they  all  coincide.


Take a bicomplex module $X$. As we saw in (\ref{idempotent repr of a form}), any Hermitian  $\ast$--ses\-qui\-li\-near    form  $\mathcal B$  has ``an  idempotent  representation"
$$\begin{array}{rcl}
\mathcal B (x,y)    &  =   &    \e  \, \mathcal B ( \e x , \e y )  + \edag \mathcal B  ( \edag x , \edag y )
\\  &  &   \\   & = &
 \e  \,  \pi_{1 , \i }  (  \mathcal B ( \e x , \e y )  )   + \edag   \, \pi_{2 , \i }  (   \mathcal B  ( \edag x , \edag y ) )   .
 \end{array}$$

On   $X_\e  = \e X  $  set
$$  \mathfrak b_{1 , \i }      ( \e x , \e y )  :=   \pi_{1 , \i }  (  \mathcal B  ( \e x , \e y ) ) .$$

We will show  that   $  \mathfrak b_{1 , \i }   $   is  a  Hermitian  $ \C ( \i) $--sesquilinear  form  on  $ X_{ \e , \i }$. Indeed,

\medskip

\begin{enumerate}

\item[1)]    $ \overline{   \mathfrak b_{1 , \i }      ( \e y , \e x ) }  \;  = \;    \overline{  \,   \pi_{1 , \i }  (\mathcal B      ( \e y , \e x ) ) } $
\\   \\
\hbox{}   \quad   \qquad \qquad  $ = \;    \overline{ \,  \pi_{1 , \i }  (\mathcal B      ( \e x , \e y )^\ast ) }  $
\\    \\
\hbox{}     \quad    \qquad  \qquad   $  \displaystyle    =  \;    \overline{     \pi_{1 , \i }   \left(  \left(  \e \,   \pi_{1 , \i }  (\mathcal B     ( \e x , \e y ) )   +   \edag \,      \pi_{2 , \i }  (\mathcal B     ( \edag  x , \edag  y )   )  \right)^\ast   \right)   \,     }   $
\\    \\
\hbox{}     \quad    \qquad  \qquad   $  \displaystyle    =  \;    \overline{     \pi_{1 , \i }   \left(   \e \, \overline{  \pi_{1 , \i }  (\mathcal B     ( \e x , \e y ) ) }   +   \edag \,  \overline{    \pi_{2 , \i }  (\mathcal B     ( \edag  x , \edag  y )   ) }      \right)   \,     }   $
\\    \\
\hbox{}     \quad    \qquad  \qquad   $  \displaystyle    =  \;    \overline{   \overline{       \pi_{1 , \i }   (   \mathcal B     ( \e x , \e y ) ) }  }  \; = \;     \pi_{1 , \i }   (   \mathcal B     ( \e x , \e y ) ) $   \\   \\
\hbox{}   \quad  \qquad  \qquad   $  =      \mathfrak b_{1 , \i }      ( \e x , \e y )     $.

\medskip

\item[2)]      $   \mathfrak b_{1 , \i }      (  \mu \,  \e x   + \e z   , \e y )   \;  = \;             \pi_{1 , \i }  (    \mathcal B      (   \mu \,  \e x   + \e z  ,    \e y  ) )  $
\\   \\
\hbox{}   \qquad      \qquad   \qquad \qquad  $ = \;      \pi_{1 , \i }  (    \mu \,   \mathcal B      ( \e x , \e y )   +    \mathcal B      ( \e z , \e y )     )   $
\\    \\
\hbox{}     \qquad    \qquad    \qquad      \qquad   $  \displaystyle    =  \;      \mu \,    \mathfrak b_{1 , \i }      ( \e x , \e y )   +    \mathfrak b_{1 , \i }      ( \e z , \e y )         $.

\end{enumerate}

\medskip

The two  properties  together  imply  that
$$        \mathfrak b_{1 , \i }      ( \e x ,      \mu \,    \e y )   =   \overline{ \mu } \,       \mathfrak b_{1 , \i }      ( \e x , \e y )  .$$

In the same  way we can define three more forms
$$  \mathfrak b_{2 , \i }      ( \edag  x , \edag  y )  :=   \pi_{2 , \i }  (  \mathcal B  ( \edag x , \edag y ) ) ,  $$

$$  \mathfrak b_{1 , \,  \j }      ( \e x , \e y )  :=   \pi_{1 ,  \,  \j }  (  \mathcal B  ( \e x , \e y ) ) ,  $$

$$  \mathfrak b_{2 ,   \,  \j }      ( \edag  x , \edag  y )  :=   \pi_{2 ,    \,  \j }  (  \mathcal B  ( \edag x , \edag y ) ) ,  $$

which turn out to be Hermitian  complex  sesquilinear    forms  on   $X_{ \edag , \i } $,   $X_{ \e ,  \,  \j } $,    $X_{ \edag , \,  \j } $   respectively.


One can also obtain polarization   identities  for  matrices,  which directly descend from those we have described above. Indeed, any  matrix  $A \in \bc^{n \times n}$,   generates two  $\ast$--Hermitian  matrices
$$  A_1 :=  A  +  A^{\ast \, t }  \qquad  \quad  {\rm  and}  \qquad  \quad  A_2  :=  \i  \, ( A - A^{\ast \, t } ) .$$

It  is  easy to  see  that  the  formulas
$$  \mathcal B_1 (x , \, y )  : = x^t \cdot A_1 \cdot  y^{\ast  } $$

and
$$  \mathcal B_2 (x , \, y )  :=    x^t   \cdot A_2 \cdot  y^{\ast  } $$

define  $\ast$--Hermitian  and  $\ast$--sesquilinear  bicomplex  forms  on  the  bicomplex  module  $\bc^n$,  with
$ \mathcal Q_1 (x)  := \mathcal B_1 ( x , \, x)  =  x^t  \cdot A_1 \cdot  x^{\ast  } $   and      \;  $ \mathcal Q_2 (x)  := \mathcal B_2 ( x , \, x)  =  x^t  \cdot A_2 \cdot  x^{\ast  } $.   Thus  any  of  the  formulas   (\ref{1i})--(\ref{ik})  implies  the  corresponding  formula  for   $\mathcal B_1$,   $\mathcal Q_1$  and    $\mathcal B_2$,     $\mathcal Q_2$.    For  instance,  formula  (\ref{1i})  applied  to the    bicomplex  form  $ \mathcal B_2$  gives,  for  any  matrix  $A$,
\begin{equation}\label{formula}    \begin{array}{l}
\displaystyle  x^t \cdot  \left(   A - A^{ \ast \, t } \right) \cdot y^{\ast  }   =
\\   \\
\qquad    \displaystyle =   \frac{1}{4}  \,  \left(  (x + y )^t \cdot ( A - A^{ \ast \, t } ) \cdot  ( x + y )^{ \ast  }  \right.
\\   \\    \qquad  \qquad   -  \;     (x - y )^t   \cdot ( A - A^{ \ast \, t } ) \cdot  ( x - y )^{ \ast  }  \;  +
\\    \\
  + \,  \i  \,  \displaystyle  \left.   \left(  (x + \i  y )^t   \cdot ( A - A^{ \ast \, t } ) \cdot  ( x + \i  y )^{ \ast  }   -   (x - \i  y )^t    \cdot ( A - A^{ \ast \, t } ) \cdot  ( x -  \i  y )^{ \ast  }   \right)  \right).
\end{array}
\end{equation}

Assume  additionally  that  for   any  column     $b \in \bc^n$  the  matrix  $A$  verifies:
$$  b^t   \cdot A \cdot b^{\ast  }  \in \mathbb D ;$$

the  hyperbolic  numbers  seen  as  $1 \times 1 $  matrices  are  invariant  under  transposition  and  under  $\ast$--involution,  hence
$$ b^t  \cdot A \cdot b^{\ast } =  (  b^t   \cdot A \cdot b^{\ast }  )^{\ast \, t }  =    (     b^{ \ast \, t }    \cdot A^\ast   \cdot b )^t   =     b^t    \cdot A^{ \ast \, t }   \cdot b^{\ast  }  ,$$

which  implies  the  vanishing of the right--hand   side  of  (\ref{formula})  for  such  a  matrix.  Thus
$$x^t   \cdot   (  A - A^{ \ast \, t } )  \cdot  y^{ \ast }  =0  $$

for  any  $x , \, y $  and  so
$$  A   =  A^{ \ast \, t }  . $$

\medskip

\section{Linear operators on $\bc$--modules.}

Let  $X$  and  $Y$  be  two  $\bc$--modules. A  map    $T:  X  \to Y$  is said to be a  linear  operator  if  for  any   $x , \, z  \in X$  and  for  any  $\lambda \in \bc$  it is
$$  T [   \lambda x + z ]  = \lambda T [x]  + T [z]  .$$

Defining  as  usual  the  sum  of  two  operators and the multiplication of a linear  operator    by  a  bicomplex  number,  one  can see  that  the  set  of   all  linear  operators  forms  a   $\bc$--module.  In case when  $Y=X$  the  same  set  equipped  additionally  with  the composition  becomes  a  non--commutative   ring  and  an  algebra--type  $\bc$--module  (as before we  use  this  term  instead  of  the  term  ``algebra"  for  a ring  which is  a  module,  not  a  linear  space).

\medskip

Let us now see  how the idempotent decompositions of $X$ and $Y$  manifest themselves in the properties of a linear operator. Write  $X= X_\e + X_\edag $,  $Y = Y_\e + Y_\edag$ so that any   $x \in X$  is of the  form
$$  x = \e \cdot \e x + \edag \cdot \edag x =  :  \e \cdot x_1  + \edag \cdot  x_2 .$$

The linearity  of $T$ gives:
$$ \begin{array}{rcl}
T[x]  &=   &   T [  \e \cdot \e x + \edag \cdot \edag x  ]
\\  & & \\ & = &
\e \cdot T[\e x ]  +  \edag \cdot T [ \edag x ]
\\ & & \\ & = &
\displaystyle \e \cdot \left( \e \cdot T [ \e x ]  \right)  +   \edag \cdot \left( \edag \cdot T [ \edag x ]  \right) ,
\end{array}  $$

and introducing  the operators $T_1$ and $T_2$ by
$$ T_1 [x]  := \e \cdot T [\e x ] , \qquad   T_2 [x]  := \edag \cdot T [\edag x ]  $$

we obtain  what   can  be  called  the  idempotent  representation of a linear  operator:
$$  T= \e \cdot T_1  + \edag \cdot T_2 = T_1 + T_2 .$$

Although $T_1$ and $T_2$ are defined on the whole $X$  their ranges do not fill all of $Y$:  the  range  $T_1 (X) $ lies inside $Y_\e$  and  the  range  $T_2 (X) $ lies inside $Y_\edag$; in addition, the action of $T_1$ is determined by the component  $x_1 = \e x$ of $x$ and the action of $T_2$ is determined by the component  $x_2 = \edag x$ of $x$.  Thus, with a somewhat imprecise terminology, we can  say that    $T_1$ is an operator from $X_\e$ into $Y_\e$, and  $T_2$ is an operator from $X_\edag$ into $Y_\edag$. Since $T_1$ and $T_2$ are obviously  $\bc$--linear  then one may  consider them as $\C ( \i)$--  or  $\C ( \, \j)$--linear  on $X_\e$ and $X_\edag$  respectively.

\medskip

Let   $T: X  \longrightarrow Y$  be  a  $\bc$--linear operator. Assume   additionally  that  $X$  and  $Y$  have  $\D$--norms  $\| \cdot \|_{\D , \, X}$  and   $\| \cdot \|_{\D , \, Y}$
respectively (as well as the  bicomplex norms   $\| \cdot \|_X$  and   $\| \cdot \|_Y$).  Recall that  we can always write $X= X_\e + X_\edag$,  $Y = Y_\e + Y_\edag$,  that  the  norms   $\| \cdot \|_{X_\e}$,   $\| \cdot \|_{X_\edag}$  are the restrictions of  $\| \cdot \|_X$  to  $X_\e$  and  to  $  X_\edag$,  and   that      $\| \cdot \|_{Y_\e}$,   $\| \cdot \|_{Y_\edag}$  are the restrictions of  $\| \cdot \|_Y$   to  $Y_\e$  and  to  $  Y_\edag$.

\begin{subDn}
The operator  $T: X  \longrightarrow Y$  is called  $\D$--bounded  if  there  exists  $\Lambda \in \D^+$ such that for any $x \in X$  one has
\begin{equation}\label{def op acotado}
 \| Tx \|_{\D, \, Y}  \lessdot \Lambda \cdot  \| x \|_{\D , \, X} .
 \end{equation}

The  infimum  of  these  $\Lambda$ is  called  the  norm  of  the  operator  $T$.
\end{subDn}

It is clear  that the set  of  $\D$--bounded  $\bc$  operators  is  a  $\bc$--module.     Indeed,  let   $T_1$ and $T_2$  be   two   $\D$--bounded   operators:
$$ \displaystyle  \|  T_1 x \|_{\D ,Y}  \lessdot  \Lambda_1 \cdot  \| x \|_{\D,X} ,   \qquad    \|  T_2 x \|_{\D ,Y}  \lessdot  \Lambda_2 \cdot  \| x \|_{\D,X}   \forall x \in X ,$$
then we have
$$ \begin{array}{rcl}
\displaystyle  \|  \left( T_1  +  T_2 \right) [ x ]  \|_{\D ,Y} &  \lessdot   &  \|  T_1 x \|_{\D ,Y}   +   \|  T_2 x \|_{\D ,Y}
\\  & & \\ & \lessdot &
 \Lambda_1 \cdot  \| x \|_{\D,X}   +   \Lambda_2 \cdot  \| x \|_{\D,X} .
 \end{array}  $$
 If   $\Lambda_1$ and $\Lambda_2$  are comparable  with respect to the partial order $\lessdot$  then
$$ \displaystyle  \|  \left( T_1  +  T_2 \right) [ x ]  \|_{\D ,Y}   \lessdot    \max \{ \Lambda_1 , \, \Lambda_2 \} \cdot     \| x \|_{\D,X} . $$
But the  hyperbolic numbers  $\Lambda_1$ and $\Lambda_2$  may be  non comparable.  In  this  case  the following   bound  works.  If  $\Lambda_1 = a_1 \e + b_1 \edag$,    $\Lambda_2 = a_2 \e + b_2 \edag$  then  set   $a := \max \{ a_1 , a_2 \}$,   $b := \max \{ b_1 , b_2 \}$,   and   $\Lambda:=  a \e + b \edag$ which gives
$$ \displaystyle  \|  \left( T_1  +  T_2 \right) [ x ]  \|_{\D ,Y}   \lessdot     \Lambda  \cdot     \| x \|_{\D,X}   \;   \; \forall x \in X , $$
thus,  $T_1 + T_2$  is  $\D$--bounded.   The  case  of the multiplication  by bicomplex  scalars  is even  simpler  since
$$   \displaystyle    \|  \mu  Tx \|_{\D ,Y}   =  \| \mu \|_\k  \cdot  \| Tx\|_{\D, X}  \lessdot   \| \mu \|_\k  \cdot \Lambda  \cdot   \| x \|_{\D,X}   \;   \; \forall x \in X .$$

  \medskip

One can finally prove, following  the usual argument from  classical   functional  analysis,    that  the  map
$$  T  \mapsto  {\rm inf}_\D   \,   \{ \Lambda \mid  \Lambda \; {\rm satisfies \; eq.  (\ref{def op acotado}) }  \, \}  =:   \| T \|_\D $$

defines  a  hyperbolic    norm  on  the  $\bc$--module   of  all    bounded  operators.

\bigskip
  
Note that inequality  (\ref{def op acotado})  is  equivalent  to asking that, for  any  $x  \in X$,
$$  \begin{array}{rcl}
\| Tx\|_{ \D , Y} & = & \|  T_1 x_1 + T_2 x_2 \|_{\D . Y}
\\  &  & \\ & = &
\e \cdot  \| T_1 x_1 \|_{Y_\e }  +  \edag \cdot  \| T_2 x_2 \|_{Y_\edag }  \; \lessdot
\\ & & \\ &    \lessdot  &
\e  \cdot \Lambda_1 \cdot \| x_1 \|_{ X_\e }  +   \edag  \cdot \Lambda_2 \cdot \| x_2 \|_{ X_\edag }   ,
\end{array} $$

with  $\Lambda = \e \Lambda_1 + \edag \Lambda_2 \in \D^+$.  This means in particular that $T_1$ and $T_2$ are bounded and  $\|T_1 \| = \inf  \{  \Lambda_1 \} $,  $\|T_2 \| = \inf  \{  \Lambda_2 \} $. It is clear now that   if  $T_1$ and $T_2$ are bounded   then $T$ is $\D$--bounded.
Thus, the operator $T$ is   bounded if and only if the operators $T_1 : X_\e \to Y_\e$  and  $T_2 : X_\edag \to Y_\edag $  are  bounded.  What  is  more,  the  respective  norms are connected by equation
\begin{equation}\label{nov_19}
\| T \|_\D = \e \cdot \| T_1 \| + \edag \cdot \| T_2 \| ,
\end{equation}

which is true because 
$$ \begin{array}{rcl}
{\rm inf}_\D \{ \Lambda \} & = &  {\rm inf}_\D \{ \e \Lambda_1 + \edag \Lambda_2 \}
\\ & & \\ & = &
\e \cdot \inf \{ \Lambda_1 \} + \edag  \cdot   \inf  \{ \Lambda_2 \}  \, = \, \e \cdot \| T_1 \|  + \edag \| T_2 \| .
\end{array}  $$

As  always,  the  operator  $T$ is (sequentially) continuous at a point $x \in X$ if for any $\{ x_n \} \subset X$ such   that  $ \displaystyle    \lim_{ n \to \infty }    \| x_n - x \|_{ \D , X }  =0$  one has that
$ \displaystyle    \lim_{ n \to \infty }    \|  T x_n -   T  x \|_{ \D , Y  }  =0$.  It is obvious  that  $T$ is continuous at $x$ if and only if   $T_1$ is continuous at $ \e x$ and
$T_2$ is continuous at $ \edag x$.  This implies immediately  that  $T$  is   $\D$--bounded if and only if  $T$ is continuous.

\medskip

As in the classical case, the  norms of bounded  ope\-ra\-tors on complex  spaces can be represented in different ways:
\begin{equation}\label{star_Dec}
 \begin{array}{rcl}
\| T_1 \|  & = &   \displaystyle  \sup  \left\{  \, \frac{  \| T_1 z \|_{Y_\e } }{  \|  z \|_{X_\e } }  \mid  \| z \|_{X_\e }  \neq 0 \right\}
\\ & & \\ & = &
 \displaystyle  \sup  \left\{  \,   \| T_1 z \|_{Y_\e }    \mid  \| z \|_{X_\e }  \leq 1 \,   \right\}
\\ & & \\ & = &
 \displaystyle  \sup  \left\{  \,   \| T_1 z \|_{Y_\e }           \mid  \| z \|_{X_\e }  =1  \,   \right\}  ;
\end{array}
\end{equation}

(and similarly  for  $T_2$).

\medskip

By using these equalities, formula (\ref{nov_19})  gives
$$ \begin{array}{l}
\| T \|_\D \,  = \,    \displaystyle   \e \cdot  \sup  \left\{    \frac{ \| T_1 z \|_{Y_\e } }{  \|   z \|_{X_\e } }  \mid  z \in X_\e \setminus \{ 0 \}  \,    \right\} \; +
\\    \\
\qquad   \qquad \quad   + \;    \displaystyle   \edag   \cdot  \sup  \left\{    \frac{ \| T_2 u \|_{Y_\edag } }{  \|   u \|_{X_\edag } }  \mid  u \in X_\edag \setminus \{ 0 \}  \,    \right\}  \;  =
\\  \\
\quad   =  \,    \displaystyle {\rm sup}_\D  \,   \left\{    \frac{    \e \cdot   \| T_1 z \|_{Y_\e }     + \edag \cdot  \| T_2 u \|_{Y_\edag }      }{     \e \cdot    \|   z \|_{X_\e }    +  \edag \cdot  \| u \|_{X_\edag }    }    \mid        \e \| z \|_{ X_\e } + \edag \| u \|_{X_\edag} \not\in \mathfrak S_0 \right\} \;
\\      \\  \quad  = \,
  \displaystyle  {\rm sup}_\D \,   \left\{    \frac{     \|   \e \cdot   T_1 z      + \edag \cdot      T_2 u     \|_{\D , Y }         }{      \|    \e \cdot  z    +  \edag \cdot u    \|_{   \D ,  X }    }    \mid         \|   \e   z
   + \edag  u \|_{   \D ,  X} \not\in \mathfrak S_0 \right\} \;
   \\      \\ \quad   =  \,
  \displaystyle {\rm sup}_\D \,   \left\{    \frac{     \|    Tw     \|_{\D , Y }         }{      \|   w  \|_{   \D ,  X }    }    \mid         \|  w   \|_{   \D ,  X} \not\in \mathfrak S_0 \right\}  ,
\end{array}  $$

that is,  we  have  proved  for  $\D$--bounded  operators  an  exact  analogue  of  the  first  formula  in   (\ref{star_Dec}).

\medskip

Similarly, we can prove the  analogues  of  the   other two formulas  in   (\ref{star_Dec}).
$$   \begin{array}{rcl}
  \| T \|_\D   &  = &  \displaystyle {\rm sup}_\D \,   \left\{       \|    Tw     \|_{\D , Y }        \mid         \|  w   \|_{   \D ,  X}     \lessdot  1  \,   \right\}
  \\  &  &  \\  &  = &  \displaystyle {\rm sup}_\D \,   \left\{       \|    Tw     \|_{\D , Y }        \mid         \|  w   \|_{   \D ,  X}     =  1  \,   \right\} .
\end{array}  $$

\medskip

Assume now that  $X$ and $Y$  are bicomplex  Hilbert modules with  inner  products  $\langle \cdot , \cdot \rangle_X$  and   $\langle \cdot , \cdot \rangle_Y$  respectively.  Recall  that, with  obvious  meanings of the symbols,
$$  \langle x , z \rangle_X = \e \cdot \langle \e x , \e z \rangle_{X_\e}  +    \edag \cdot \langle \edag x , \edag z \rangle_{X_\edag}.  $$

The  adjoint   operator  $T^\sharp : Y \to X$  for  a  bounded  operator  $T : Y \to X$ is defined  by the  equality
$$  \langle T[x] , y \rangle_Y = \langle x , T^\sharp [y] \rangle_X  .$$

Consider the left--hand side:
$$  \begin{array}{rcl}
\langle T[x] , y \rangle_Y  & = &  \displaystyle    \langle \left(  T_1  + T_2 \right) \left[  \e x + \edag x \right], \e y  + \edag y \rangle_Y
\\ & & \\  & =  &
\displaystyle  \langle \e \cdot T_1 [ \e x ] + \edag \cdot  T_2 [ \edag x] , \e y + \edag y \rangle_Y
\\  & & \\ &= &
\displaystyle    \e \cdot \langle T_1 [ \e x ]   , \e y   \rangle_{Y_\e}    +      \edag   \cdot  \langle  T_2 [ \edag x]  , \edag y \rangle_{Y_\edag}
\\  & & \\ &= &
\displaystyle    \e \cdot \langle  \e x    ,  T_1^\sharp [ \e y ]     \rangle_{Y_\e}    +      \edag   \cdot  \langle   \edag x ,     T_2^\sharp  [ \edag y]   \rangle_{Y_\edag}  ;
\end{array}  $$

similarly for the right--hand side:
$$  \begin{array}{rcl}
\langle x ,   T^\sharp  [ y  ]   \rangle_X  & = &  \displaystyle    \langle \e x + \edag x ,    \left(  (T^\sharp)_1 +  (T^\sharp )_2  \right)  [ \e y  + \edag y  ]   \rangle_X
\\ & & \\  & =  &
\displaystyle   \e \cdot   \langle \e x  ,   (T^\sharp)_1     [ \e y  ]  \rangle_X      +    \edag \cdot  \langle   \edag x ,      (T^\sharp )_2 [   \edag y  ]   \rangle_X .
\end{array}  $$

Thus,  one concludes that the bicomplex adjoint  $T^\sharp$  always  exists  and its idempotent  components   $  (T^\sharp)_1 $  and  $ (T^\sharp)_2 $  are  the complex  adjoints of the operators $ T_1 $  and  $   T_2$, that is,
\begin{equation}\label{Tsharp}
T^\sharp = \e \cdot T_1^\sharp + \edag \cdot T_2^\sharp   .
\end{equation}

\bigskip


\chapter{Schur analysis}

This last chapter is an application of the previous results, and  it  begins
a study of Schur analysis in the bicomplex setting. Further
results, and in particular the associated theory of linear
systems, will be presented elsewhere.

\section{A  survey of classical Schur analysis}\label{section survey}
Under the term {\sl Schur analysis} one understands problems
associated to the class of functions $s$ that are holomorphic and whose modulus  is  bounded  by one in the open
unit disk $\mathbb K$ (we will call such functions Schur
functions, or Schur multipliers, and denote their class by
$\mathscr S$). They can be  interpreted as  the transfer functions of time-invariant
dissipative linear systems and, as such, play also an important
role in the theory of linear systems    (see \cite{MR2002b:47144}
for a survey).   It   was in the early years of  the   $XX^{th}$ century that I. Schur  wrote two particular works  \cite{schur1,schur2}  in 1917 and  1918; they were  motivated by  the works of Carath\'eodory and Fej\'er  \cite{CF}, Herglotz \cite{H} and Toeplitz \cite{T}  (in particular in  the trigonometric moment problem). In those works  Schur
associated to a Schur function  $s$   a sequence
$\rho_0,\rho_1,\ldots$ of complex numbers in the open unit disk:  one sets $s_0(z)=s(z)$, $\rho_0=s_0(0)$;  then  one defines
recursively a sequence of Schur functions $s_0,s_1,\ldots$ as
follows: if $|s_n(0)|=1$ the process ends. If $|s_n(0)|<1$, one
defines
\[
s_{n+1}(z)=\begin{cases}\,\,    \displaystyle      \frac{s_n(z)-s_n(0)}{z(1-\overline{s_n(0)}s_n(z))},\quad
{\rm if}\,\, z\not=0, \\  \\
\,\,   \displaystyle      \frac{s_n^\prime(0)}{(1-|s_n(0)|^2)},\quad\hspace{0.5cm}{\rm
if}\,\, z   =0.
\end{cases}
\]

Then $s_{n+1}$ is still a Schur function in view   of    Schwartz' lemma, and one sets $\rho_{ n + 1}    :  =  s_{ n +1}  (0)$. The coefficients $\rho_0,\rho_1,\ldots$
are called the Schur coefficients of $s$   and  they characterize in a unique way  the function $s$.   The construction of the sequence  $s_0$, $s_1, \ldots $  is  called  Schur's algorithm.    Since the function $s$  is  holomorphic, it  admits  its  Taylor   expansion:
$$   \displaystyle   s(z)  =  a_0  + a_1 \, z + a_2 \, z^2  + \cdots   +  a_n \, z^n  +  \cdots   ,$$

and the sequence of coefficients  $\{a_n\}$  also  characterize in a unique way  the function  $s$. It is remarkable   that  the sequence  of  Schur's coefficients  $ \{ \rho_n \}$    has  proved  to  give  a  better    characterization  of  $s$  and  makes possible  to  link this  type  of  functions  with, for instance,  problems  of  interpolation of  holomorphic functions,  filtration  of  stationary  stochastic  processes, etc.     Examples of  Schur functions  are the  so called  Blaschke  factors:
$$  \displaystyle   b_\omega (z)  =  \frac{  z - \omega }{ 1 - z \overline{\omega } } ,$$

with  $\omega \in \mathbb K$,  and more   general  finite  products:
\begin{equation}\label{Blaschke factors product}
s(z) =  c \cdot  \prod_{\ell = 1}^k  b_{\omega_\ell} ,
\end{equation}

where $c$ is a complex number with modulus   one and  $\omega_1, \omega_2 , \ldots , \omega_k$  are in $\mathbb K$.       I. Schur  proved that  the  sequence of  coefficients   $ \{ \rho_n \}$   is finite if and  only  if   $s$ is   a  function of the  type   (\ref{Blaschke factors product}), i.e., if $s$  is   a finite Blaschke
product.\\

The role of Schur functions in the theory of time-invariant
dissipative linear systems has led to the extension of the definition of
Schur functions to other settings.   For instance, the case of several complex
variables corresponds to multi-indexed ($ND$) systems (recall that a linear  system can be seen as   a  holomorphic or a rational function defined in   the unit  disk  and  whose values  are  matrices). Analogues
of Schur functions were introduced by Agler, see
\cite{agler-hellinger,agler-carthy}. Besides the Schur-Agler classes, we
mention that counterparts of Schur functions have been studied in
the time-varying case, \cite{MR93b:47027}, the stochastic case,
\cite{alpay-2008},  the Riemann surface case \cite{av3}, and the
multiscale case, \cite{a_m,a_m_ieot}.   A number of other cases also
are of importance. In the setting of hypercomplex analysis,
counterparts of Schur functions were studied in
\cite{asv-cras,MR2124899,MR2240272,MR2275397} in the setting of Fueter
series, and in \cite{acs1} in the setting of slice
hyperholomorphic functions. In this chapter we study the
analogue of Schur functions in the setting of bicomplex numbers
and $\bc$-holomorphic functions. \\

We begin with the following characterization of Schur functions
of  one  complex  variable.

\begin{subtheorem}
\label{thmSchurdisk}
Let $s$ be a function holomorphic in the unit disk $\mathbb K$. The following are equivalent:\\

\begin{enumerate}
\item[{\rm (i)}]    $s$ is a Schur function.\\

\item[{\rm (ii)}]   The kernel
\begin{equation}
\label{ks}
k_s(\lambda,\mu)=\frac{1-  s(\lambda)   \,  \overline{
s(\mu)}   }{1-\lambda  \, \overline{ \mu} }
\end{equation}
is positive definite in $\mathbb K$,   that  is,  given  any  collection   $\lambda_1 ,  \ldots , \lambda_n$    of  points in  $\mathbb K$,  the  matrix    $ \displaystyle   \left(   k_s  (  \lambda_\ell  , \, \lambda_m  )  \right)_{ \ell , \, m }  $  is  positive  definite.      \\

\item[{\rm (iii)}]  One can write $s$ in the form
\begin{equation}
\label{eq:real} s(\lambda)=D+\lambda C(I_{\mathcal H (s) }-\lambda
A)^{-1}B,
\end{equation}
where $\mathcal H (s) $ is a Hilbert space and where the operator
matrix
\[
\begin{pmatrix}
A&B\\ C&D\end{pmatrix}\,\,\, :\,\,\, \mathcal H (s)  \oplus   \mathbb
C\quad\longrightarrow \quad\mathcal H (s)  \oplus\mathbb C
\]
is coisometric    (i.e.,  $\displaystyle  \begin{pmatrix}
A&B\\ C&D\end{pmatrix}  \cdot    \begin{pmatrix}
A &  B   \\   C  &   D    \end{pmatrix}^\sharp    = \,    I_{\mathcal H (s)  \oplus  \C  } $,  where we are using the symbol $\sharp$  instead of $\ast$  to denote the dual of the operator matrix)   or unitary.  \\

\item[{\rm (iii)}] The multiplication operator    $M_s$  defined  by   $M_s [  f  ]  : =  sf$ is a contraction
from the Hardy space $\mathbf H^2(\mathbb K)$ into itself.
\end{enumerate}
\end{subtheorem}

Expression \eqref{eq:real} is called a realization of $s$.  Here we need to recall  that  a   Hilbert space  $\mathcal H$  whose elements are functions  defined on   a set $\Omega$  and with  values in $\C^n$   is called  a  reproducing kernel  Hilbert  space  if   for  any  $ c \in \C^n$  and  any  $\omega \in \Omega$   the   functional       $f  \mapsto  \overline{c}^{\, t}  f(\omega)$  is bounded;   by  Riesz' Theorem  (see for  instance   \cite{Rudin} Theorem 4.12, p. 77)  there exists a function  
$$ K (z , \omega ):    \Omega \times \Omega \to \C^{n \times n } $$

such that  for  any  $ c \in \C^n$,    $\omega \in \Omega$  and $f \in \mathcal H$, the function   $ z \mapsto  K (z, \omega) \cdot c$  belongs  to  $\mathcal H$  and   
\begin{equation}\label{reproducing property}
\langle f , \, K ( \cdot , \omega) \, c \, \rangle_{\mathcal H} =  \overline{c}^{\, t}    f(\omega) ,
\end{equation}

where   $\langle \cdot , \cdot \rangle_{\mathcal H}$  is the  scalar product in the  Hilbert space $\mathcal H$;  it   is the equation (\ref{reproducing property})      that  justifies  the term ``reproducing kernel"; it is proved that the function  $K(z, \omega)$ is unique and is called the  reproducing kernel of the  Hilbert space $\mathcal H$; it is direct to  prove that   the function $K$ is  positive  definite  and  also  is    Hermitian (in the usual complex  setting), that is:
$$ \overline{  K(z , \omega)  }   =  K  (\omega ,z) .$$

We need to recall also that  the reciprocal  of the  above  is also true, that is,    given  a  positive definite function  $K(z, \omega): \Omega \times \Omega \to  \C^{n \times n}$ there exists a unique Hilbert space  whose elements  are  functions   from $\Omega$ to $\C^n$  and  with  reproducing kernel  $K$.  An important example for the case $n=1$   that we will use later is the function
\begin{equation}\label{example positive function}
\displaystyle  K(z, w) :=   \frac{1}{1 - z \overline{w} }  \quad  z , w \in \mathbb K,
\end{equation}  
which is a positive definite function and it is known that  its  re\-pro\-du\-cing kernel Hilbert space  is  the Hardy space  ${\mathbf H}^2 (\mathbb K)$.  \\

In the previous Theorem we have the particular case $n=1$ and we  have  denoted by $\mathcal H(s)$ the reproducing kernel Hilbert space of
functions holomorphic in $\mathbb K$ and with reproducing kernel
$k_s(\lambda,\mu)$. This space provides a coisometric
realization of $s$ called the backward shift realization
of $s$, which is already present (maybe in an implicit way) in the
work of de Branges and Rovnyak \cite{dbr2}. See for instance
Problem 47, p. 29 in that book. 

\medskip

\begin{subtheorem}
The operators
\begin{equation}
\begin{split}
A [  f  ]  (z)&   : =\frac{f(z)-f(0)}{z},\\  \\  
(B  [  c  ]  )(z)&  : =  \frac{s(z)-s(0)}{z}c,\\   \\ 
C  [  f  ]  &  : =   f(0),\\    \\
D  [  c  ]   &=s(0)c,
\end{split}
\end{equation}
where $f\in\mathcal H(s)$ and $c\in \mathbb C$, define a
coisometric realization of $s$.
\end{subtheorem}

Extending the notion of Schur function to the case of several
complex va\-ria\-bles is not so simple. We will focus on the cases
of the polydisk and of the ball. We  collect here in two theorems, \cite{agler-hellinger},
the analogues of Theorem \ref{thmSchurdisk} for the polydisk and the
ball. 

\medskip

\begin{subtheorem}
Let $s$ be holomorphic in the polydisk $\mathbb K^N$, let
$\la=(\la_1,\ldots, \la_N)   \in \C^N$ and let $\Lambda={\rm
diag}~(\lambda_1,\lambda_2,\ldots, \lambda_N)$.
Then, the fo\-llo\-wing are equivalent:\\

\begin{enumerate}

\item[{\rm (i)}]  There exist   $\C$--valued   functions   $k_1(\lambda,\mu),\ldots,
k_N(\lambda, \mu)$,  positive definite in $\mathbb K^N$ and such
that
\begin{equation}
\label{schur-agler1}
1-s(\lambda)  \,   \overline{ s(\mu)}
=\sum_{n=1}^N(1-\la_n   \overline{ \mu}_n  )k_n(\la, \mu).
\end{equation}

\item[{\rm (ii)}]   There exist Hilbert spaces ${\mathcal H}_1,\ldots
,\mathcal H_N$ and a unitary operator matrix
\[
\begin{pmatrix}
A&B\\ C&D\end{pmatrix}\,\,\, :\,\,\, \mathcal H\oplus\mathbb
C\quad\longrightarrow \quad\mathcal H\oplus\mathbb C
\]
such that
\[
s(\la)=D+C(I_{\mathcal H}-\Lambda A)^{-1}\Lambda B  ,
\]
where $\mathcal H=\oplus_{n=1}^N\mathcal H_n$.\\
\end{enumerate}
\end{subtheorem}

When $N=1$, $k_1$ is uniquely determined, and is in fact the function $k_s$
introduced in Theorem \ref{thmSchurdisk}. Note that in general the decomposition \eqref{schur-agler1} is not unique.
Dividing both sides of \eqref{schur-agler1} by $  \displaystyle
\prod_{n=1}^N   (1-\la_n    \,   \overline{\mu}_n)$ we see that
in particular $M_s$ is a contraction from the Hardy space of the
polydisk $\mathbf H^2(\mathbb K^N)$ into itself. When $N>2$ the
class of Schur-Agler functions of the polydisk is strictly
included in the class of functions holomorphic and contractive in
the polydisk, or equivalently such that  the multiplication operator   $M_s [ f ]   = s  f $  is a contraction
from  the Hardy space 
$\mathbf H^2(\mathbb K^N)$ into itself. The two classes
coincide when $N=1$ or $N=2$.\\

We now consider the case of the open unit ball $\mathbb B_N$ of
$\mathbb C^N$. The function
\[
\displaystyle   \frac{1}{1-\la \, \overline{
\mu}}=\frac{1}{1-    \displaystyle  \sum_{n=1}^N\la_n   \,   \overline{\mu}_n}
\]
is positive-definite there. Its associated reproducing kernel
Hilbert space is called the Drury-Arveson space, and will be
denoted by $\mathcal A$. The space $\mathcal A$  is strictly
included in the Hardy space $\mathbf H^2(\mathbb B_N)$ of the
ball when $N>1$, and  the class of Schur-Agler functions of the
ball is strictly included in the class of functions holomorphic and
contractive in $\mathbb B_N$, or equivalently, in the class of
functions holomorphic in $\mathbb B_N$ and such that $M_s$ is a
contraction from $\mathbf H^2(\mathbb B_N)$ into itself.

\begin{subtheorem}
\label{schur-agler-ball}
Let $s$ be holomorphic in $\mathbb B_N$. The following are
equi\-va\-lent:\\

\begin{enumerate}

\item[{\rm (i)}]  The kernel
\[
\frac{1-s(\la)   \, \overline{ s(\mu)}  }{1-\la \,
\overline{ \mu} }
\]
is positive definite in $\mathbb B_N$.\\

\item[{\rm (ii)}]  There exists a Hilbert space ${\mathcal H}$ and a unitary
operator matrix
\[
\begin{pmatrix}
A&B\\ C&D\end{pmatrix}\,\,\, :\,\,\, \mathcal H    \oplus   \C   \quad\longrightarrow \quad\mathcal H  \oplus\mathbb C
\]
such that
\[
s(\la)=D+C(I_{\mathcal H}-\lambda A)^{-1}\lambda B  .
\]

\item[{\rm (iii)}]   The operator $M_s$ is a contraction from the
Drury-Arveson space into itself.
\end{enumerate}
\end{subtheorem}



\section{The bicomplex Hardy space}\label{bicomplex Hardy space}


We now want to show how the ideas   and concepts of  Section \ref{section survey}   can be translated to the case of bicomplex analysis, and bicomplex
holomorphic functions. Let
\begin{equation}
\label{ok2}
\Omega_{\mathbb K^2}   :=\left\{Z=z_1+z_2\mathbf
j   =  \e \beta_1 + \edag \beta_2     \,\,|  \, \,   (\beta_1,\beta_2)\in\mathbb K^2\right\}  ,
\end{equation}

with  $\mathbb K$, as in the previous section,   being  the unit disk in the complex  plane  and  $\mathbb K^2 =  \mathbb K  \times \mathbb K  $.  The reader may consider  redundant   the notation   $ \Omega_{\mathbb K^2}  $, but recall  that  there  are many different ways of  identify the set $\bc$ with $\C^2$  and we want to  make emphasis  that  now   we are  making  this identification via the  idempotent  representation.  \\

           We define   the  bicomplex  Hardy  space    $\mathbf H^2 (\Omega_{\mathbb K^2} )$ to be  the set of functions  $f: \Omega_{\mathbb K^2}  \to  \bc$  such  that  for  any  $Z  \in   \Omega_{\mathbb K^2} $:
$$\displaystyle  f(Z)=   \sum_{ n=0}^\infty  f_n \cdot Z^n,$$

where  $\forall n \in \N$,     $f_n \in \bc$   and  the series   of  hyperbolic  numbers     $\displaystyle \sum_{ n=0}^\infty  | f_n |_\k^2 $  is  convergent.  Setting  $f_n =  \e   f_{n1} +   \edag     f_{n2}$  one gets:
$$ | f_n |_\k^2 =  \e    |  f_{n1}  |^2 +     \edag     |  f_{n2}  |^2 ,$$

which means  that  both  complex  series
$$\displaystyle \sum_{ n=0}^\infty  | f_{n1} |^2  ,  \qquad \quad  \sum_{ n=0}^\infty  | f_{n2} |^2   $$

are convergent and thus both  functions
$$\displaystyle    f_1 (\beta_1)  :=    \sum_{ n=0}^\infty   f_{n1}     \beta_1^n  ,  \qquad \quad    f_2(\beta_2):=   \sum_{ n=0}^\infty   f_{n2}     \beta_2^n   $$

belong to the  Hardy  space of the unit  disk  ${\mathbf H}^2 (\mathbb K)$;   of  course
$$  f (Z)  =  \e  \cdot  f_1 (\beta_1)+  \edag    \cdot  f_2 (\beta_2) ,$$

or equivalently,
$$  f (z_1 + \j   z_2)  = \e   \cdot  f_1 (z_1 - \i z_2  ) +  \edag   \cdot  f_2 (z_1 + \i  z_2)  .$$

This means that the bicomplex  Hardy space  can be written as 
\begin{equation}\label{bc Hardy space}
\mathbf H^2 (\Omega_{\mathbb K^2} )  =  \e \cdot  \mathbf H^2 (\mathbb K )   +   \edag \cdot \mathbf H^2 (\mathbb K ).
\end{equation}

The  $\bc$--valued  inner  product   on  $\mathbf H^2 (\Omega_{\mathbb K^2} )$     is  given  by
$$ \begin{array}{rcl}
\langle  f , g \rangle_{  \mathbf H^2 (\Omega_{\mathbb K^2} ) }   &   :=  &   \displaystyle   \sum_{n=0}^\infty  f_n \cdot g_n^\ast
\\  & & \\ & = &
 \displaystyle   \e     \cdot   \sum_{n=0}^\infty  f_{n 1 }   \cdot   \overline{ g_{n 1} }   +   \edag     \cdot   \sum_{n=0}^\infty  f_{n 2 }   \cdot   \overline{ g_{n 2} } 
\\  & & \\ & = &
 \e \cdot   \langle  f_1  , g_1    \rangle_{  \mathbf H^2 ( \mathbb K ) }  + \edag  \cdot    \langle  f_2  , g_2    \rangle_{  \mathbf H^2 ( \mathbb K ) }    ,
\end{array}$$

and  it  generates  the  $\D$--valued   norm:
$$  \begin{array}{rcl}
\| f \|^2_{ \D , \,    \mathbf H^2 (\Omega_{\mathbb K^2} ) }  &  =  &     \langle  f , f \rangle_{  \mathbf H2 (\Omega_{\mathbb K^2} ) }
\; = \;  
\displaystyle    \sum_{n=0}^\infty    |   f_n   |^2_\k
\\  & &  \\  & = &
\displaystyle  \e   \cdot   \sum_{n=0}^\infty    |   f_{n  1}   |^2   + \edag      \cdot   \sum_{n=0}^\infty    |   f_{n  2}   |^2 
\\  & &  \\  & = &
\e    \cdot  \| f_1   \|^2_{   \mathbf H^2 ( \mathbb K) }   +   \edag    \cdot  \| f_2   \|^2_{   \mathbf H2 ( \mathbb K) }  .
\end{array}  $$

Thus one can say  that  we  have  obtained  an analogue  of the classic  Hardy space on the unit disk for the bicomplex setting where we deal with the bidisk,   but  now   in the idempotent,  not cartesian,    coordinates;  the inner  product in this Hardy space is $\bc$--valued  and  (what is even more remarkable)  the ``adequate"  norm  takes values in $\D^+$,  the positive hyperbolic  numbers.

\bigskip

\section{Positive definite functions}
\label{sec:rkhs}


 Let $\Omega$ be some  set. The $\bc ^{n\times n}$-valued function
 $K(z,w)$ defined for $z,w\in\Omega$ is said to be positive
 definite if it is $\ast$--Hermitian:
 \[
 K(z,w)=K(w,z)^{\ast \, t},\quad\forall z,w\in\Omega,
 \]
 and if for every choice of $N\in\mathbb N$, of  columns  $c_1,\ldots ,
 c_N\in\mathbb C^n (\i)    $ and      of   $z_1,\ldots, z_N\in\Omega$ one has  that  
 \[
 \sum_{\ell,  \,  j=1}^N c_j^{*t}  \,  K(z_j,z_\ell)   \, c_\ell
 \in \mathbb D^+.
 \]

 \medskip

 \begin{subPn}
The $\bc ^{n\times n}$-valued function
 \[
 K(z,w)=K_1(z,w)+\mathbf j K_2(z,w)
 \]
 is positive definite if and only if $K_1(z,w)$ is positive
 definite and $\mathbf i K_2(z,w)$ is Hermitian and such that
 \begin{equation}
 -K_1(z,w)\le \mathbf i K_2(z,w)\le K_1(z,w).
 \end{equation}
 \end{subPn}
 {\bf Proof:} This follows from \eqref{ineq} applied to the matrices
 \[
 \begin{pmatrix}K(z_1,z_1)&K(z_1,z_2)&\cdots& K(z_1,z_N)\\
 K(z_2,z_1)&K(z_2,z_2)&\cdots& K(z_2,z_N)\\
 \vdots&\vdots & &\vdots\\
 K(z_N,z_1)&K(z_N,z_2)&\cdots& K(z_N,z_N)
 \end{pmatrix} .
 \]
\mbox{}\qed\mbox{}\\

\medskip

\begin{subCy}
\label{cor-positive}
Let $K(z,w)$ be $\bc ^{n\times
n}$-valued function defined in $\Omega ^2$. Write
$K(z,w)=k_1(z,w)\e+k_2(z,w)\edag$, where $k_1$ and $k_2$ are
$\mathbb C^{n\times n}$-valued functions. Then, $K$ is
positive-definite in $\Omega$ if and only if the functions $k_1$
and $k_2$ are positive definite in $\Omega$.
\end{subCy}

\begin{subtheorem}
\label{rkhsbc}
Let $K$ be a $\bc^{n\times n}$-valued function
positive definite on  the  set    $\Omega $. There exists a unique
reproducing kernel Hilbert space   $\mathscr H (K)$     of functions defined on $\Omega$
with reproducing kernel $K$. Namely, for all $w\in\Omega$,
$c\in\bc^n$ and
$f\in\mathscr H(K)$ one has that\\

\begin{enumerate}

\item[(i)]      The function $z\mapsto K(z,w)c$ belongs to $\mathscr H(K)$,
and

\item[(ii)]  The reproducing kernel property
\[
\langle f, K(\cdot, w)c\rangle_{\mathscr H(K)}=c^{\ast \, t}  f(w)
\]
holds.
\end{enumerate}
\end{subtheorem}

{\sf Proof:}  \\   Given $K=K_1+\j K_2  = \e   \cdot ( K_1 - \i  K_2 )  + \edag   \cdot  ( K_1 + \i K_2 ) $, a $\bc^{n\times n}$--valued
positive
 definite function on a   set   $\Omega$, the associated reproducing
 kernel Hilbert space consists of the functions of the form
\[
 F(z)=   f_1 (z)  \,  \mathbf e+    f_2(z)  \, \mathbf e^\dag,
 \]
 where $f_1   \in\mathscr H(K_1-\i K_2)$ and $f_2   \in\mathscr H(K_1+\i
 K_2)$, and with $\bc$--valued inner product
 \begin{equation}
 \label{eq:form}
 \langle  F,G  \rangle_{\mathscr H (K) }      : =\mathbf e\,\langle f_1,  g_1  \rangle_{\mathscr H(K_1-\i K_2)}+
 \mathbf e^\dag\,\langle f_2  ,   g_2   \rangle_{\mathscr H(K_1+\i K_2)}.
 \end{equation}
\mbox{}\qed\mbox{}\\

In the previous theorem, the $\mathbb C^{n\times n}$-valued
kernels are in particular po\-si\-ti\-ve definite
 in $\Omega$. We denote  by $\mathscr H(K_1-\i K_2)$ and $\mathscr H(K_1+\i K_2)$ the associated reproducing kernel
 Hilbert spaces of $\mathbb C^n$-valued functions on $\Omega$.\\

\begin{subCy}
The function   $K$ is positive definite in $\Omega$ if and only if there is a
$\mathbb{BC}$-Hilbert space $\mathcal H$ and a function $f$ from
$\Omega$ into $\mathcal H$ such that
\begin{equation}
K(z,w)=\langle f(w),f(z)\rangle_{\mathcal H} \quad \forall \; z , w  \in \Omega.
\end{equation}
\end{subCy}

{\sf Proof:}   \\    Write $K(z,w)=k_1(z,w)\e+k_2(z,w)\edag$, where $k_1$
and $k_2$ are complex-valued positive definite functions. Then
there exist complex Hilbert spaces $\mathcal H_1$ and $\mathcal
H_2$ and functions $f_1$ and $f_2$, defined on $\Omega$ and with
values in $\mathcal H_1$ and $\mathcal H_2$ respectively, such
that for any  $z, \, w \in \Omega$ there holds:
\[
k_1(z,w)=\langle f_1(w),f_1(z)\rangle_{\mathcal H_1}\quad {\rm
and} \quad k_2(z,w)=\langle f_2(w),f_2(z)\rangle_{\mathcal H_2}.
\]
Then using the process described in Section \ref{construction}, define  $\mathcal H$ to be
\[
\mathcal H=  \e  \cdot    \mathcal H_1+ \edag  \cdot    \mathcal H_2  .
\]

Applying to this  $\bc$--module  the results  from Subsection   \ref{inner product on the sum of two complex spaces} and in particular  the formula (\ref{avestruz}),  the inner product  on $\mathcal H$ is given by
$$ \begin{array}{rcl}
\langle  h , g \rangle_{\mathcal H}   &  =  &        \langle   \e \, h_1 + \edag \, h_2 , \,   \e \, g_1 + \edag \, g_2  \rangle_{\mathcal H}    
\\  &  &   \\   & = &
\e \,   \langle    h_1  ,  \,   g_1   \rangle_{\mathcal H_1}    \; + \;     \edag \,     \langle    h_2 , \,    g_2  \rangle_{\mathcal H_2}  ,
\end{array}   
$$  

and  for  any   $ z \in \Omega$:
\[
f(z)  :  =   \e \cdot      f_1(z)   + \edag \cdot  f_2(z)  ;
\]

thus one has:

$$  \begin{array}{rcl}
\langle f(w),f(z)\rangle_{\mathcal H}   & = &   \langle \e \,  f_1(w)    +  \edag \, f_2 (w)  ,   \e \,  f_1  (z)   + \edag  f_2 (z) \rangle_{\mathcal H} 
\\  &  &   \\  & = & 
\e \,  \langle f_1(w),f_1(z)\rangle_{\mathcal H_1}   + \edag   \langle f_2(w),f_2(z)\rangle_{\mathcal H_2}
\\  &  &   \\  & = &
\e \, k_1 (z,w)  + \edag \, k_2 (z ,w)   \;  =  \;   K  (z,w)    .
\end{array}
$$

This   concludes  the proof.
\mbox{}\qed\mbox{}\\

\medskip

 In the complex case, there is a one-to-one correspondence between
 po\-si\-ti\-ve definite functions on a set $\Omega$ and reproducing
 kernel Hilbert spaces of functions defined on $\Omega$. In the
 present case, one has the following result:

 \begin{subtheorem} There is a one-to-one
 correspondence between $\bc^{n\times n}$--valued positive definite
 functions on a set $\Omega $ and reproducing kernel
 Hilbert spaces of $\bc^n$--valued functions defined on $\Omega $.\\
 \end{subtheorem}

We  conclude this section with an example.

\begin{subtheorem}\mbox{}

\begin{enumerate}

\item[(1)] The function $  \displaystyle  \frac{1}{1-ZW^\ast}$ is positive definite in
$\Omega_{\mathbb K^2}$.
\mbox{}\\

\item[(2)]   Its associated $\bc$ reproducing kernel Hilbert space   is    $\bf
H^2(\Omega_{\mathbb K^2})$.   \\

\item[(3)]  For every $f\in\bf H^2 (\Omega_{\mathbb K^2}  )  $, and every $a\in \Omega_{\mathbb
K^2}$, we have:
 \[
 \langle  f   ,   \,  (1-Za^*)^{-1}   \rangle_{ \bf H^2 ( \Omega_{\mathbb K^2} ) }   =f(a)   .
 \]

 \end{enumerate}
 \end{subtheorem}

 {\sf Proof:}\\
 $(1)$ Writing $Z=\beta_1\e+\mu_1\edag$ and
$W=\beta_2\e+\mu_2\edag$, where $\beta_1,   \beta_2,   \mu_1$,   $\mu_2$ are
in $\mathbb K$,    we have:
\[
\frac{1}{1-ZW^*}=\frac{1}{1-\beta_1  \overline{\mu}_1}\e+
\frac{1}{1-\beta_2    \overline{\mu}_2}\edag.
\]
The claim follows then from Corollary \ref{cor-positive}.\\

Items  (2) and  (3)    are  consequences of Corollary \ref{rkhsbc},  and  of   (\ref{bc Hardy space}). 
\mbox{}\qed\mbox{}\\

\bigskip


\section{Schur multipliers and their characterizations in the
$\bc$-case}

\begin{subDn}
A $\bc^{n\times m}$-valued    $\bc$-holomorphic function   $s$   on
$\Omega_{\mathbb K^2}$      is called a Schur function if it can be
written as
\[
s(Z)=s_1(\beta_1)\e+s_2(\beta_2)\edag,
\]
where $s_1$ and $s_2$ belong to $\mathscr S^{n\times m}(\mathbb
K)$.
\end{subDn}

We will use the notation $\mathscr S^{n\times m}(\Omega_{\mathbb
K^2})$ for these   Schur   functions $s$.

\begin{subtheorem}
Let $s$ be a $\bc^{n\times m}$-valued     $\bc$-holomorphic   function 
 on $\Omega_{\mathbb K^2}$. Then, the following
are equivalent:\\

\begin{enumerate}

\item[(1)]    $s$ belongs to $\mathscr S^{n\times m}(\Omega_{\mathbb
K^2})$.\\

\item[(2)]
\[
s(Z)s(Z)^{*t}\le I_n,\quad \forall Z  \in\Omega_{\mathbb K^2}.
\]

\medskip

\item[(3)]    The function
\begin{equation}
\label{ksbc}
\frac{I_n-s(Z)s(W)^{*t}}{1-ZW^*}
\end{equation}
is positive definite in $\Omega_{\mathbb K^2}$.\\

\item[(4)]    The \, operator  \,   of multiplication by $s$ is a $\D$--contraction from
$(\mathbf H^2(\Omega_{\mathbb K^2}))^m$ into $(\mathbf
H^2(\Omega_{\mathbb K^2}))^n$, meaning that

\[
\langle sf,sf\rangle_{(\mathbf H^2(\Omega_{\mathbb K^2})  )^n}     \lessdot
\langle f,f\rangle_{(\mathbf H^2(\Omega_{\mathbb
K^2}))^m},\quad\forall f\in(\mathbf H^2(\Omega_{\mathbb K^2}))^m.
\]

\medskip

\item[(5)]   $s$ admits a realization
\begin{equation}
\label{reasbc}
s(Z)=D+ZC(I_{\mathcal H}-Z  A)^{-1}B
\end{equation}
where $\mathcal H$ is a $\bc$-Hilbert space and the operator
matrix
\[
\begin{pmatrix}
A&B\\ C&D\end{pmatrix}\,\,:\,\,\mathcal H\oplus
\bc^m\,\,\longrightarrow \mathcal H\oplus \bc^n
\]
is coisometric.

\end{enumerate}
\end{subtheorem}

{\sf Proof:}  \\    Assume that $(1)$ holds. Then, with
$W=\gamma_1\e+\gamma_2\edag$, we can write
\[
\begin{split}
I_n-s(Z)s(W)^{*t}&=(I_n-s_1(\beta_1)s_1(\gamma_1)^*)\e
+(I_n-s_2(\beta_2)s_2  (\gamma_2)^*)\edag.
\end{split}
\]
Since $s_1$ and $s_2$ are classical Schur functions we have:
\[
I_n-s_1(\beta_1)s_1(\gamma_1)^*\ge 0\quad{\rm and}\quad
I_n-s_2(\beta_2)s_2  (\gamma_2)^*\ge 0
\]

for all $\beta_1,\beta_2\in\mathbb K$. From Proposition
\ref{pnmat} we get that $(2)$ holds. When $(2)$ holds
both the kernels
\begin{equation}
\label{kern12}
\frac{I_n-s_1(\beta_1)  \overline{ s_1(\gamma_1) }^{\, t } }{1-\beta_1
\overline{ \gamma_1}   }   \quad{\rm and}
\quad\frac{I_n-s_2(\beta_2)   \overline{ s_2 (\gamma_2)}^{\, t} }{1-\beta_2   \overline{ \gamma_2}}
\end{equation}

are positive definite in the open unit disk. Writing
\[
\begin{split}
\frac{I_n-s(Z)   s(W)^{*t}}{1-ZW^*}   &       =\frac{I_n-s_1(\beta_1)  \overline{ s_1 (\gamma_1)}^{\, t} }{1-\beta_1
\overline{ \gamma_1}   }    \e + \frac{I_n-s_2(\beta_2)   \overline{  s_2 (\gamma_2)}^{\, t}  }{1-\beta_2
\overline{  \gamma_2} }\edag.
\end{split}
\]

it follows from Corollary \ref{cor-positive} that the kernel
\eqref{ksbc} is positive definite in $\Omega_{\mathbb K^2}$. Thus
$(3)$ holds. When $(3)$ holds, the kernels \eqref{kern12} are in
particular positive definite in the open unit disk. Thus the
operators of multiplication by $s_1$ and $s_2$ are contractions
from $(\mathbf H^2(\mathbb K))^m$ into $(\mathbf H^2(\mathbb
K))^n$, that is
\[
\langle s_j   f_j,   s_j     f_j\rangle_{(\mathbf H^2(\mathbb K))^n}\le
\langle f_j   ,   f_j    \rangle_{(\mathbf H^2(\mathbb K))^m},\quad j=1,2.
\]

Thus
$$ \begin{array}{rcl}
\langle sf , sf \rangle_{(\mathbf H^2(   \Omega_{\mathbb K^2} )  )^n}   & =  &
\langle s_1f_1,s_1f_1\rangle_{(\mathbf H^2(\mathbb K))^n}\e+
\langle s_2f_2,s_2f_2\rangle_{(\mathbf H^2(\mathbb K))^n}\edag \; \lessdot
\\  &  &   \\   & \lessdot  &
\langle f_1,f_1\rangle_{(\mathbf H^2(\mathbb K))^m}\e+ \langle
f_2,f_2\rangle_{(\mathbf H^2(\mathbb K))^m}\edag
\\   &   &  \\   & = &
\langle f , f \rangle_{(\mathbf H^2(   \Omega_{\mathbb K^2} )  )^m}    ,
\end{array}  $$

that is,     $(4)$ holds. Assume now that $(4)$ holds. By definition
of the inequality, $s_1$ and $s_2$ are classical Schur functions.
Using Theorem \ref{thmSchurdisk} we can write
\[
s_j(\beta_j)=D_j+\beta_jC_j(I_{\mathscr H_j}-\beta_j
A_j)^{-1}B_j,\quad j=1,2,
\]

where in these expressions, $\mathscr H_1$ and $\mathscr H_2$ are
complex Hilbert spaces, and the operator matrices
\[
\begin{pmatrix}
A_j&B_j\\ C_j&D_j\end{pmatrix}\,\,:\,\,\mathcal H_j\oplus \mathbb C
^m\,\,\longrightarrow \mathcal H_j\oplus \mathbb C^n
\]
are coisometric. Realization \eqref{reasbc} follows by taking   into  account  again  the  process  described  in  Section  \ref{construction}  and  setting
\[
\mathcal H  :     =    \e  \cdot  \mathcal H_1+     \edag    \mathcal
H_2   .
\]

Then the operators
\[
\begin{split}
A&: =   \e  A_1    +    \edag   A_2   ,    \\
B&  : =     \e  B_1      +    \edag   B_2   ,   \\
C&   : =   \e     C_1  +   \edag  C_2  ,   \\
D  &   : =    \e   D_1    +     \edag  D_2  ,
\end{split}
\]
are such that
\[
s(Z )=D+Z C(I_{\mathcal H}-ZA)^{-1}B
\]
and  they    form a coisometric matrix operator, that is,   $(5)$ holds. When $(5)$ holds we have for every $Z,   W  \in\Omega_{\mathbb
K^2}$:
\[
\frac{I_n-s(Z)s(W)^{*t}}{1-ZW^*}  =    C(I_{\mathcal
H}-ZA)^{-1}(I_{\mathcal H}-WA)^{-1}C^\sharp.
\]

Setting $Z = W$ we see that $(2)$ holds, which implies $(1)$
as one can see by considering the  idempotent    components.

\bigskip


\section{An example: bicomplex     Blaschke factors}
\label{sec:blaschke}

We conclude this work by presenting an example of     bicomplex   Schur multiplier, namely,    a   bicomplex  Blaschke
factor. Let $a=a_1+\j a_2\in\Omega_{\mathbb K^2}$. We set
\begin{equation}
\label{eq:bl}
b_a(Z) : =    \frac{Z-a}{1-Za^*}.
\end{equation}

Clearly,   the function $b_a$ is defined for all
$Z$ such that
 \begin{equation}
 (z_1-  \i z_2)(   \overline{ a_1}  +   \i   \, \overline{ a_2} \, )   \neq  1,\quad{\rm and}\quad
 (z_1+  \i  z_2)(   \overline{a_1}   -   \i \,   \overline{ a_2}  \,  ) \neq  1,
 \label{eq:not1}
 \end{equation}

 and in particular for all $Z \in\Omega_{\mathbb K^2}$.
 We will call $b_a$ the   bicomplex   Blaschke factor with zero $a$ when $aa^*\not= 1$.
 As in the complex case,
 \[
 b_a(Z)  =  -a+  Z  \frac{1-aa^*}{1-Z a^*},
 \] so that a realization of
 $b_a(Z)$ is given by
 \[
\begin{pmatrix}a^*&1\\1-aa^*&-a
 \end{pmatrix}  .
 \]
This realization is not coisometric, but this is readily fixed
as follows. Let $a=\eta_1\e  +   \eta_2\edag$, and
\[
u=   \sqrt{1-   |  \eta_1 |^2\, }  \,  \e  +    \sqrt{1-|\eta_2|^2  \, } \,  \edag  .
\]
Then $u\in\mathbb D^+$ and $u^2=1-aa^*$. The matrix
\[
\begin{pmatrix}A&B\\ C&D\end{pmatrix}=
\begin{pmatrix}a^*&u\\ u&-a
 \end{pmatrix}
\]
also defines a realization of $s$, and it is unitary (and in
particular co-isometric).

\begin{subPn}
It holds that
 \begin{equation}
 \label{eq:unit}
 b_a(Z)   b_a(Z)^*  =   \begin{cases}\,\,\lessdot 1,\quad {\rm for}\,\,   Z   \in\Omega_{\mathbb
 K^2},\\     \\
 \,\, 1,\quad\hspace{4mm} {\rm for}\,\,ZZ^*=1.
 \end{cases}
 \end{equation}
\end{subPn}

 {\sf Proof:} \\   Let $Z  \in \Omega_{\mathbb
 K^2}$. Then $1-Za^*$ is invertible and we have:
\[
\begin{split}
1-b_a(Z)   b_a(Z)^*  &    =    1-\frac{(  Z-a)(Z-a)^*}{(1-Za^*)(1-Z^*a)}\\
&=\frac{(1-Za^*)(1-Z^*a)-(Z-a)(Z-a)^*}{(1-Za^*)(1-Z^*a)}\\
&=\frac{(1-ZZ^*)(1-aa^*)}{(1-Za^*)(1-Z^*a)},
\end{split}
\]
which belongs to $\mathbb D^+$ as a product and quotient of
elements in $\mathbb D^+$.\\

Let now $ZZ^*=1$. The element $Z$ is in particular
 invertible and we can write for such $Z$
 \[
 b_a(Z)=\frac{Z-a}{ Z  (Z^*-a^*)},
 \]
 Thus,  for such $Z$,
 \[
 b_a(Z)b_a(Z)^*   =    \frac{(Z-a)(Z-a)^*}{ZZ^*(Z^*-a^*)(Z-a)}=1.
 \]
 \mbox{}\qed\mbox{}\\

 We note that both
 \begin{equation}
 \label{b1b2}
 b_{a, 1}  (\beta_1)  =     \frac{\beta_1-   \eta_1   }{1-\beta_1     \overline{  \eta_1}   }\quad
 {\rm and}
 \quad
 b_{a, 2} (\beta_2)   =     \frac{\beta_2-   \eta_2  }{1- \beta_2   \overline{  \eta_2}   }   \quad
 \end{equation}
 are classical Blaschke factors.\\

 We will study in a future publication interpolation problems in
 the present setting. Here we only mention the following result:

\begin{subPn}
Let      $a\in\Omega_{\mathbb K^2}$. Then:\\

\begin{enumerate}

\item[(1)]    For every $f\in\mathbf H^2(\Omega_{\mathbb K^2})$ it holds
that
 \[
   \langle b_af,b_af \rangle_{  \mathbf H^2(\Omega_{\mathbb K^2}) }   =     \langle f,f \rangle_{  \mathbf H^2(\Omega_{\mathbb K^2})  }  ,
 \]

 \noindent
 that is,  the operator of multiplication by $b_a$ is  $\D$--isometric.

\medskip

\item[(2)]    A function  $f$    belongs to $\mathbf H^2(\Omega_{\mathbb K^2})$
and vanishes at the point $a$ if and only if it can be written as
$f=b_ag$ with $g\in\mathbf H^2(\Omega_{\mathbb K^2})$. In this case,
\[
 \langle f,f  \rangle_{  \mathbf H^2(\Omega_{\mathbb K^2}) }   =    \langle g,g \rangle_{ \mathbf H^2(\Omega_{\mathbb K^2}) }    .
\]

\noindent
i.e.,  their  $\D$--norms  coincide.

\end{enumerate}
\end{subPn}

 {\sf Proof:}\\
$(1)$ We first note that for $Z  \in\Omega_{\mathbb K^2}$ we have:
 \[
 b_a(Z)=-a+(1-aa^*)\sum_{n=1}^\infty Z^n(a^*)^{n-1},
 \]
 where the convergence is in $\mathbf H^2(\Omega_{\mathbb K^2})$, and also
 pointwise, since we are in a reproducing kernel Hilbert space.
 Thus, for $N,M\in\mathbb N$ we have:
 \begin{equation}
 \label{eq:ba}
 b_a(Z)  Z^N=-aZ^N+(1-aa^*)\sum_{n=1}^\infty
 Z^{n+N}(a^*)^{n-1}.
 \end{equation}

Denote   $\mathbb N_0 := \mathbb N  \cup  \{ 0 \}$.    Thus, for $N\in\mathbb N_0$,
 \[
 \begin{split}
 \langle
 b_a(Z)  Z^N,b_a(Z)  Z^N \rangle_{  \mathbf H^2(\Omega_{\mathbb K^2}) }    &=aa^*+(1-aa^*)^2\sum_{n=1}^\infty(aa^*)^{n-1}\\
 &=aa^*+(1-aa^*)^2(1-aa^*)^{-1}\\
 &=1.
 \end{split}
 \]
 Let now $N,M\in\mathbb N_0$ with $N<M$, and let $M=N+H$. We have
 in view of \eqref{eq:ba}:
$$
\begin{array}{l}
  \langle b_a(Z)  Z^N,b_a(Z)  Z^M \rangle_{ \mathbf H^2(\Omega_{\mathbb K^2}) }  \;   =
\\  \\  
\qquad  \quad  \displaystyle = \;      -a^*(1-aa^*)(a^*)^{H-1}+(1-aa^*)\sum_{n=0}^\infty
 (a^*)^{H+n-1}a^{n-1}  \;   \\     \\
\qquad  \quad  \displaystyle     = \;  -a^*(1-aa^*)(a^*)^{H-1}+(1-aa^*)^2(a^*)^H(1-aa^*)^{-1}   \\    \\
 \qquad  \quad  \displaystyle    =   \;  0.
 \end{array}
 $$

 Hence, for any polynomial $p$ we have
 \[
 \langle b_ap,b_ap\rangle_{ \mathbf H^2(\Omega_{\mathbb K^2})  }    =\langle p,p\rangle_{ \mathbf H^2(\Omega_{\mathbb K^2}) }  ,
 \]
 and the result follows by continuity for every element in
 $\mathbf H^2(\Omega_{\mathbb K^2})$.\\

 $(2)$ One direction is clear. If $f=b_ag$ with
 $g\in\mathbf H^2(\Omega_{\mathbb K^2})$, then
 $f\in\mathbf H^2(\Omega_{\mathbb K^2})$ (this can  be  verified  component-wise) and vanishes at the point $a$. Conversely, assume that
 $f(a)=0$, and let $a=a_1+\j a_2$. Then $f_1(\eta_1)=f_2(\eta_2)=0$, where
\[
\eta_1=a_1-\i a_2\quad{\rm and}\quad \eta_2=a_1+\i a_2.
\]
By the classical theory we have (with $b_{a, 1}$ and $b_{a, 2}$ defined by
\eqref{b1b2}):
\[
f_j=b_{a, \, j}  g_j,\quad j=1,2,
\]
with $g_j\in\mathbf H^2(\mathbb K)$. The result follows by
regrouping the components.
\mbox{}\qed\mbox{}\\

As we have noted already in the Introduction  the classic,  complex  Schur  analysis  has  connections and applications to many interesting problems, thus we expect that  our brief study of some basic facts of  bicomplex  Schur analysis will become a  starting  point  for  deeper  developments of it together with many applications.

\bigskip

\bigskip

\bibliographystyle{amsalpha}

\end{document}